\documentclass{article}

\usepackage[letterpaper,margin=1in]{geometry}
\usepackage{graphicx} 
\usepackage{amsmath}
\usepackage{amsthm}
\usepackage{amssymb}
\usepackage{booktabs}
\usepackage{longtable,pdflscape}
\usepackage{multirow}
\usepackage{mathrsfs}
\usepackage{makecell}
\usepackage{array}
\usepackage{caption}
\usepackage[authoryear,round,sort]{natbib}
\ifdefined\cacesubmission
\usepackage{xr-hyper}
\fi
\usepackage{hyperref}
\ifdefined\cacesubmission
\fi
\usepackage{xurl}
\usepackage{authblk}

\usepackage{placeins}
\usepackage{float}

\usepackage{relsize}              
\newcommand{\bigveebar}{%
  \mathop{\mathlarger{\mathlarger{\veebar}}}\displaylimits%
}

\usepackage{tikz}
\usetikzlibrary{arrows.meta}
\usepackage{threeparttable}
\definecolor{lime}{HTML}{A6CE39}
\newcommand{\figbox}[3]{\parbox[t][#1][t]{#2}{\raggedright\hyphenpenalty=10000\exhyphenpenalty=10000\setlength{\baselineskip}{9pt}#3}}
\definecolor{oiblue}{HTML}{0072B2}
\definecolor{oisky}{HTML}{56B4E9}
\definecolor{oiorange}{HTML}{E69F00}
\definecolor{oigreen}{HTML}{009E73}
\definecolor{oivermilion}{HTML}{D55E00}
\definecolor{oipurple}{HTML}{CC79A7}
\DeclareRobustCommand{\orcidicon}{
	\begin{tikzpicture}
	\draw[lime, fill=lime] (0,0)
	circle [radius=0.16]
	node[white] {{\fontfamily{qag}\selectfont \tiny ID}};
	\draw[white, fill=white] (-0.0625,0.095)
	circle [radius=0.007];
	\end{tikzpicture}
	\hspace{-2mm}
}
\foreach \x in {A, ..., Z}{\expandafter\xdef\csname orcid\x\endcsname{\noexpand\href{https://orcid.org/\csname orcidauthor\x\endcsname}
			{\noexpand\orcidicon}}
}

\newcommand{\david}[1]{\textcolor{red}{D: #1}}
\newcommand{\albert}[1]{\textcolor{blue}{A: #1}}
\newcommand{\claude}[1]{\textcolor{teal}{C: #1}}

\newcommand{\papertitle}{Generalized Disjunctive Programming with Derivative-Free Configuration Search}
\newcommand{\paperfirstauthor}{Albert Joon Lee}
\newcommand{\papersecondauthor}{Woosik Kim}
\newcommand{\papercorrespondingauthor}{David E. Bernal Neira}
\newcommand{\papercorrespondingemail}{dbernaln@purdue.edu}
\newcommand{\paperauthors}{\paperfirstauthor, \papersecondauthor, and \papercorrespondingauthor}
\title{\papertitle}
\author[1]{\paperfirstauthor\orcidB}
\author[2]{\papersecondauthor\orcidC}
\author[1,*]{\papercorrespondingauthor\orcidA}

\affil[1]{Davidson School of Chemical Engineering, Purdue University, West Lafayette, IN 47907, USA}
\affil[2]{Elmore Family School of Electrical and Computer Engineering, Purdue University, West Lafayette, IN 47907, USA}

\date{*Email: \papercorrespondingemail}

\begin{document}

\maketitle

\begin{abstract}
Generalized disjunctive programming (GDP) describes process design alternatives with Boolean variables, disjunctions, and logical propositions, and fixing one configuration leaves a subproblem in the continuous variables.
Conventional GDP methods require explicit algebraic formulations that simulators may not expose.
Even when such formulations are available, nonconvex models can be intractable.
The external variables define an integer lattice that derivative-free optimization (DFO) searches using subproblem values.
We compare surrogate DFO methods and the direct-search method \texttt{NOMAD} with discrete steepest descent (\texttt{D-SDA}), logic-based solvers, mixed-integer nonlinear programming (MINLP) solvers, and enumeration.
The benchmarks include nine process models, 21 test functions, a pharmaceutical flowsheet simulated in \texttt{PharmaPy}, and a catalytic distillation column with 100,947 configurations.
On the largest reactor network, the surrogate methods \texttt{MISO} and \texttt{RBFOpt} reach the reference in nine or ten of ten runs with about 30 subproblem solves, whereas logic-based methods solve nearly one subproblem for each admissible configuration and most \texttt{D-SDA} descents end at one local minimum.
In the flowsheet example, where GDP solvers do not apply, \texttt{RBFOpt} and \texttt{MISO} find the best of 36 simulated configurations in every run with 20 evaluations.
For the catalytic column, \texttt{D-SDA} and \texttt{MISO} reach the best known design with 100 evaluations, whereas MINLP solvers find no design within 15\% of it and global solvers return contradictory certificates.
This search pays off when the configuration space is too large to enumerate, has deceptive local minima, or has MINLP reformulations that fail, and enumeration or a direct reformulation is simpler for small algebraic models.
\end{abstract}


\section{Introduction}
In chemical process design, evaluating a single candidate configuration can require costly optimization or simulation, making the choice of which configuration to explore next a central computational challenge.
Process design involves discrete decisions, such as equipment selection, the number of processing stages, and operating modes, along with continuous variables such as flow rates, temperatures, and equipment sizes.
Generalized disjunctive programming (GDP) provides a framework to represent these decisions through Boolean variables, logical propositions, and conditional constraints~\citep{raman1994modelling}.
In logic-based GDP solution methods, fixing the Boolean variables associated with a selected configuration determines which constraints remain active and yields a reduced optimization subproblem.
When all discrete decisions are fixed, this subproblem can be formulated as a nonlinear programming (NLP) problem.
Exploring multiple configurations may therefore require repeatedly solving NLP subproblems, which can be computationally expensive.
Consequently, identifying high-quality configurations with fewer evaluations of NLP subproblems is an important objective in solving GDP problems.

Conventional GDP solution methods generally rely on mixed-integer nonlinear programming (MINLP) reformulations or logic-based decomposition.
The Big-\(M\) and Hull reformulations transform disjunctive constraints into algebraic formulations~\citep{raman1994modelling,grossmann2003hull}, while logic-based approaches, such as logic-based outer approximation (LOA) and the relaxation with integer cuts (RIC) of \texttt{GDPopt}, select discrete configurations and solve the corresponding reduced subproblems~\citep{turkay1996logic,bergamini2005logic,pyomoGDPopt}.
In logic-based approaches, the strategy used to select discrete configurations can significantly affect the number of subproblem evaluations required.
Both families need the algebraic equations of every disjunct.
When a unit is described by a process simulator, these equations are not available, and when the model is nonconvex, a reformulation keeps the equations of units that a configuration does not contain, whose terms can be poorly scaled or undefined at the values that absent units take.

The Discrete-Steepest Descent Algorithm (\texttt{D-SDA}) was previously applied to process superstructure optimization by \citet{linan2020optimal}.
Building on this approach, our previous work introduced the Logic-Based Discrete-Steepest Descent Algorithm (\texttt{LD-SDA}), which represents ordered Boolean decisions using external integer variables~\citep{ovalle2025logic}.
With this external-variable reformulation, the GDP becomes an optimization over an integer lattice whose objective is the value of the subproblem of each configuration.
The GDP model still describes the alternatives, their logic, and the subproblem of each configuration, but the method that searches the lattice needs only the value of each configuration.
\texttt{LD-SDA} explores the integer space through neighborhood search and discrete line search, evaluating selected configurations by solving their corresponding reduced optimization subproblems.
However, as a local search method, \texttt{LD-SDA} may converge to different locally optimal solutions depending on the starting point and the selected neighborhood~\citep{peng2024dynamic}.

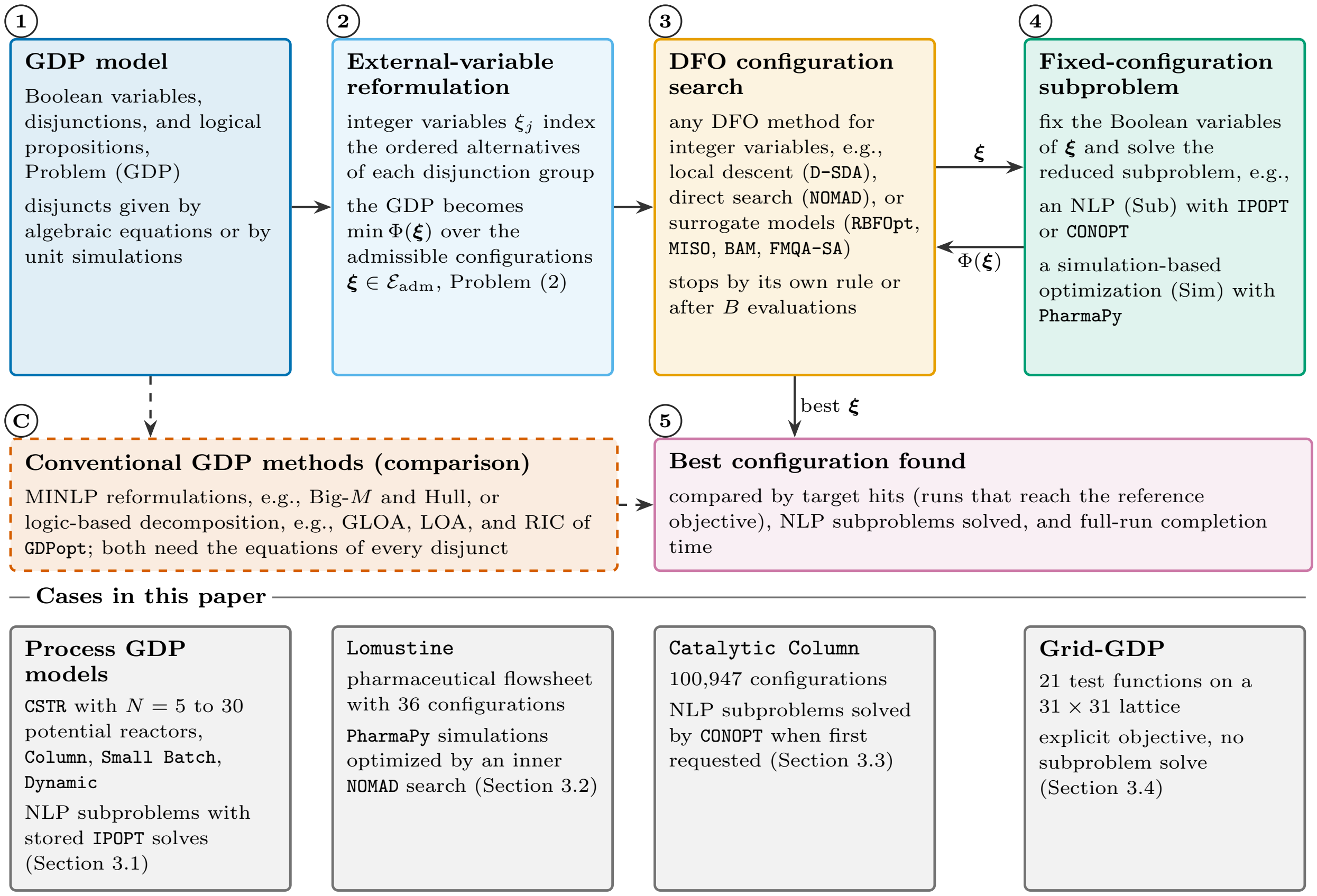
\begin{figure}[tbp]
\centering
\begin{tikzpicture}[font=\scriptsize,
  box/.style={draw=#1, fill=#1!12, line width=0.9pt, rounded corners=2pt, inner sep=5pt, anchor=north west},
  case/.style={box=gray, draw=black!55, fill=black!5},
  badge/.style={circle, draw=black!85, fill=white, line width=0.6pt, inner sep=1.2pt, minimum size=4mm, font=\scriptsize\bfseries},
  arr/.style={-{Stealth[length=2.2mm]}, line width=0.8pt, draw=black!80},
  lab/.style={font=\scriptsize, inner sep=1.5pt}]
\node[box=oiblue] (gdp) at (0,0) {\figbox{3.85cm}{3.15cm}{{\footnotesize\bfseries GDP model}\\[3pt]
  Boolean variables, disjunctions, and logical propositions, Problem~\eqref{eqn:gdp.general}\\[3pt]
  disjuncts given by algebraic equations or by unit simulations}};
\node[box=oisky] (ext) at (4.03,0) {\figbox{3.85cm}{3.15cm}{{\footnotesize\bfseries External-variable\\reformulation}\\[3pt]
  integer variables \(\xi_j\) index the ordered alternatives of each disjunction group\\[3pt]
  the GDP becomes \(\min\Phi(\boldsymbol\xi)\) over the admissible configurations \(\boldsymbol\xi\in\mathcal E_{\mathrm{adm}}\), Problem~\eqref{eqn:outer.value}}};
\node[box=oiorange] (dfo) at (8.06,0) {\figbox{3.85cm}{3.15cm}{{\footnotesize\bfseries DFO configuration\\search}\\[3pt]
  any DFO method for integer variables, e.g., local descent (\texttt{D-SDA}), direct search (\texttt{NOMAD}), or surrogate models (\texttt{RBFOpt}, \texttt{MISO}, \texttt{BAM}, \texttt{FMQA-SA})\\[3pt]
  stops by its own rule or after \(B\) evaluations}};
\node[box=oigreen] (eval) at (12.69,0) {\figbox{3.85cm}{3.15cm}{{\footnotesize\bfseries Fixed-configuration\\subproblem}\\[3pt]
  fix the Boolean variables of \(\boldsymbol\xi\) and solve the reduced subproblem, e.g.,\\[3pt]
  an NLP~\eqref{eqn:Sub} with \texttt{IPOPT} or \texttt{CONOPT}\\[3pt]
  a simulation-based optimization~\eqref{eqn:sim-sub} with \texttt{PharmaPy}}};
\node[box=oivermilion, dashed] (conv) at (0,-5.0) {\figbox{1.3cm}{7.24cm}{{\footnotesize\bfseries Conventional GDP methods (comparison)}\\[3pt]
  MINLP reformulations, e.g., Big-\(M\) and Hull, or logic-based decomposition, e.g., GLOA, LOA, and RIC of \texttt{GDPopt}; both need the equations of every disjunct}};
\node[box=oipurple] (best) at (8.06,-5.0) {\figbox{1.3cm}{7.87cm}{{\footnotesize\bfseries Best configuration found}\\[3pt]
  compared by target hits (runs that reach the reference objective), NLP subproblems solved, and full-run completion time}};
\foreach \n/\t in {gdp/1, ext/2, dfo/3, eval/4, best/5}
  \node[badge, anchor=south west] at ([yshift=0.6mm]\n.north west) {\t};
\node[badge, anchor=south west] at ([yshift=0.6mm]conv.north west) {C};
\draw[arr] (gdp.east) -- (ext.west);
\draw[arr] (ext.east) -- (dfo.west);
\draw[arr] ([yshift=5mm]dfo.east) -- node[lab, above] {\(\boldsymbol\xi\)} ([yshift=5mm]eval.west);
\draw[arr] ([yshift=-5mm]eval.west) -- node[lab, below] {\(\Phi(\boldsymbol\xi)\)} ([yshift=-5mm]dfo.east);
\draw[arr] (dfo.south) -- node[lab, right] {best \(\boldsymbol\xi\)} (dfo.south |- best.north);
\draw[arr, dashed] (gdp.south) -- (gdp.south |- conv.north);
\draw[arr, dashed] (conv.east) -- (best.west);
\draw[black!55, line width=0.6pt] (0,-7.0) -- (16.22,-7.0);
\node[anchor=west, font=\footnotesize\bfseries, fill=white, inner sep=2pt] at (0.25,-7.0) {Cases in this paper};
\node[case] (c1) at (0,-7.35) {\figbox{2.95cm}{3.15cm}{{\footnotesize\bfseries Process GDP models}\\[2pt]
  \texttt{CSTR} with \(N=5\) to 30 potential reactors, \texttt{Column}, \texttt{Small Batch}, \texttt{Dynamic}\\[2pt]
  NLP subproblems with stored \texttt{IPOPT} solves (Section~\ref{sec:process-cases})}};
\node[case] (c2) at (4.03,-7.35) {\figbox{2.95cm}{3.15cm}{{\footnotesize\bfseries \texttt{Lomustine}}\\[2pt]
  pharmaceutical flowsheet with 36 configurations\\[2pt]
  \texttt{PharmaPy} simulations optimized by an inner \texttt{NOMAD} search (Section~\ref{sec:lomustine-case})}};
\node[case] (c3) at (8.06,-7.35) {\figbox{2.95cm}{3.15cm}{{\footnotesize\bfseries \texttt{Catalytic Column}}\\[2pt]
  100,947 configurations\\[2pt]
  NLP subproblems solved by \texttt{CONOPT} when first requested (Section~\ref{sec:catdist-case})}};
\node[case] (c4) at (12.69,-7.35) {\figbox{2.95cm}{3.15cm}{{\footnotesize\bfseries Grid-GDP}\\[2pt]
  21 test functions on a \(31\times31\) lattice\\[2pt]
  explicit objective, no subproblem solve (Section~\ref{sec:grid-gdp-formulation})}};
\end{tikzpicture}
\caption{Methodology of this study.
Top, the general method: a GDP model (1) is rewritten with external integer variables \(\boldsymbol\xi\) that index the alternatives of its disjunctions (2), so that solving it becomes a search over the admissible configurations whose objective is the subproblem value \(\Phi(\boldsymbol\xi)\).
A DFO method (3) proposes a configuration, which is evaluated by fixing its Boolean variables and solving the reduced subproblem (4), and the returned value updates the search until the method stops or the evaluation budget \(B\) is spent; the best configuration found (5) is the result.
The dashed path (C) shows the conventional GDP methods, MINLP reformulations or logic-based methods such as global logic-based outer approximation (GLOA), logic-based outer approximation (LOA), and relaxation with integer cuts (RIC), which solve the same models when the equations of every disjunct are available.
Bottom: the four cases of this paper and the subproblem evaluation each one uses.}
\label{fig:methodology}
\end{figure}

Because the lattice search needs only configuration values, any derivative-free optimization (DFO) method that handles integer variables can perform it, and the GDP model can be solved even when its subproblems are simulations or its MINLP reformulations fail.
We pair GDP models in their external integer representation with DFO methods that search the configurations.
Figure~\ref{fig:methodology} summarizes this approach and the conventional GDP methods with which we compare it.
We include the mesh adaptive direct-search solver \texttt{NOMAD}~\citep{audet2022nomad} and four surrogate-based methods, which use accumulated objective observations to guide new evaluations~\citep{larson2019derivative}.
\texttt{RBFOpt} builds surrogate models from radial basis functions and balances predicted objective values with exploration~\citep{costa2018rbfopt}, and \texttt{MISO} combines radial-basis-function surrogate models with sampling strategies designed for integer variables~\citep{muller2016miso}.
\texttt{NOMAD} and \texttt{MISO} performed best on mixed-integer problems in the comparison of \citet{ploskas2022review}.
The other two are Branch-and-Model (\texttt{BAM})~\citep{ma2023bam} and a variant of the factorization machine with quantum annealing (\texttt{FMQA}) framework~\citep{kitai2020designing}, in which classical simulated annealing replaces quantum annealing; we denote this variant \texttt{FMQA-SA}.
\texttt{FMQA} represents its surrogate as a quadratic unconstrained binary optimization problem, the input format of quantum annealers, which makes it a candidate for a possibly quantum-assisted configuration search.
We ask three questions.
First, do these methods find high-quality configurations with fewer NLP subproblem solves than \texttt{D-SDA} and conventional GDP methods, and when does this advantage translate into shorter computing time?
Second, does the same search solve GDP models whose subproblems are evaluated by a simulator?
Third, does it solve GDP models whose MINLP reformulations defeat conventional solvers?

Candidate selection has its own cost, and evaluation costs vary between configurations~\citep{more2009benchmarking}.
We therefore report three measures for every method: how often it reaches the reference objective, how many NLP subproblems it solves, and its estimated completion time.
The complete enumeration of each configuration domain and uniform sampling serve as references, and all methods share initial configurations, evaluation budgets, and the responses of each configuration.
The contributions are as follows:
\begin{itemize}
\item A use of the external-variable reformulation as the interface between a GDP model and a derivative-free configuration search, extended to subproblems defined by a process simulator and demonstrated on a lomustine manufacturing flowsheet simulated in \texttt{PharmaPy} and written as a GDP.
\item A comparison of surrogate-based and direct-search configuration methods with local descent, logic-based and direct GDP methods, and full enumeration on four process GDP families and 21 Grid test functions, measured by target attainment, NLP subproblems solved, and completion time.
\item Evidence from execution records that the logic-based methods RIC and LOA solve one subproblem for each admissible configuration of most of these models, whereas surrogate search reaches near-optimal configurations with a small fraction of those subproblems.
\item A benchmark with live NLP evaluations on a catalytic distillation column whose 100,947 configurations are too many to enumerate, together with valid Big-\(M\) and Hull reformulations of its GDP model, on which conventional MINLP solvers fail, and global solvers return contradictory certificates.
\item An analysis of landscape features that determines when the surrogate search outperforms local descent, with recommendations for choosing between enumeration, direct reformulation, local descent, and surrogate search.
\end{itemize}

\section{Background}

\subsection{Generalized Disjunctive Programming and Fixed-Configuration Evaluation}
\label{sec:gdp-background}

GDP represents operational and design alternatives using Boolean variables and associated algebraic constraints~\citep{raman1994modelling}.
The original continuous variables have domain \(\mathcal X=[\underline{\mathbf x},\overline{\mathbf x}]\subseteq\mathbb R^{n_x}\).
The original integer variables have domain \(\mathcal Z=[\underline{\mathbf z},\overline{\mathbf z}]\cap\mathbb Z^{n_z}\).
Let \(K\) index the disjunctions and \(D_k\) index the disjuncts in disjunction \(k\in K\).
Each disjunct \(i\in D_k\) has a Boolean variable \(Y_{ik}\).
The GDP formulation is given by
\begin{equation}\label{eqn:gdp.general}
\begin{aligned}
\min_{\mathbf x,\mathbf z,\mathbf Y}\quad & f(\mathbf x,\mathbf z)\\
\text{s.t.}\quad &\mathbf g(\mathbf x,\mathbf z)\leq\mathbf 0,\\
&\bigveebar_{i\in D_k}\left[
\begin{gathered}
Y_{ik}\\
\mathbf h_{ik}(\mathbf x,\mathbf z)\leq\mathbf 0
\end{gathered}\right],\qquad k\in K,\\
&\Omega(\mathbf Y)=\mathrm{True},\\
&\mathbf x\in\mathcal X,\quad\mathbf z\in\mathcal Z,\quad
\mathbf Y\in\{\mathrm{False},\mathrm{True}\}^{n_y}.
\end{aligned}\tag{GDP}
\end{equation}
Here, \(f\) is the objective function, \(\mathbf g\) represents the global constraints, and \(\mathbf h_{ik}\) represents the constraints associated with each disjunct.
The operator \(\veebar\) requires that exactly one disjunct be selected within each disjunction.
Logical propositions \(\Omega(\mathbf Y)\) impose additional relationships between Boolean variables.

Following the external-variable reformulation of \citet{ovalle2025logic}, each ordered group of \(m_j\) Boolean alternatives is represented by an external integer variable \(\xi_j\in\{1,\ldots,m_j\}\), which identifies the selected alternative.
Ordered alternatives, such as the position of a feed along a column, are one of the conditions that \citet{ovalle2025logic} use because they make neighboring values of \(\xi_j\) physically similar configurations.
The representation itself does not require this order: unordered alternatives, such as a choice among equipment types, can be indexed in any fixed order, as in Section~\ref{sec:lomustine-case}, although the order then shapes the neighborhoods and surrogate models of the search (Section~\ref{sec:lomustine-results}).
Let \(n_E\) denote the number of external integer variables, and \(\mathcal Y_{\mathrm{logic}}\) the set of Boolean assignments that satisfy the conditions \texttt{exactly-one}, and \(\Omega(\mathbf Y)\).
The external-variable lattice and its logically admissible subset are
\begin{equation}\label{eqn:external.lattice}
\begin{aligned}
\mathcal L_E&=\prod_{j=1}^{n_E}\{1,\ldots,m_j\},\\
\mathcal E_{\mathrm{adm}}&=\{\boldsymbol\xi\in\mathcal L_E:
\mathbf Y(\boldsymbol\xi)\in\mathcal Y_{\mathrm{logic}}\}.
\end{aligned}
\end{equation}
For each logically admissible configuration \(\boldsymbol\xi\), the corresponding Boolean variables are fixed, yielding the following reduced optimization subproblem:

\begin{equation}\label{eqn:Sub}
\begin{aligned}
\Phi(\boldsymbol\xi)=\inf_{\mathbf x,\mathbf z}\quad &f(\mathbf x,\mathbf z)\\
\text{s.t.}\quad &\mathbf g(\mathbf x,\mathbf z)\leq\mathbf 0,\\
&\mathbf h_{ik}(\mathbf x,\mathbf z)\leq\mathbf 0
\quad\text{for all }(i,k)\text{ with }Y_{ik}(\boldsymbol\xi)=\mathrm{True},\\
&\mathbf x\in\mathcal X,\quad\mathbf z\in\mathcal Z.
\end{aligned}\tag{Sub}
\end{equation}
Here, \(\Phi(\boldsymbol\xi)\) represents the optimal objective value of the reduced subproblem.
The reduced subproblem is an NLP when all remaining decision variables are continuous.
If the reduced subproblem is infeasible, we define \(\Phi(\boldsymbol\xi)=+\infty\).

When the external-variable mapping represents all logically admissible Boolean configurations, the original GDP can be expressed as the following upper-level optimization problem~\citep{ovalle2025logic}:
\begin{equation}\label{eqn:outer.value}
\min_{\boldsymbol\xi\in\mathcal E_{\mathrm{adm}}}\Phi(\boldsymbol\xi).
\end{equation}
A logically admissible configuration may still yield an infeasible reduced subproblem.
Problem~\eqref{eqn:outer.value} separates the two tasks of solving a GDP: the model defines the admissible configurations and the subproblem of each, and the method that solves Problem~\eqref{eqn:outer.value} decides which configurations to evaluate.

Problem~\eqref{eqn:outer.value} uses only the value \(\Phi(\boldsymbol\xi)\) of each configuration, so it also applies when a configuration is evaluated by a process simulator.
Fixing \(\boldsymbol\xi\) then selects the unit models of a flowsheet, and the simulator \(\mathcal S_{\boldsymbol\xi}\) computes its states \(\mathbf s\), such as stream compositions and crystal size distributions, from operating inputs \(\mathbf u\in\mathcal U_{\boldsymbol\xi}\).
The reduced subproblem becomes
\begin{equation}\label{eqn:sim-sub}
\Phi(\boldsymbol\xi)=\min_{\mathbf u\in\mathcal U_{\boldsymbol\xi}}\bigl\{f(\mathbf u,\mathbf s):\ \mathbf s=\mathcal S_{\boldsymbol\xi}(\mathbf u),\ \mathbf g(\mathbf u,\mathbf s)\leq\mathbf 0\bigr\},\tag{Sim}
\end{equation}
whose objective and constraints are known only through simulation outputs.
In \texttt{PharmaPy}, for example, each unit is described by differential-algebraic equations (DAEs), and the simulator integrates them unit by unit in the order of the flowsheet.
A DAE model can instead be discretized by orthogonal collocation on finite elements into algebraic equations, which turns its subproblem into an NLP; the \texttt{Dynamic} benchmark of Section~\ref{sec:process-cases} is solved this way with \texttt{pyomo.dae}~\citep{nicholson2018pyomo}.
A derivative-free method in \(\mathbf u\) solves Problem~\eqref{eqn:sim-sub}, so \(\Phi(\boldsymbol\xi)\) is approximated by the best value that this inner search finds within its simulation budget.
Logic-based and direct GDP methods require the equations of every disjunct and do not apply to Problem~\eqref{eqn:sim-sub}.
The configuration-search methods compared here use it through the same evaluation interface as an NLP subproblem, which Section~\ref{sec:lomustine-case} illustrates with a pharmaceutical flowsheet simulated in \texttt{PharmaPy}~\citep{casasorozco2021pharmapy}.

\subsection{Conventional GDP Solution Approaches}

Big-\(M\) and Hull reformulations transform GDP models into MINLP formulations by associating Boolean variables \(Y_{ik}\) with binary variables \(y_{ik}\) and translating logical propositions into algebraic constraints~\citep{raman1994modelling,grossmann2003hull}.
In the Big-\(M\) reformulation, the disjunctive constraints are expressed as \(\mathbf h_{ik}(\mathbf x,\mathbf z)\leq\mathbf M_{ik}(1-y_{ik})\), where the inequality is applied componentwise.
The Big-\(M\) values must be sufficiently large to render the corresponding constraints redundant when \(y_{ik}=0\), but excessively large values can weaken the continuous relaxation and cause numerical difficulties.

Hull reformulations introduce disaggregated variables, linking constraints, and bounds scaled by binary indicators~\citep{grossmann2003hull}.
Disjunctive constraints are imposed on the disaggregated variables when \(y_{ik}=1\), while \(y_{ik}=0\) forces these variables to zero, making the corresponding disjunct inactive.
For convex disjunctive sets, this formulation can provide a stronger continuous relaxation than Big-\(M\), although it introduces additional variables and constraints.
For nonconvex GDP problems, the Hull reformulation does not guarantee a convex reformulated model or global optimality~\citep{chen2022pyomogdp}.

Logic-based methods address GDP problems directly by exploiting their logical and disjunctive structure to select discrete configurations and construct the corresponding optimization subproblems.
LOA alternates between a mixed-integer master problem and fixed-configuration NLP subproblems~\citep{turkay1996logic}.
Under suitable convexity and regularity assumptions, the master problem uses linearization cuts to obtain lower bounds, while feasible solutions to the NLP subproblem provide upper bounds.
For nonconvex GDP problems, these linearizations may not provide valid global lower bounds, motivating global logic-based outer approximation (GLOA), which requires globally valid relaxations and subproblem bounds.
\citet{bergamini2005logic} developed a GLOA approach for process networks using piecewise under- and overestimators.

\texttt{GDPopt} provides implementations of logic-based methods within \texttt{Pyomo.GDP}~\citep{chen2022pyomogdp,pyomoGDPopt}.
Its RIC method uses a mixed-integer linear programming (MILP) relaxation to select configurations and adds no-good cuts to exclude evaluated discrete assignments, without generating the successive outer-approximation cuts used by LOA~\citep{pyomoGDPopt}.

\texttt{LD-SDA} is an existing alternative that searches ordered external configurations by improving neighborhood moves and discrete line searches~\citep{ovalle2025logic}.
It selects the best improving neighbor, continues along the resulting direction while improvement persists, and resumes neighborhood exploration.
The present work uses this discrete-search structure as a reference and compares other configuration-selection strategies and their computational costs.

\subsection{Derivative-Free Optimization Solvers}

DFO uses function observations without requiring derivatives of the objective with respect to the search variables~\citep{larson2019derivative}.
Black-box configuration search selects configurations using their evaluation responses without accessing the subproblem equations.
For GDP configuration search, the outer method can be derivative-free, even if its evaluation routine uses a derivative-based nonlinear solver.
\texttt{D-SDA} also uses this derivative-free interface, although it does not construct a surrogate model.
The upper-level search operates over external integer variables that represent the logical decisions in the GDP. 
Some logical restrictions are embedded in the external-variable reformulation, while additional constraints can be imposed on the configuration space.
One approach to DFO uses surrogate models to guide the search.
Surrogate methods seek to reduce expensive evaluations by fitting models to observed data, with computational cost arising from both objective evaluations and model-based candidate selection.

Several DFO methods handle integer or categorical variables directly, and \citet{ploskas2022review} compare many of them on mixed-integer test problems.
Mesh adaptive direct search, implemented in \texttt{NOMAD}, extends pattern search to mixed variables through user-defined discrete neighborhoods~\citep{abramson2009mads,audet2022nomad}, and \texttt{MISO} combines radial basis function surrogates with sampling strategies for integer variables~\citep{muller2016miso}.
Bayesian optimization has also been adapted to combinatorial spaces, for example by fitting sparse second-order models of binary variables in \texttt{BOCS}~\citep{baptista2018bocs}, which are closely related to the factorization-machine surrogate used here, or by placing kernels on graphs of discrete configurations in \texttt{COMBO}~\citep{oh2019combo}.
The direct-search method \texttt{NOMAD} polls trial configurations on an adaptive mesh and keeps no global surrogate.
We also compare four surrogate-based methods that represent different modeling choices: \texttt{RBFOpt} and \texttt{MISO} interpolate all observations with global radial basis function models but select candidates differently, \texttt{BAM} partitions the domain and fits local algebraic models, and \texttt{FMQA-SA} fits a quadratic model over a binary encoding.

\subsubsection{\texttt{NOMAD}}
\texttt{NOMAD} implements mesh adaptive direct search (MADS)~\citep{audet2022nomad}.
Each poll evaluates trial points on a mesh around the incumbent, and the mesh is refined after unsuccessful polls and coarsened after successful ones; for integer variables, it cannot become finer than the integer spacing~\citep{abramson2009mads}.
By default, \texttt{NOMAD}~4 complements the poll with search steps based on local quadratic models, but it does not retain a global surrogate of the objective.

\subsubsection{\texttt{RBFOpt}}
\texttt{RBFOpt} is an open-source surrogate optimizer for bounded black-box problems, with support for continuous, integer, and categorical variables~\citep{costa2018rbfopt,nannicini2021implementation}.
Given \(N\) observed objective values in external integer configurations \(\boldsymbol\xi^{(1)},\ldots,\boldsymbol\xi^{(N)}\), \texttt{RBFOpt} constructs a radial basis function (RBF) surrogate model of the form
\begin{equation}\label{eqn:rbf.surrogate}
\widehat\Phi_N(\boldsymbol\xi)
=\sum_{r=1}^{N}\lambda_r\varphi\!\left(\|\boldsymbol\xi-\boldsymbol\xi^{(r)}\|_2\right)
+p(\boldsymbol\xi),
\end{equation}
where \(\widehat\Phi_N(\boldsymbol\xi)\) approximates the objective value for a candidate configuration \(\boldsymbol\xi\), \(\lambda_r\) are interpolation coefficients determined from the observed values, \(\varphi\) is the RBF kernel evaluated at the Euclidean distance between configurations, and \(p(\boldsymbol\xi)\) is a polynomial tail depending on the selected kernel.

Candidate selection balances low predicted values against exploration away from evaluated points.
The metric stochastic response surface method (MSRSM) combines scaled surrogate predictions and distances to existing samples, changing their relative weights during the search~\citep{regis2007stochastic}.
New objective evaluations update the model.

\subsubsection{\texttt{MISO}}
\texttt{MISO} uses the radial basis function surrogate of Equation~\eqref{eqn:rbf.surrogate} but selects candidates with strategies designed for integer and mixed-integer variables~\citep{muller2016miso}.
Its coordinate perturbation and target value (CPTV) strategy first perturbs randomly chosen coordinates of the best configuration found and scores the candidates by surrogate value and distance to the evaluated points, as MSRSM does.
When this stage stops improving, it switches to a target value step, which proposes the configuration at which the surrogate could attain a prescribed target value with the least added bumpiness.

\subsubsection{\texttt{BAM}}
\texttt{BAM} combines domain partitioning with local surrogate construction~\citep{ma2023bam}.
The evaluated points define boxes, and a potentially optimal selection procedure identifies regions for further modeling and search.
Within selected boxes, Automated Learning of Algebraic Models (\texttt{ALAMO}) constructs sparse algebraic models from nearby observations~\citep{cozad2014learning}.
\texttt{BARON} optimizes these local models to propose candidates, while box-size exploration maintains broader coverage~\citep{ma2023bam,bamManual,sahinidis1996baron}.
Within \texttt{BAM}, \texttt{BARON} optimizes these learned models and never receives the GDP itself.
Its support for integer restrictions allows for the configuration domains considered here~\citep{bamManual}.

\subsubsection{Factorization-Machine-Based Optimization}

A second-order factorization machine (FM) models individual feature effects and pairwise interactions through shared latent parameters~\citep{rendle2010factorization}.
For a binary encoding \(\mathbf b\in\{0,1\}^{m}\) of an external integer configuration \(\boldsymbol\xi\), the FM fit surrogate is
\begin{equation}\label{eqn:fm.surrogate}
\widehat\Phi_{\mathrm{FM}}(\mathbf b)
=\beta_0+\sum_{i=1}^{m}\beta_i b_i
+\sum_{i=1}^{m-1}\sum_{j=i+1}^{m}
\langle\mathbf a_i,\mathbf a_j\rangle b_i b_j,
\end{equation}
where \(\widehat\Phi_{\mathrm{FM}}(\mathbf b)\) denotes the predicted training response, \(m\) is the number of binary features, \(\beta_0\) is the intercept, \(\beta_i\) is the linear coefficient of \(b_i\), and \(\mathbf a_i\in\mathbb R^r\) is its latent vector of dimension \(r\).
The inner product \(\langle\mathbf a_i,\mathbf a_j\rangle\) determines the pairwise interaction coefficient.
The model parameters are learned from the evaluated configurations and their corresponding training responses.

For binary inputs, \(b_i^2=b_i\), the fitted model is a binary quadratic model.
Excluding the constant term \(\beta_0\) gives the following quadratic unconstrained binary optimization (QUBO) problem:
\begin{equation}\label{eqn:fm.qubo}
\min_{\mathbf b\in\{0,1\}^{m}}\mathbf b^{\mathsf T}\mathbf Q\mathbf b,
\qquad
Q_{ij}=\begin{cases}
\beta_i,&i=j,\\
\langle\mathbf a_i,\mathbf a_j\rangle,&i<j,\\
0,&i>j.
\end{cases}
\end{equation}
The QUBO models the learned surrogate, and a binary string may decode to a configuration outside the domain or one that violates the logical constraints.

The \texttt{FMQA} framework combines FM surrogate fitting with annealing-based candidate generation~\citep{kitai2020designing}.
Based on observed evaluations, the method fits the FM, samples binary candidates from the resulting quadratic model, and evaluates a selected decoded configuration through the objective-evaluation routine before updating the training data.

\texttt{FMQA-SA} uses classical simulated annealing to sample binary candidates from the learned QUBO.
Temperature-dependent acceptance of uphill moves supports exploration, while gradual cooling favors lower-energy solutions~\citep{kirkpatrick1983optimization}.
In this study, no quantum hardware was used.

\section{Experimental Methods}
\label{sec:methodology}\label{sec:cases}

\subsection{Process GDP benchmark problems}
\label{sec:process-cases}
The process benchmarks include a network of continuous stirred-tank reactors (CSTRs)~\citep{linan2020optimal}, a 17-stage distillation column~\citep{ghouse2018comparative}, and a small batch plant~\citep{kocis1988global}.
We use their GDP formulations reported by \citet{ovalle2025logic}, and the nine-stage nonlinear mode-switching model of \citet{peng2024dynamic}.
Their objectives minimize the total reactor volume (\texttt{CSTR}), condenser and reboiler duties, plus the active tray cost (\texttt{Column}), the investment cost (\texttt{Small Batch}), and the negative integral of the squared state (\texttt{Dynamic}).
The external variables represent the reactor/feed cutoff and recycle location (\texttt{CSTR}), reflux and boilup positions (\texttt{Column}), parallel unit counts (\texttt{Small Batch}), and final stages of modes 1 and 2 (\texttt{Dynamic}).
For these process benchmarks, fixing the external choices leaves only continuous variables, so the resulting subproblems are NLPs.
The \texttt{CSTR} superstructure contains \(N\) potential reactor positions, where \(N\in\{5,10,15,20,25,30\}\).
We denote each \texttt{CSTR} instance by its size parameter \(N\), such as \texttt{CSTR30} for \(N=30\).
The DAEs of the \texttt{Dynamic} model are discretized with \texttt{pyomo.dae} by orthogonal collocation, with 30 finite elements per stage and three Radau collocation points~\citep{nicholson2018pyomo}, so each of its fixed-configuration subproblems is an NLP.
The reference configurations in external variable coordinates are \((N,N)\) for \texttt{CSTR}, \((13,4)\) for \texttt{Column} (physical trays \((14,5)\)), \((2,2,1)\) for \texttt{Small Batch}, and \((1,4)\) for \texttt{Dynamic}.
The full GDP formulations of the \texttt{CSTR}, \texttt{Column}, and \texttt{Small Batch} benchmarks are provided in Appendices A.1-A.3 of \citet{ovalle2025logic}, respectively, while the \texttt{Dynamic} formulation is described in Section IV.C of \citet{peng2024dynamic}.
External-variable mappings for the three design cases are provided in Appendices B.1, B.2, and B.4 of \citet{ovalle2025logic}.
The \texttt{CSTR} volume is reported in liters and the \texttt{Small Batch} investment in the source model's dollars.
The \texttt{Column} uses the source model's weighted cost proxy based on duties and active trays.
The model sizes and unit conventions are given in the Supplementary Section~\ref{sec:supp-models} (Table~\ref{tab:supp-model-sizes}).
\texttt{Dynamic} is an abstract benchmark with no calibrated physical units for state or time.
Minimizing its negative squared-state integral favors greater integrated state magnitude and tests mode sequencing coupled to dynamic equations.

Table~\ref{tab:process-cases} summarizes the original GDP and DFO candidate domains, the number of configurations whose objective values the interior-point NLP solver \texttt{IPOPT}~\citep{wachter2006implementation} returns with an accepted status, the number of configurations with measured evaluation times, and the evaluation budgets.
The 980 primary process runs, seven methods with ten seeds in each of the 14 model-budget conditions, evaluate 2,106 distinct configurations over all models and sizes.
The process searches use frozen objective tables from fixed-configuration \texttt{IPOPT} evaluations; only accepted finite objectives are used as valid lookup objectives.
The objective, acceptance status, and point cost are taken from the same recorded evaluation attempt for the same model instance and configuration.
Measured coverage includes all configurations in each DFO candidate box, including unsuccessful and structural-rejection evaluations.
The evaluation budget \(B\) limits the total number of objective evaluations, including initialization, repeated calls, and unsuccessful evaluations.
The full GDP box precedes the logical and algebraic feasibility restrictions.
\texttt{Column} searches \(\{8,\ldots,15\}\times\{1,\ldots,8\}\) (64 coordinates), the logical domain of its original \(\{1,\ldots,15\}^{2}\) box.
Other cases use their full boxes.
The \texttt{Column} bounds encode the reflux and boilup placement restrictions; algebraic constraints are still checked during configuration evaluation.
All DFO methods use this box, and the \texttt{Column} conclusions apply only to it.
For \texttt{CSTR}, recycle/feed ordering admits \(N(N+1)/2\) logical candidates, but not all yield objective values satisfying the solver acceptance criteria defined in Section~\ref{sec:protocol}.
The \texttt{CSTR} comparisons do not isolate scalability: evaluation budgets change with size.
\texttt{Dynamic} admits only configurations with \(\xi_1<\xi_2\) or \((\xi_1,\xi_2)=(9,9)\), which leaves 37 logically admissible configurations.
\texttt{Small Batch} uses all 27 configurations in \(\{1,2,3\}^{3}\), whose decisions count parallel mixer, reactor, and centrifuge units.
Of these configurations, 15 are rejected during feasibility-based bound tightening (FBBT) preprocessing.

\begin{table}[htbp]
\centering
\caption{Process GDP models: configuration domains, evaluation budgets, reference objectives, and the cost of full enumeration.
The full GDP box lists the external-variable values before logical restrictions, and the DFO domain is the part of it that the DFO methods search; it equals the full box except for \texttt{Column}, whose DFO domain is its logical domain \(\{8,\ldots,15\}\times\{1,\ldots,8\}\), and every configuration in it has a measured evaluation cost.
Accepted configurations are those whose fixed-configuration \texttt{IPOPT} solve returns a finite objective, violates no constraint or bound by more than \(10^{-5}\), and ends with an accepted solver status (Section~\ref{sec:protocol}), and the reference objective \(f_{\mathrm{ref}}\) is the lowest accepted objective of the DFO domain.
The enumeration time is the summed measured cost of evaluating every configuration of the DFO domain once.}
\label{tab:process-cases}
\footnotesize
\setlength{\tabcolsep}{3.5pt}
\begin{threeparttable}
\begin{tabular}{@{}llcrrrrrcr@{}}
\toprule
Model & \(N\) & Full GDP box & \makecell{Full box\\size} & \makecell{DFO domain\\size} & Accepted & \(B\) & \(f_{\mathrm{ref}}\) & Unit & \makecell{Enumeration\\time (s)} \\
\midrule
\multirow{6}{*}{\texttt{CSTR}} & 5 & \multirow{6}{*}{\(\{1,\ldots,N\}^{2}\)} & 25 & 25 & 15 & 8 & 3.06201 & L & 30.27 \\
 & 10 &  & 100 & 100 & 54 & 15 & 2.88954 & L & 122.06 \\
 & 15 &  & 225 & 225 & 120 & 22 & 2.82939 & L & 298.42 \\
 & 20 &  & 400 & 400 & 210 & 29 & 2.79941 & L & 553.38 \\
 & 25 &  & 625 & 625 & 325 & 36 & 2.78148 & L & 899.39 \\
 & 30 &  & 900 & 900 & 464 & 43 & 2.76957 & L & 1{,}719.00 \\
\midrule
\texttt{Column} & -- & \(\{1,\ldots,15\}^{2}\) & 225 & 64 & 29 & 13 & 19{,}346.1 & none & 137.24 \\
\midrule
\texttt{Small Batch} & -- & \(\{1,2,3\}^{3}\) & 27 & 27 & 12 & 10 & 167{,}428 & \$ & 10.18 \\
\midrule
\texttt{Dynamic} & -- & \makecell{\(\{1,\ldots,9\}\)\\\(\times\{2,\ldots,9\}\)} & 72 & 72 & 37 & \makecell{12,15,20,\\30,37,50} & -212.009 & none & 47.25 \\
\bottomrule
\end{tabular}
\begin{tablenotes}[flushleft]\footnotesize
\item \(N\): number of potential reactors of \texttt{CSTR}; \(B\): evaluation budget; \(f_{\mathrm{ref}}\): reference objective; --: not applicable.
Unit: L, liters of reactor volume; \$, dollars of investment; none, no physical unit (the weighted cost of duties and active trays of \texttt{Column} and the negative squared-state integral of \texttt{Dynamic}).
\end{tablenotes}
\end{threeparttable}
\end{table}
\FloatBarrier

\subsection{Simulation-based benchmark: lomustine manufacturing}
\label{sec:lomustine-case}
The \texttt{Lomustine} benchmark selects the flowsheet of a manufacturing process for the anticancer drug lomustine, whose model is a \texttt{PharmaPy} simulation~\citep{casasorozco2021pharmapy,casasorozco2023techno}.
The superstructure was introduced by \citet{laky2022superstructure}, related flowsheet synthesis frameworks are described by \citet{barhate2025hybrid}, and \citet{park2026ising} studied its discrete landscape with Ising-based solvers; we use the corrected-model data of \citet{secquoiaPdIsing}.
The superstructure contains two reactors in series, each a batch reactor, a CSTR, or a plug-flow reactor (PFR), followed by a batch vaporizer and a crystallizer that is either an unseeded batch cooling crystallizer or a train of one, two, or three mixed-suspension, mixed-product-removal (MSMPR) crystallizers.
The external variables \(\xi_1\) and \(\xi_2\) select the two reactors in the order (batch, CSTR, PFR), and \(\xi_3\) selects the crystallizer in the order (batch, one, two, three MSMPR stages), which gives \(3\times3\times4=36\) configurations, all logically admissible.
These choices are categorical, so this case does not meet the ordering condition of the external-variable reformulation; we index the alternatives in the order of the source enumeration, and the neighborhoods of \texttt{D-SDA} and \texttt{NOMAD} and the surrogate models act on this order.

Written as a GDP, the flowsheet is
\begin{equation}\label{eqn:lomustine-gdp}
\begin{aligned}
\min_{\mathbf u,\mathbf s,\mathbf Y}\quad & f(\mathbf u,\mathbf s)\\
\text{s.t.}\quad
&\bigveebar_{a\in\mathcal R}\left[\begin{gathered}Y^{1}_{a}\\ \mathbf s_1=\mathcal S^{\mathrm R}_{a}(\mathbf u,\mathbf s_0)\end{gathered}\right],\qquad
\bigveebar_{a\in\mathcal R}\left[\begin{gathered}Y^{2}_{a}\\ \mathbf s_2=\mathcal S^{\mathrm R}_{a}(\mathbf u,\mathbf s_1)\end{gathered}\right],\\
&\mathbf s_3=\mathcal S^{\mathrm V}(\mathbf u,\mathbf s_2),\qquad
\bigveebar_{c\in\mathcal C}\left[\begin{gathered}Y^{\mathrm C}_{c}\\ \mathbf s_4=\mathcal S^{\mathrm C}_{c}(\mathbf u,\mathbf s_3)\end{gathered}\right],\\
&\mathbf g(\mathbf u,\mathbf s)\leq\mathbf 0,\qquad \mathbf u\in\mathcal U,\qquad
\mathbf Y\in\{\mathrm{False},\mathrm{True}\}^{10},
\end{aligned}\tag{Lom}
\end{equation}
where \(\mathcal R=\{\text{batch},\text{CSTR},\text{PFR}\}\) and \(\mathcal C=\{\text{batch},\text{one},\text{two},\text{three MSMPR stages}\}\).
Each disjunct contains the \texttt{PharmaPy} model of one unit, which maps the operating inputs \(\mathbf u\) and the state \(\mathbf s_{k-1}\) of the incoming stream to the outlet state \(\mathbf s_k\), and \(\mathbf s_0\) is the reagent feed.
The vaporizer \(\mathcal S^{\mathrm V}\) is present in every configuration, \(f\) is the objective defined below, and \(\mathbf g\) is the solubility path constraint.
The external variables \(\xi_1\), \(\xi_2\), and \(\xi_3\) index the three disjunctions, and fixing them selects one unit model in each, which gives Problem~\eqref{eqn:sim-sub}.
The unit models are available only as simulations, so Problem~\eqref{eqn:lomustine-gdp} is a GDP whose disjuncts have no algebraic form for the conventional GDP methods.

For each configuration, \texttt{NOMAD} optimizes 12 to 16 normalized operating inputs, including reactor temperatures and residence times, the reagent concentration, the cooling profile, and the crystallizer residence times, with 60 simulations, which solves Problem~\eqref{eqn:sim-sub} approximately.
The objective \(-(6\times10^{4}\,\dot m+d_{10})\), with the production rate \(\dot m\) in kg/s and the number-weighted mean crystal size \(d_{10}\) in \(\mu\)m, rewards both quantities, and a solubility path constraint prevents crystallization before the crystallizer; all 36 configurations satisfy it.
These inner optimizations take 91 to 1,019~s per configuration, so evaluating every configuration takes 5.3~h.
All methods read the stored results of these inner optimizations, so they evaluate the same \(\Phi(\boldsymbol\xi)\), and each evaluation is charged its recorded inner-optimization time.
The reference configuration is a PFR followed by a batch reactor and the batch crystallizer, \(\boldsymbol\xi=(3,1,1)\), with \(f_{\mathrm{ref}}=-487.842\); the second-best configuration replaces the PFR by a CSTR and lies 0.8\% above it.
The budgets are \(B=10\), 15, and 20, with the shared initialization, seeds, unaccepted-response values, and target rule of the process benchmarks, and \texttt{D-SDA} uses the same unit \(L_\infty\) neighborhood.
Supplementary Section~\ref{sec:supp-lomustine} lists the 36 objective values and the data provenance.

\subsection{Catalytic distillation column with live NLP evaluations}
\label{sec:catdist-case}
The \texttt{Catalytic Column} benchmark designs a catalytic distillation column that produces ethyl \textit{tert}-butyl ether (ETBE) from ethanol and a mixture of butenes, for which \citet{linan2020optimalII} formulated a GDP model and solved it with \texttt{D-SDA}.
The column has a condenser, 20 internal stages, and a reboiler, and the reflux returns to the top internal stage.
The GDP decisions place the ethanol feed at stage \(f_1\), catalyst on three stages \(c_1<c_2<c_3\), the butenes feed at stage \(f_2\), and the boil-up return at stage \(b\), with \(2\le f_1\le c_1<c_2<c_3\le f_2\le b\le 21\).
Internal stages below the boil-up return do not exist, so \(b\) also sets the number of trays.
Each existing stage is an equilibrium stage with tray hydraulics, and each catalytic stage adds the reaction kinetics and a 120~\(^\circ\)C limit on the catalyst temperature.
Every decision chooses a stage from an ordered list, so this model meets the ordering condition of the external-variable reformulation~\citep{ovalle2025logic}.
We use a corrected version of the published model in which thermodynamic, property, hydraulic, and cost errors are fixed and the pure-component data are fitted to open reference data (Supplementary Section~\ref{sec:supp-catdist}).
The objective is the negative annual profit, in \$/yr, of a column whose bottoms product contains at least 95~mol\% ETBE, at 2010 prices; it is positive because the butenes left in the distillate earn no credit.
With the configuration fixed, the model is an NLP with about 2,600 variables and 2,700 equations, which \texttt{CONOPT}~\citep{drud1994conopt} solves from the initialization of the published model.

The search methods use the six gaps \(\mathbf g=(f_1-2,\;c_1-f_1,\;c_2-c_1-1,\;c_3-c_2-1,\;f_2-c_3,\;b-f_2)\in\{0,\ldots,6\}^6\) as external variables, so that every point of the box satisfies the order of the stages.
The GDP admits 100,947 configurations, and the bound of six on each gap, which keeps every design of a pilot within 5\% of the pilot's best, leaves 54,159 of them in the box.
Of the 117,649 points of this box, the 54,159 with \(\sum_i g_i\le 17\) are admissible; the others would place the boil-up return below the last internal stage and are charged to the budget without an NLP solve, like the logical rejections of \texttt{CSTR} and \texttt{Dynamic}.
Unlike in the other process cases, no stored table of responses exists: each evaluation solves the NLP of its configuration when the configuration is first requested.
Because these responses are deterministic, they are stored and reused when another run requests the same configuration, and each evaluation is charged the recorded solve time of its configuration.
A configuration is accepted when \texttt{CONOPT} returns an optimal or locally optimal solution with normal completion within 300~s.
The 3,499 configurations solved in the campaign took a median of 4.8~s and a mean of 7.9~s, so solving the 54,159 admissible configurations of the box would take about 119~h, and all 100,947 configurations of the GDP about 220~h.

The budgets are \(B=25\), 50, and 100 evaluations, and each seed has seven shared initial configurations drawn uniformly from the admissible configurations, one more than the number of external variables, which replace the three Latin hypercube points of the other benchmarks.
\texttt{D-SDA} starts from the best of them and uses the separable neighborhood, which changes one external variable by \(\pm1\) and corresponds to the 2-norm neighborhood of the original \texttt{D-SDA} study of this column~\citep{linan2020optimalII}.
The unit \(L_\infty\) neighborhood of the other benchmarks would contain 728 configurations in six dimensions, more than seven times the largest budget, whereas the separable neighborhood contains at most 12.
Because \texttt{D-SDA} moves in the gap coordinates, one move shifts every stage after the changed gap, so its neighbors differ from the neighbors in stage positions.
\texttt{RBFOpt}, \texttt{BAM}, \texttt{NOMAD}, and \texttt{MISO} keep the settings of the process benchmarks, \texttt{FMQA-SA} has a time limit of 7,200~s instead of 1,800~s because the solve times of its evaluations count toward it, and uniform sampling draws without replacement from its own random stream.
No global optimum is known, so the reference is the best accepted configuration found by any run or by a pilot of 202 configurations, \((f_1,c_1,c_2,c_3,f_2,b)=(8,8,12,17,17,21)\) with \(f_{\mathrm{ref}}=24{,}895.55\)~\$/yr, and the target rule is that of the process benchmarks.
Its first gap lies on the bound of the box, so we also solved the 747 configurations outside the box within two stages of it: none is better, and the best lies 0.14\% above it (Supplementary Section~\ref{sec:supp-catdist}).
The 210 runs ran as three concurrent single-core jobs on the server described in Section~\ref{sec:protocol} with the General Algebraic Modeling System (\texttt{GAMS})~\citep{gams2025} and \texttt{CONOPT}.

We also solved three MINLP forms of the column with \texttt{DICOPT}~\citep{viswanathan1990combined}, \texttt{SBB}~\citep{gamsSBB}, \texttt{BARON}~\citep{sahinidis1996baron}, and \texttt{SCIP}~\citep{bolusani2024scip} in \texttt{GAMS}, each on one core of the same server with a limit of 3,600~s, from the initialization and initial configuration of the published model.
These forms contain all 100,947 configurations of the GDP.
The first is the MINLP form published with the \texttt{D-SDA} study, which multiplies the switched equations by binary variables~\citep{linan2020optimalII}.
The other two are the Big-\(M\) and Hull reformulations that the \texttt{GAMS} extended mathematical programming framework~\citep{ferris2009emp} generates from a disjunctive form of the corrected model with one disjunction for each stage decision.
Both reformulations are valid only after two modeling steps.
Each nonlinear equation of a disjunct defines a residual variable on every stage, so that every disjunct constraint is linear and the Hull reformulation never evaluates terms such as \(1/\rho_L\) at zero.
Each disjunct constraint is also divided by an interval bound on its residual, so that \(M=1\) relaxes every constraint whose disjunct is inactive (Supplementary Section~\ref{sec:supp-catdist}).
These interval bounds are finite only if the liquid is denser than the vapor, so the two reformulations also bound the mass densities of the liquid (at least 200~kg/m\(^3\)) and the vapor (5 to 150~kg/m\(^3\)), which the GDP evaluations and the published form do not.
For any vapor composition, the vapor bounds cannot bind within the model's bounds on pressure, temperature, and compressibility factor, and every validated solution has a liquid density of at least 529~kg/m\(^3\).
With the binary variables fixed, the Big-\(M\) and Hull reformulations reproduce the GDP objective at the six feasible configurations we tested, the published form reproduces it at four of them and reaches other local optima at the remaining two, and all forms agree that three other configurations are locally infeasible.

\subsection{Controlled benchmark: Grid-GDP}
\label{sec:grid-gdp-formulation}
To study search behavior without NLP subproblems, we also consider 21 two-dimensional test functions adapted from the Quantum Hamiltonian Descent benchmark collection~\citep{qhdBenchmarks}.
Each is sampled at 31 equally spaced values per coordinate, giving 961 configurations.
For each benchmark function, let \([\ell_x,u_x]\) and \([\ell_y,u_y]\) denote the coordinate bounds, and let \(I=J=\{0,\ldots,30\}\) be the corresponding index sets.
The Grid optimizers receive \((i,j)\in\{0,\ldots,30\}^{2}\) directly.
Relative to the one-based convention in Section~\ref{sec:gdp-background}, \((\xi_1,\xi_2)=(i+1,j+1)\) is only a relabeling of alternatives.
\texttt{FMQA-SA} encodes each Grid index in five unsigned bits and rejects decoded values of 31.
The physical coordinates are defined as
\begin{equation}\label{eqn:grid-coordinates}
a_i=\ell_x+\frac{i(u_x-\ell_x)}{30},\qquad
b_j=\ell_y+\frac{j(u_y-\ell_y)}{30}.
\end{equation}

Let \(f\) denote the original analytical objective function and let \(L\) and \(U\) be its minimum and maximum over the \(31\times31\) lattice.
Both constants are obtained by exhaustive preprocessing of all 961 lattice points prior to any search.
A normalized target tolerance of \(10^{-5}\) corresponds to \(10^{-5}(U-L)\) in raw objective units, so equal normalized tolerances do not need to define equally large target sets.
Supplementary Table~\ref{tab:grid-domains} gives target-set sizes for all functions.
Full-precision normalization constants are available in the repository (Supplementary Section~\ref{sec:supp-data}).
Each benchmark function is represented by the following GDP:
\begin{equation}\label{eqn:grid-gdp}
\begin{aligned}
\min_{x,y,Y^x,Y^y}\quad & F(x,y)=\frac{f(x,y)-L}{U-L}\\
\mathrm{s.t.}\quad &
\bigvee_{i\in I}\left[\begin{array}{c}Y_i^x\\x=a_i\end{array}\right],\qquad
\bigvee_{j\in J}\left[\begin{array}{c}Y_j^y\\y=b_j\end{array}\right],\\
&\operatorname{ExactlyOne}(Y_0^x,\ldots,Y_{30}^x),\\
&\operatorname{ExactlyOne}(Y_0^y,\ldots,Y_{30}^y),\\
&\ell_x\leq x\leq u_x,\qquad \ell_y\leq y\leq u_y,\\
&Y_i^x,Y_j^y\in\{\mathrm{False},\mathrm{True}\}.
\end{aligned}
\end{equation}
The two exclusive disjunctions contain 62 coordinate-selection Boolean indicators.
The objective retains the explicit mathematical expression.
The conventional GDP methods therefore receive the expression and its disjunctive structure, whereas the DFO methods observe only the objective values at the coordinates they select.

Fixing both coordinate disjunctions fixes \(x\) and \(y\), leaving no remaining original continuous degrees of freedom.
The measured fixed-point \texttt{IPOPT}~\citep{wachter2006implementation} evaluations therefore consist of model preparation, solver invocation, and verification.
The unfixed objective functions remain nonlinear and generally nonconvex, so Big-\(M\) and Hull produce MINLP representations.
The complete mathematical definitions and coordinate domains of the 21 benchmark functions are provided in Supplementary Sections~\ref{sec:supp-domains}-\ref{sec:supp-objectives} (Table~\ref{tab:grid-domains})~\citep{rbfoptTestArchive,qhdBenchmarks}.

\subsection{Search protocol and compared algorithms}

\subsubsection{Shared initialization and evaluation budgets}
In the process, \texttt{Lomustine}, and Grid-GDP benchmarks, all configuration-search methods and sampling controls use precomputed configuration responses during the search; the \texttt{Catalytic Column} solves each configuration when it is first requested (Section~\ref{sec:catdist-case}).
For process cases, accepted responses are objective values from \texttt{IPOPT} under the initialization and acceptance rules in Section~\ref{sec:protocol}.
Because \texttt{IPOPT} returns local solutions, these values can exceed the subproblem value \(\Phi(\boldsymbol\xi)\).

Each Grid-GDP benchmark uses ten seeds, each defining a different initial design, and evaluation budgets \(B\in\{10,30,50,70,100\}\).
The primary Grid comparison comprises six search implementations and two sampling controls, giving \(21\times10\times5\times8=8{,}400\) runs.
The process benchmarks comprise 14 model-budget conditions, each evaluated with ten paired seeds, giving 140 runs per method and 980 runs across the six search implementations and uniform sampling (Table~\ref{tab:process-cases}).
The \texttt{Lomustine} benchmark adds three budgets with ten seeds each, 210 runs in total, and the \texttt{Catalytic Column} another 210 runs with seven initial configurations per seed (Section~\ref{sec:catdist-case}).
In all other benchmarks, for each problem-budget-seed combination, three distinct initial configurations are generated using seeded Latin hypercube sampling, mapped to integer coordinates, and shared in the same order across all methods and sampling controls.
This common design controls the first three evaluations; subsequent sampling, response handling, and stopping rules follow the individual implementations.
All process cases use seeds 0-9, with no replacement of failed starts.

Methods may terminate before exhausting \(B\); only Grid \texttt{D-SDA} receives the reference target as a stopping criterion.
Detailed settings and worker time limits are provided in the Supplementary Section~\ref{sec:supp-protocol}.
Every objective evaluation consumes budget, including initialization, repeated configurations, and unaccepted evaluations; internal cache lookups that do not request an evaluation are free.
Target attainment uses the recorded acceptance status and original objective values, while each method receives the following values for unaccepted evaluations.
In all process cases, \texttt{D-SDA}, \texttt{RBFOpt}, \texttt{NOMAD}, and \texttt{MISO} assign a penalty of \(10^6\) to unaccepted evaluations.
\texttt{BAM} represents these evaluations by \texttt{NaN}, indicating that no accepted objective value is available.
For surrogate training, \texttt{FMQA-SA} scales the accepted objective values to \([0,1]\) and assigns the unaccepted evaluations a value of \(2\).

\subsubsection{Budget-limited \texttt{D-SDA}} 
We use a budget-limited implementation of the previously developed \texttt{D-SDA}~\citep{linan2020optimal}, applying its neighborhood search and discrete line search over the GDP external integer lattice as in \texttt{LD-SDA}~\citep{ovalle2025logic}.
It evaluates the unit \(L_\infty\) neighborhood in fixed enumeration order, selects the best improving neighbor, and continues along the resulting direction through a discrete line search; on the six-dimensional \texttt{Catalytic Column}, it uses the separable neighborhood instead (Section~\ref{sec:catdist-case}), and Supplementary Section~\ref{sec:supp-dsda-neighborhood} repeats every other \texttt{D-SDA} run with the separable neighborhood.
Neighborhood and line-search moves require an absolute objective improvement greater than \(10^{-5}\).
As in the original \texttt{D-SDA} and \texttt{LD-SDA} implementations, a configuration is never evaluated twice: when a neighbor or line-search point was already evaluated in the same run, including the shared initial configurations, \texttt{D-SDA} reuses its stored value, and no evaluation is charged.
If the evaluation budget \(B\) is exhausted, the best solution found is retained.

For the process benchmarks, \texttt{D-SDA} starts with the best accepted configuration among the three initial points evaluated and stops when this descent ends.
For Grid-GDP benchmarks, initial configurations are ranked by objective value and successive descents are attempted while evaluation budget remains.
After each completed descent, the search stops if the best observed objective satisfies the common target criterion in Equation~\eqref{eqn:target}; otherwise, it tries the next ranked initial configuration.
For all process benchmarks, \texttt{D-SDA} searches within the declared integer bounds without precomputed feasibility filtering, using the evaluation accounting and unaccepted-response values stated above.
To check that the single process descent does not handicap \texttt{D-SDA}, we repeated the 140 process runs with the Grid restart rule, and with further descents from uniformly drawn unvisited configurations until the budget is exhausted (Supplementary Section~\ref{sec:supp-dsda-restarts}).

\subsubsection{\texttt{NOMAD}}
\label{sec:nomad-settings}
\texttt{NOMAD}~\citep{audet2022nomad}, called through \texttt{PyNomadBBO}, runs mesh adaptive direct search with all external variables declared integer.
It starts from the best of the three shared initial configurations, whose stored value it reuses without charge, with other settings at their defaults.
With integer variables, \texttt{NOMAD} stops when its mesh reaches the integer resolution, so, like \texttt{D-SDA}, it can end before the budget.
\texttt{NOMAD} sometimes terminated with a segmentation fault after its last new evaluation.
Such a run is repeated with larger evaluation caps, and a rerun is accepted only if it evaluates exactly the same configurations in the same order; 123 of the 1,190 runs completed this way.
Twenty-one Grid runs at \(B=100\) crashed on every attempt after identical evaluations; they are scored on the evaluations completed before the crash, and 16 of them had already reached the target (Supplementary Section~\ref{sec:supp-nomad-miso}).
The \texttt{Lomustine} harness repeats a crashed run once; 12 runs of its option-order and capital-cost repetitions, none in the source order of Table~\ref{tab:lomustine}, crashed on both attempts after identical evaluations and are scored in the same way (Supplementary Section~\ref{sec:supp-lomustine}).

\subsubsection{\texttt{RBFOpt}}
\texttt{RBFOpt} uses MSRSM, automatic RBF selection, and the default genetic candidate-search option~\citep{costa2018rbfopt,nannicini2021implementation}.
All external variables are integers, and the three shared initial configurations are supplied with the additional package initialization disabled.
Search stops at the evaluation budget or the native iteration and stall limits.
Other surrogate-minimization operations use standalone \texttt{Bonmin}~\citep{bonami2008algorithmic} and \texttt{IPOPT}~\citep{wachter2006implementation} executables.
These calls contribute to search time, but not the objective-evaluation budget.

\subsubsection{\texttt{MISO}}
\texttt{MISO} uses the CPTV strategy of \citet{muller2016miso} as implemented in \texttt{soogo}, with a cubic radial basis function model with a linear polynomial tail, integer coordinates, and a median low-pass filter on large objective values.
The three shared configurations and their responses initialize the model.
Because \texttt{soogo} does not accept pure-integer domains in its coordinate-perturbation step, our adapter applies the same dynamically dimensioned search perturbation on the integer lattice and disables the continuous local search, which has no continuous coordinates to act on.
A linear tail needs \(d+1\) affinely independent points, so when the shared initial configurations do not span the space, as with the three initial configurations of the three-dimensional \texttt{Small Batch} and \texttt{Lomustine} cases, uniformly drawn admissible configurations, charged to the budget, are added until they do.
A run stops at the budget, or after ten consecutive iterations without a new evaluation when no unevaluated candidate remains.

\subsubsection{\texttt{BAM}}
\texttt{BAM} uses local algebraic surrogate models optimized by \texttt{BARON} to generate candidate configurations~\citep{ma2023bam,bamManual}.
All external variables are integers.
The three shared configurations and their responses are supplied without an additional starting point.
The remaining evaluation budget is \(B-3\), with other settings at their defaults.
On Grid-GDP, the next ten planned evaluations use a Sobol sequence, which distributes points across the domain, subject to the remaining budget.
In the process experiments, \texttt{BAM} stops at its own evaluation limit after the \(B\)-th evaluation.

\subsubsection{\texttt{FMQA-SA}}
\texttt{FMQA-SA} uses the FM training and QUBO construction code of \citet{secquoiaBboFmqa}.
At each iteration, a new second-order rank-eight FM is fitted to all evaluated configurations.
Finite precomputed responses are min-max scaled to \([0,1]\) (or zero when all finite values are identical), while non-finite evaluations receive a training target of \(2.0\).

Integer coordinates are encoded in unsigned binary, and the fitted FM is converted to a QUBO.
Binary candidates are generated on CPU using D-Wave's classical simulated annealing sampler~\citep{dwaveSamplers}.
The lowest-energy decoded candidate within the coordinate bounds that has not been previously evaluated is selected; if none is available, a candidate is drawn uniformly from the remaining unevaluated configurations.
Candidates are not screened for the feasibility of logical or nonlinear constraints before evaluation.
The new observation is appended, and FM fitting and candidate generation repeat until the evaluation budget \(B\) or the 1,800~s method time limit is reached.

\subsubsection{Sampling controls and full enumeration}
\label{sec:sampling-controls}
Uniform sampling evaluates the same three initial configurations and then draws uniformly without replacement from the remaining declared integer domain.
It uses a separate deterministic random stream, no objective-dependent selection, and no logical or NLP-feasibility screening before evaluation.
The Grid \texttt{BAM}-prefix control repeats \texttt{BAM}'s ten predetermined sampling points in the same order after the three shared starts, then samples uniformly without replacement, stopping at the budget.
This control isolates the contribution of the fixed initial design of \texttt{BAM}.
Full enumeration evaluates every configuration of the DFO domain once.
It reaches the target by construction, because \(f_{\mathrm{ref}}\) is the best accepted value among these configurations, and its cost is the sum of their measured point costs.

\subsection{Execution environment and conventional GDP comparisons}
\label{sec:protocol}
All results reported here come from a server with an AMD EPYC 7643 processor (48 cores, 2.3~GHz) and 1~TB RAM; two earlier runs of \texttt{Catalytic Column} checks on a workstation serve only as comparisons (Supplementary Section~\ref{sec:supp-catdist}).
All optimization methods and numerical libraries were configured for single-threaded execution.
The software comprised \texttt{Python} 3.11.16, \texttt{Pyomo} 6.10.0~\citep{bynum2021pyomo}, \texttt{RBFOpt} 4.3.0, \texttt{BAM} 5.0.0, \texttt{soogo} 2.1.0, \texttt{NOMAD} 4.6.0 called through \texttt{PyNomadBBO} 4.6.0, \texttt{PyTorch} 2.8.0~\citep{paszke2019pytorch}, \texttt{Bonmin} 1.8.9, \texttt{IPOPT} 3.14.19 with \texttt{MUMPS} 5.8.1~\citep{amestoy2001mumps}, and \texttt{CONOPT} 4.39.1.
The process and Grid comparisons called \texttt{IPOPT} and \texttt{SCIP} 9.2.4~\citep{bolusani2024scip} through \texttt{GAMS} 52.2.0.
The \texttt{Catalytic Column} used \texttt{GAMS} 54.3.1 with \texttt{SCIP} 10.0.3 and \texttt{BARON} 26.5.27 (Section~\ref{sec:catdist-case}).

\paragraph{Fixed-disjunction NLP evaluations.}
We fix the original Boolean GDP decisions and evaluate the resulting subproblems using \texttt{IPOPT} with the \texttt{MUMPS} linear solver.
For the Grid-GDP benchmarks, all 20,181 configurations were evaluated individually, while the process benchmarks used fixed-configuration evaluations of the corresponding process models.
\texttt{IPOPT} used a tolerance of \(10^{-8}\) and a solver time limit of 60~s per configuration.
Process acceptance requires a finite objective, a maximum absolute violation of any constraint or variable bound not greater than \(10^{-5}\), and the accepted solver and termination statuses specified in the Supplementary Section~\ref{sec:supp-acceptance}.
Time-limited candidates remain unaccepted even when feasible.
\texttt{CSTR} sizes 5-30 use the mass-balance/zero-recycle initialization.
\texttt{Column} and \texttt{Small Batch} use their saved initialization files, and \texttt{Dynamic} uses fresh clones of the unsolved discretized template.

\paragraph{Conventional GDP comparisons.}
Conventional comparisons include the \texttt{GDPopt} methods GLOA, LOA, and RIC, as well as the Big-\(M\) and Hull reformulations.
\texttt{GDPopt} uses \texttt{SCIP} for MILP master problems and \texttt{IPOPT} for NLP subproblems~\citep{chen2022pyomogdp,pyomoGDPopt}.
Because fixing both Grid disjunctions also fixes \(x\) and \(y\), the Grid subproblems are solved as linear problems by \texttt{SCIP}, and \texttt{IPOPT} is not called in those runs.
The Big-\(M\) and Hull reformulations are solved directly using \texttt{SCIP}.
All conventional GDP methods use the original GDP domains, a time limit of 1,800~s per instance, and a relative optimality gap of \(10^{-5}\); the MINLP solvers of the \texttt{Catalytic Column} use a limit of 3,600~s and a relative gap of zero (Section~\ref{sec:catdist-case}).

\paragraph{Artificial intelligence assistance.}
\texttt{ChatGPT} and \texttt{Codex} (OpenAI, including GPT-6 Astra) and \texttt{Claude Code} (Anthropic, Claude Opus 5.5) assisted with developing and revising experiment, analysis, and plotting code.
The authors reviewed and revised the resulting code and validated the computational analyses.
Data visualizations were generated by the documented plotting programs from the archived numerical results.
The authors remain responsible for the computational methods and reported findings.

\subsection{Target attainment and cumulative success profiles}
\label{sec:metrics}
Let \(f_{\mathrm{best}}\) denote the lowest feasible objective value observed so far within the evaluation budget, and let \(f_{\mathrm{ref}}\) denote the reference objective value.
For process benchmarks, \(f_{\mathrm{best}}\) is computed from evaluations that meet the solver acceptance criteria in Section~\ref{sec:protocol}.
A run reaches the reference target when
\begin{equation}\label{eqn:target}
 |f_{\mathrm{best}}-f_{\mathrm{ref}}|
 \leq 10^{-5}\max\{1,|f_{\mathrm{ref}}|\}.
\end{equation}
For the process experiments, \(f_{\mathrm{ref}}\) is the lowest accepted objective value among the precomputed fixed-configuration NLP evaluations, so reaching the target means matching the best configuration that full enumeration of the stored responses would find.
Runs that terminate with an execution failure are counted as unsuccessful, even if the target was observed earlier; the exception is a \texttt{NOMAD} run that crashes on every attempt after evaluating the same configurations, which is scored on the evaluations completed before the crash (Section~\ref{sec:nomad-settings}).
For process runs, the reported solution is the best accepted configuration recorded within budget, including the final evaluation.
Reported objectives and bounds of the process benchmarks use six significant figures; all target and bound checks use unrounded values.
For the full-process-GDP comparison, a feasible objective at or below \(f_{\mathrm{ref}}+10^{-5}\max\{1,|f_{\mathrm{ref}}|\}\) counts as achievement, so solutions better than the reference are also included.

For each evaluation budget \(B\), let \(\mathcal R_B\) denote the set of problem-seed instances compared across methods, and let \(s\) identify a solution method.
For each instance \(r\in\mathcal R_B\), \(H_{sBr}\) equals one if the run reaches the target and zero otherwise.
Let \(E_{sBr}\) be the index of the first objective evaluation at which the target is reached, with \(E_{sBr}=+\infty\) if the target is not reached.
At a threshold of \(N\) objective evaluations, where \(0\leq N\leq B\), the absolute evaluation-count profile is
\begin{equation}\label{eqn:absolute-evaluations}
 C^{\mathrm{eval}}_{sB}(N)
 =\frac{1}{|\mathcal R_B|}\sum_{r\in\mathcal R_B}
 \mathbf{1}\{H_{sBr}=1,\ E_{sBr}\leq N\}.
\end{equation}
Here, \(\mathbf{1}\{\cdot\}\) equals one when the stated conditions hold and zero otherwise, and \(|\mathcal R_B|\) is the number of problem-seed instances.
Thus, \(C^{\mathrm{eval}}_{sB}(N)\) gives the fraction of runs that reach the target within \(N\) evaluations.
For the Grid-GDP benchmarks, \(|\mathcal R_B|=210\) at each budget; for each process benchmark and budget, \(|\mathcal R_B|=10\).
All runs, including target misses, remain in the denominator.

Experiments at different evaluation budgets are analyzed separately, so a larger-budget trajectory need not extend a smaller-budget one.
These profiles count raw evaluations, without the best-solver normalization of performance profiles~\citep{dolan2002benchmarking} or the dimension normalization of data profiles~\citep{more2009benchmarking}.

\paragraph{NLP subproblems solved.}
An evaluation counts as an NLP subproblem solve when its configuration passes logical and FBBT preprocessing and \texttt{IPOPT} is called, whether or not the result is accepted.
For \texttt{GDPopt}, we count the \texttt{IPOPT} calls confirmed in its execution logs.
Internal NLP counts are unavailable for the direct Big-\(M\) and Hull solves.

\subsection{Runtime accounting and uncertainty}
\label{sec:cost-accounting}
For a recorded lookup-based run with \(n\) objective evaluations and measured lookup-call times, the estimated full-run completion time is
\begin{equation}\label{eqn:reconstructed-cost}
 \widehat T_{\mathrm{complete}}
 =T_{\mathrm{recorded}}-\sum_{q=1}^{n}t_q^{\mathrm{lookup}}
  +\sum_{q=1}^{n}c(\boldsymbol\xi^{(q)})+t_{\mathrm{prep}}.
\end{equation}
Here, \(T_{\mathrm{recorded}}\) is the measured runtime of the search based on lookup, \(t_q^{\mathrm{lookup}}\) is the measured lookup-processing interval of its \(q\)th objective call, and \(\boldsymbol\xi^{(q)}\) is the external integer configuration evaluated in that call.
The point cost \(c(\boldsymbol\xi^{(q)})\) is measured separately for the same benchmark model, problem size and configuration, including configuration-specific preparation, preprocessing, invocation of the solver when required, and validation.
The term \(t_{\mathrm{prep}}\) accounts for the preparation of a single model that is not already included in the point costs.
Search time includes surrogate fitting, candidate generation, auxiliary optimization, and recording performed within the method clock.
Throughout this manuscript, full-run completion denotes completion of this recorded search stage with reconstructed evaluation costs, excluding final result extraction and export.
Every recorded objective call contributes its point cost, including initialization, repeated visits, and preprocessing rejections; the same stored measurement is used each time a configuration is revisited.
Completion times therefore combine recorded search paths with one stored cost measurement per configuration.

\paragraph{Grid-GDP timers.}
Grid-GDP enumeration and objective normalization are excluded from the timing analysis, and \(t_{\mathrm{prep}}=0\) because configuration-specific preparation is included in point costs.
Recorded lookup-call times are subtracted from the method runtime for all compared Grid methods.
The 21 \texttt{NOMAD} runs that crashed on every attempt have no recorded times and are left out of the Grid timing summaries.

\paragraph{Process timers.}
Structurally rejected configurations incur measured preparation and rejection costs without an NLP solution (Supplementary Section~\ref{sec:supp-acceptance}).
For \texttt{Dynamic}, the original model-template construction and discretization cost is included once per run through \(t_{\mathrm{prep}}\); additional campaign-level preparation is not charged repeatedly.
All in-budget evaluations have measured point costs, and estimated full-run completion times are available for all 980 primary process runs.

Conventional GDP completion times are measured separately and are not augmented with the DFO point costs.
For Grid-GDP, conventional GDP timings include per-run preparation, reformulation, solution, and validation, but exclude shared-template preparation and lattice enumeration.
Full-process-GDP timings include model-template preparation.

\paragraph{Statistical comparisons.}
With ten seeds per process condition, we compare two methods with the two-sided exact McNemar test on their seed-paired target outcomes, and we call differences that this test does not resolve at the 0.05 level comparable (Supplementary Table~\ref{tab:supp-paired}).
The tests of \texttt{RBFOpt} against \texttt{D-SDA} are descriptive and not adjusted for multiple comparisons, and the 70 exploratory tests that involve \texttt{NOMAD} or \texttt{MISO} use a Holm correction (Supplementary Section~\ref{sec:supp-nomad-miso}).
Grid comparisons use a percentile bootstrap that resamples the 21 functions (Supplementary Section~\ref{sec:supp-additional}).

\section{Computational Results and Discussion}
\label{sec:results}

\subsection{Process GDP case studies}
\label{sec:process-results}

\subsubsection{Target attainment}
Table~\ref{tab:process-targets} reports target attainment for the 14 process conditions.
Across the 60 paired \texttt{CSTR} runs, \texttt{RBFOpt} reaches the target in 45, \texttt{MISO} in 43, \texttt{NOMAD} in 21, and \texttt{D-SDA} in 16.
\texttt{RBFOpt} succeeds in 30 runs in which \texttt{D-SDA} misses, and \texttt{D-SDA} in one in which \texttt{RBFOpt} misses (exact paired test \(p<10^{-7}\); Supplementary Table~\ref{tab:supp-paired}).
With ten seeds per size, the difference is significant at the 0.05 level only for \texttt{CSTR5} and \texttt{CSTR30}, but it has the same sign at every size.
On \texttt{Dynamic} at \(B=12\), \texttt{RBFOpt} reaches all ten targets and \texttt{D-SDA} four (\(p=0.031\)); \texttt{D-SDA} reaches all ten from \(B=20\) onward.

\texttt{MISO} shows the same advantage over \texttt{D-SDA} with a different profile across sizes.
It succeeds in 28 \texttt{CSTR} runs in which \texttt{D-SDA} misses, against one in the other direction (\(p<10^{-6}\)), and it reaches all 30 targets on \texttt{CSTR20} through \texttt{CSTR30}, but only 13 of 30 on \texttt{CSTR5} through \texttt{CSTR15}, where \texttt{RBFOpt} reaches 21.
Pooled over the \texttt{CSTR} sizes, \texttt{MISO} and \texttt{RBFOpt} disagree in 16 runs, seven in favor of \texttt{MISO} (\(p=0.80\)).
\texttt{NOMAD} reaches more \texttt{CSTR} targets than \texttt{D-SDA} (7 against 2 discordant runs, \(p=0.18\)) but fewer than \texttt{RBFOpt} (5 against 29, \(p<10^{-4}\)).
Like \texttt{D-SDA}, it stops before the budget in many runs (77 of 140) because its mesh reaches the integer resolution.
Exploratory per-condition tests of \texttt{NOMAD} and \texttt{MISO} against \texttt{D-SDA}, \texttt{RBFOpt}, and each other find no contrast that remains significant after a Holm correction over the 70 tests (Supplementary Section~\ref{sec:supp-nomad-miso}).

On the two smallest design models, \texttt{D-SDA} and \texttt{RBFOpt} are comparable.
\texttt{D-SDA} reaches five \texttt{Column} targets and \texttt{RBFOpt} four (\(p=1\)), and on \texttt{Small Batch} they reach eight and five (\(p=0.38\)).
\texttt{MISO} reaches all ten \texttt{Column} targets and \texttt{NOMAD} seven; on \texttt{Small Batch} they reach seven and six.
\texttt{BAM} reaches all ten \texttt{CSTR5} targets, compared with eight for \texttt{RBFOpt} (\(p=0.50\)).
The only \texttt{Column} hits of \texttt{BAM} and uniform sampling occur in seed 6, whose shared initial design contains the target.

The single-descent \texttt{D-SDA} stops before exhausting its budget in 87 of the 140 process runs.
Restarting it from the remaining shared initial configurations, or from uniformly drawn unvisited configurations until the budget is spent, adds one target, raising target attainment from 79 to 80 of 140 runs; both restart variants improve the final best objective in 14 runs (Supplementary Section~\ref{sec:supp-dsda-restarts}).

\texttt{FMQA-SA} reaches 42 of 140 process targets, compared with 37 for uniform sampling.
In the small process domains, the 100 annealed samples often contain no new configuration within the bounds, so 52.9\% of the \texttt{FMQA-SA} proposals come from its uniform fallback (Supplementary Section~\ref{sec:supp-proposals}).

Supplementary Table~\ref{tab:process-search} reports the search times, which exclude configuration-evaluation costs.
Median search times stay below 0.06~s for \texttt{D-SDA} and uniform sampling and below 1~s for \texttt{BAM}.
For \texttt{RBFOpt}, they grow with the budget, from 0.06~s on \texttt{CSTR5} to 26.3~s on \texttt{Dynamic} at \(B=50\).
\texttt{MISO} needs 0.66 to 158~s because it refits its model and scores a large candidate pool at every iteration, whereas \texttt{NOMAD} stays below 2~s.
\texttt{FMQA-SA} needs 3.7 to 29.8~s because it refits its factorization machine after every evaluation.
\begin{table}[htbp]
\centering
\caption{Process GDP models: target hits and median full-run completion times of the six DFO methods and uniform sampling, with ten runs per method and condition.
A target hit is a run whose best accepted objective lies within \(10^{-5}\max\{1,|f_{\mathrm{ref}}|\}\) of the reference objective \(f_{\mathrm{ref}}\) of Table~\ref{tab:process-cases} (Equation~\eqref{eqn:target}).
The completion time of a run adds its search time to the measured cost of every evaluated configuration and, for \texttt{Dynamic}, one model preparation (Equation~\eqref{eqn:reconstructed-cost}); medians are over all ten runs, including target misses.
Table~\ref{tab:process-cases} gives the enumeration times, and Supplementary Table~\ref{tab:process-search} the search times.}
\label{tab:process-targets}
\footnotesize
\setlength{\tabcolsep}{3.5pt}
\begin{threeparttable}
\begin{tabular}{@{}lrrlrrrrrrr@{}}\toprule
Model & \(N\) & \(B\) & Measure & \texttt{D-SDA} & \texttt{NOMAD} & \texttt{RBFOpt} & \texttt{MISO} & \texttt{BAM} & \texttt{FMQA-SA} & Uniform \\\midrule
\multirow{12}{*}{\texttt{CSTR}} & \multirow{2}{*}{5} & \multirow{2}{*}{8} & Target hits (of 10) & 2 & 4 & 8 & 4 & \textbf{10} & 3 & 4 \\
 & &  & Median completion time (s) & 10.67 & 10.71 & 10.72 & 11.33 & 10.67 & 13.92 & \textbf{9.44} \\
\addlinespace[2pt]
 & \multirow{2}{*}{10} & \multirow{2}{*}{15} & Target hits (of 10) & 3 & 3 & \textbf{6} & 2 & 1 & 2 & 1 \\
 & &  & Median completion time (s) & \textbf{17.61} & 21.25 & 21.79 & 21.87 & 18.54 & 27.07 & 17.75 \\
\addlinespace[2pt]
 & \multirow{2}{*}{15} & \multirow{2}{*}{22} & Target hits (of 10) & 2 & 2 & \textbf{7} & \textbf{7} & 1 & 1 & 2 \\
 & &  & Median completion time (s) & \textbf{25.68} & 33.92 & 37.59 & 39.79 & 27.67 & 42.82 & 27.84 \\
\addlinespace[2pt]
 & \multirow{2}{*}{20} & \multirow{2}{*}{29} & Target hits (of 10) & 3 & 2 & 7 & \textbf{10} & 2 & 1 & 2 \\
 & &  & Median completion time (s) & \textbf{32.31} & 38.17 & 45.03 & 56.03 & 39.37 & 59.10 & 40.12 \\
\addlinespace[2pt]
 & \multirow{2}{*}{25} & \multirow{2}{*}{36} & Target hits (of 10) & 3 & 6 & 8 & \textbf{10} & 2 & 1 & 1 \\
 & &  & Median completion time (s) & \textbf{45.49} & 56.56 & 62.29 & 78.56 & 50.47 & 74.70 & 51.17 \\
\addlinespace[2pt]
 & \multirow{2}{*}{30} & \multirow{2}{*}{43} & Target hits (of 10) & 3 & 4 & 9 & \textbf{10} & 2 & 1 & 1 \\
 & &  & Median completion time (s) & \textbf{47.17} & 71.23 & 101.96 & 107.33 & 115.41 & 101.19 & 65.20 \\
\midrule
\multirow{2}{*}{\texttt{Column}} & \multirow{2}{*}{--} & \multirow{2}{*}{13} & Target hits (of 10) & 5 & 7 & 4 & \textbf{10} & 1 & 4 & 1 \\
 & &  & Median completion time (s) & 35.44 & 31.93 & 54.60 & 31.98 & 27.50 & 41.82 & \textbf{21.69} \\
\midrule
\multirow{2}{*}{\texttt{Small Batch}} & \multirow{2}{*}{--} & \multirow{2}{*}{10} & Target hits (of 10) & \textbf{8} & 6 & 5 & 7 & 2 & 6 & 4 \\
 & &  & Median completion time (s) & 4.23 & 4.97 & 5.12 & 5.62 & 4.23 & 9.09 & \textbf{3.49} \\
\midrule
\multirow{12}{*}{\texttt{Dynamic}} & \multirow{2}{*}{--} & \multirow{2}{*}{12} & Target hits (of 10) & 4 & 8 & \textbf{10} & \textbf{10} & 0 & 2 & 0 \\
 & &  & Median completion time (s) & 12.22 & 12.36 & 12.29 & 11.80 & 9.95 & 14.47 & \textbf{7.80} \\
\addlinespace[2pt]
 &  & \multirow{2}{*}{15} & Target hits (of 10) & 6 & \textbf{10} & \textbf{10} & \textbf{10} & 3 & 3 & 3 \\
 & &  & Median completion time (s) & 15.32 & 15.48 & 16.10 & 15.44 & 11.53 & 18.39 & \textbf{8.71} \\
\addlinespace[2pt]
 &  & \multirow{2}{*}{20} & Target hits (of 10) & \textbf{10} & \textbf{10} & \textbf{10} & \textbf{10} & 3 & 3 & 3 \\
 & &  & Median completion time (s) & 20.02 & 20.85 & 20.54 & 42.42 & 14.07 & 25.23 & \textbf{12.27} \\
\addlinespace[2pt]
 &  & \multirow{2}{*}{30} & Target hits (of 10) & \textbf{10} & \textbf{10} & \textbf{10} & \textbf{10} & 4 & 4 & 4 \\
 & &  & Median completion time (s) & 20.02 & 22.00 & 34.19 & 106.24 & 21.91 & 36.24 & \textbf{19.38} \\
\addlinespace[2pt]
 &  & \multirow{2}{*}{37} & Target hits (of 10) & \textbf{10} & \textbf{10} & \textbf{10} & \textbf{10} & 7 & 5 & 5 \\
 & &  & Median completion time (s) & \textbf{20.02} & 22.00 & 42.42 & 111.61 & 24.99 & 46.05 & 24.38 \\
\addlinespace[2pt]
 &  & \multirow{2}{*}{50} & Target hits (of 10) & \textbf{10} & \textbf{10} & \textbf{10} & \textbf{10} & \textbf{10} & 6 & 6 \\
 & &  & Median completion time (s) & \textbf{20.02} & 21.99 & 65.14 & 195.36 & 33.98 & 61.31 & 33.12 \\
\bottomrule
\end{tabular}
\begin{tablenotes}[flushleft]\footnotesize
\item \(N\): number of potential reactors of \texttt{CSTR}; \(B\): evaluation budget; --: not applicable; Uniform: uniform sampling without replacement.
Bold: most target hits and shortest median completion time in each row, including ties.
\end{tablenotes}
\end{threeparttable}
\end{table}

Supplementary Figure~\ref{fig:process-evaluations} shows the evaluation at which each run that reaches the target first reaches it.

\subsubsection{NLP subproblems solved}
\begin{figure}[htbp]
\centering\includegraphics[width=\linewidth,height=0.72\textheight,keepaspectratio]{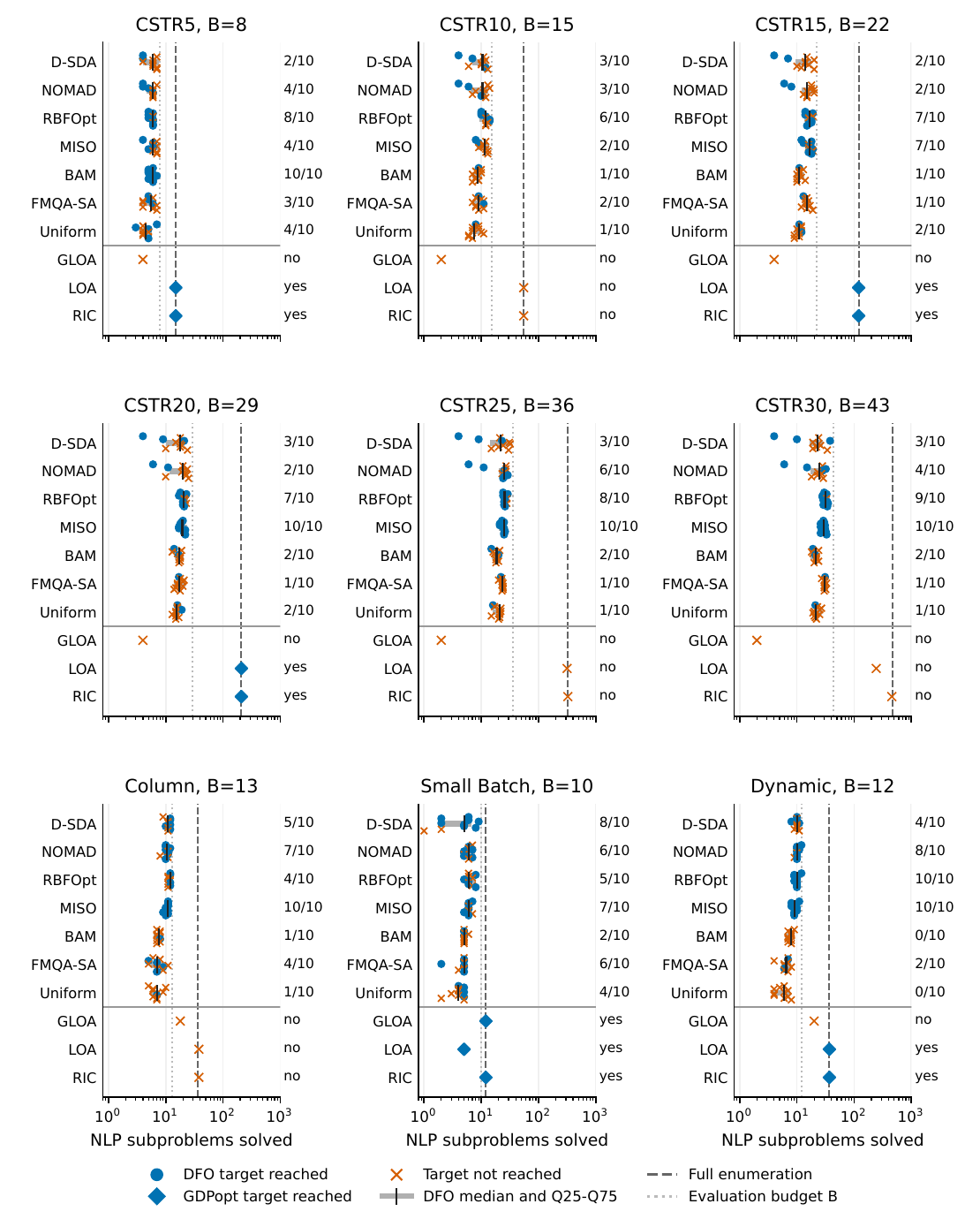}
\caption{Process GDP models: NLP subproblems solved per run by the six DFO methods and uniform sampling (ten runs each, at the evaluation budget \(B\) of each panel title) and by the single run of each \texttt{GDPopt} method: GLOA (global logic-based outer approximation), LOA (logic-based outer approximation), and RIC (relaxation with integer cuts).
Dots mark DFO runs that reach the reference target of Equation~\eqref{eqn:target}, an objective within \(10^{-5}\max\{1,|f_{\mathrm{ref}}|\}\) of the reference objective \(f_{\mathrm{ref}}\) of Table~\ref{tab:process-cases}; diamonds mark \texttt{GDPopt} runs whose feasible objective is at most \(f_{\mathrm{ref}}+10^{-5}\max\{1,|f_{\mathrm{ref}}|\}\) (Section~\ref{sec:metrics}); crosses mark the other runs; gray bars span the 25th to 75th percentiles of the ten DFO seeds, with a tick at the median.
The right-hand labels give the target hits out of ten of each DFO method and the outcome, yes or no, of each \texttt{GDPopt} run.
The dashed line marks the NLP subproblems needed to evaluate every configuration of the DFO domain, and the dotted line the evaluation budget \(B\).
Configurations rejected by logical or feasibility-based bound tightening (FBBT) preprocessing are not counted as NLP subproblems.
\texttt{GDPopt} searches the full \texttt{Column} box of 225 configurations, so its count can exceed the enumeration line of the 64-configuration DFO domain.
The direct Big-\(M\) and Hull solves are not shown because their internal NLP counts are unavailable.
Uniform: uniform sampling without replacement.}
\label{fig:process-nlp}
\end{figure}

Figure~\ref{fig:process-nlp} compares the number of NLP subproblems that each run solves with the number that full enumeration of the configuration domain requires.
On \texttt{CSTR5} through \texttt{CSTR20} and on \texttt{Dynamic}, RIC and LOA solve exactly as many subproblems as full enumeration, and their execution logs show that every subproblem fixes a different configuration (Supplementary Section~\ref{sec:supp-spread-overhead}).
RIC also enumerates all 325 admissible \texttt{CSTR25} configurations and solves 455 of the 465 \texttt{CSTR30} subproblems before the time limit (Supplementary Table~\ref{tab:joint-nlp-workload}).
The DFO methods solve at most 38 subproblems in any run, and \texttt{MISO} and \texttt{RBFOpt} reach the \texttt{CSTR30} target in ten and nine of ten runs with medians of 29.5 and 31.5, 6 to 7\% of the enumeration count.

Enumeration does not guarantee that a logic-based method reaches the target.
RIC and LOA miss the \texttt{CSTR10} target after solving all 55 of its subproblems, and RIC misses the \texttt{CSTR25} target after solving all 325.
These misses come from local NLP solutions obtained with the \texttt{GDPopt} initialization, which differs from the initialization used to build the lookup tables (Supplementary Section~\ref{sec:supp-acceptance}).
GLOA solves two to four subproblems per \texttt{CSTR} model and misses every \texttt{CSTR} target.

\subsubsection{Completion time}
In wall time, the advantage of fewer subproblem solves depends on how much each solve costs.
Configuration evaluations have median costs of 0.4 to 2.5~s in these models (Supplementary Table~\ref{tab:supp-overhead}).
Full enumeration therefore takes from 10.18~s for \texttt{Small Batch} to 1,719.00~s for \texttt{CSTR30} (Table~\ref{tab:process-cases}).
Table~\ref{tab:process-targets} gives the median completion times of the DFO methods.
On \texttt{Small Batch}, the DFO medians of 3.49 to 9.09~s are close to the enumeration time, and on \texttt{Dynamic} at \(B=50\), \texttt{RBFOpt} and \texttt{MISO} need 65.14 and 195.36~s, more than enumerating all 72 configurations.
On \texttt{CSTR30}, \texttt{MISO} and \texttt{RBFOpt} reach the target in ten and nine of ten runs with medians of 107.33 and 101.96~s, about one sixteenth of the enumeration time, while \texttt{NOMAD} and \texttt{D-SDA} reach four and three targets with medians of 71.23 and 47.17~s.

Search overhead and stopping rules also shape completion times.
For \texttt{Dynamic} at \(B=50\), search overhead contributes 42.7\% of the mean completion time of \texttt{RBFOpt}, compared with less than 0.1\% for \texttt{D-SDA}, and \texttt{MISO} spends a median of 158~s of its 195.36~s in search.
Across these ten runs, \texttt{D-SDA} and \texttt{NOMAD} evaluate 167 and 170 configurations because they stop at a local optimum or at the integer mesh resolution, whereas \texttt{RBFOpt}, \texttt{MISO}, \texttt{BAM}, and \texttt{FMQA-SA} each evaluate 500.
At \(B=30\), \texttt{RBFOpt} reaches all ten \texttt{Dynamic} targets in fewer evaluations than \texttt{D-SDA} (Supplementary Figure~\ref{fig:process-evaluations}), but it continues to the budget, and its median completion time is 34.19~s, compared with 20.02~s for \texttt{D-SDA}, 22.00~s for \texttt{NOMAD}, and 106.24~s for \texttt{MISO}.
On \texttt{Column}, where \texttt{MISO} reaches all ten targets, the median completion times are 31.93~s for \texttt{NOMAD}, 31.98~s for \texttt{MISO}, 35.44~s for \texttt{D-SDA}, and 54.60~s for \texttt{RBFOpt}.
\begin{figure}[htbp]
\centering\includegraphics[width=\linewidth,height=0.69\textheight,keepaspectratio]{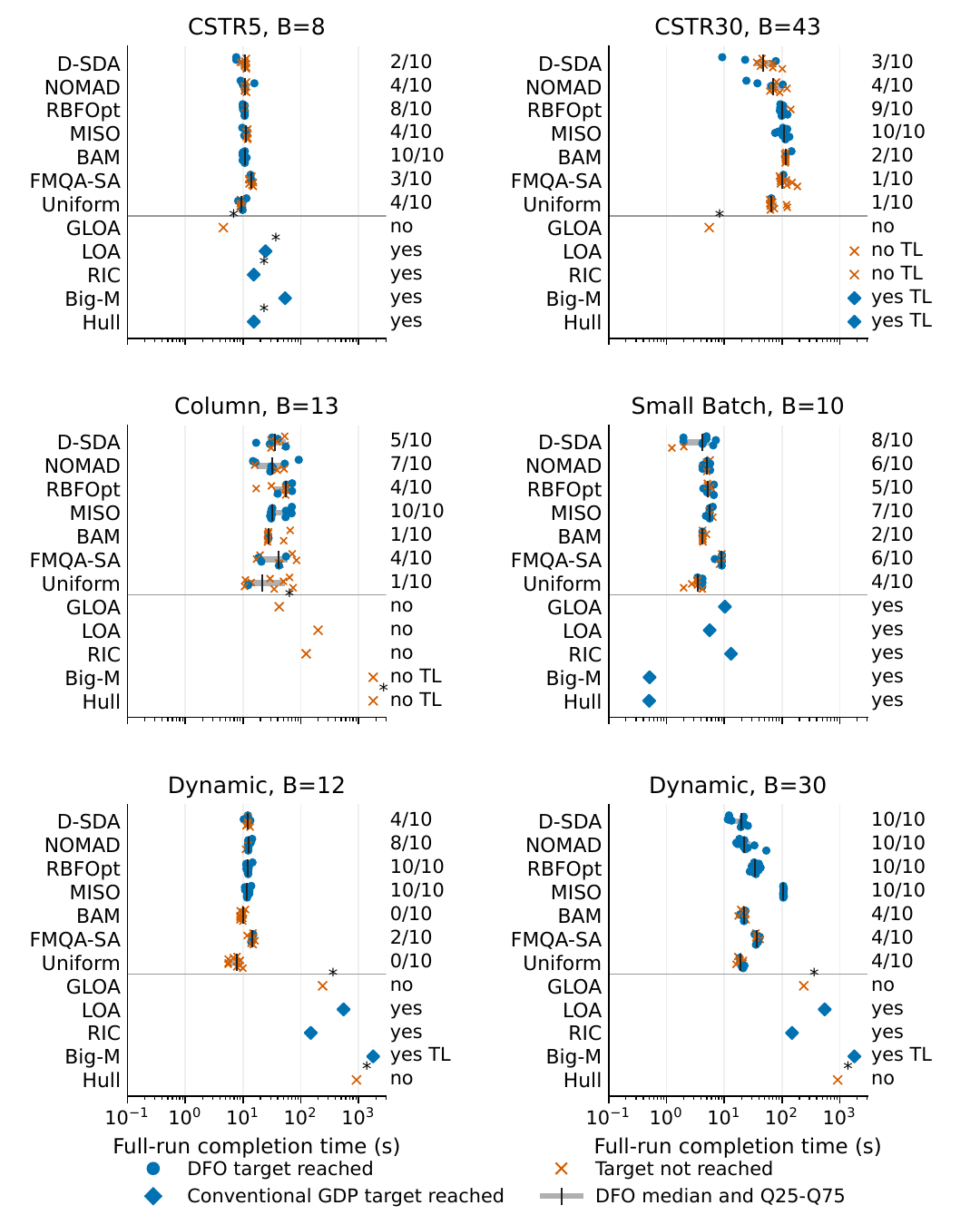}
\caption{Process GDP models: full-run completion times of the six DFO methods and uniform sampling (ten runs each, at the evaluation budget \(B\) of each panel title) and of the single run of each conventional GDP method on six conditions: GLOA (global logic-based outer approximation), LOA (logic-based outer approximation), and RIC (relaxation with integer cuts) of \texttt{GDPopt} and the Big-\(M\) and Hull reformulations solved by \texttt{SCIP}.
DFO times add the search time to the measured cost of each evaluated configuration (Equation~\eqref{eqn:reconstructed-cost}); conventional GDP times are measured and include model preparation.
Dots mark DFO runs that reach the reference target of Equation~\eqref{eqn:target}, an objective within \(10^{-5}\max\{1,|f_{\mathrm{ref}}|\}\) of the reference objective \(f_{\mathrm{ref}}\) of Table~\ref{tab:process-cases}; diamonds mark conventional GDP runs whose feasible objective is at most \(f_{\mathrm{ref}}+10^{-5}\max\{1,|f_{\mathrm{ref}}|\}\) (Section~\ref{sec:metrics}); crosses mark the other runs; gray bars span the 25th to 75th percentiles of the ten DFO seeds, with a tick at the median.
The right-hand labels give the target hits out of ten of each DFO method and the outcome, yes or no, of each conventional run, with TL for runs that report termination at the 1,800~s time limit.
Asterisks mark conventional runs whose reported lower bound exceeds a checked feasible objective (Supplementary Section~\ref{sec:bound-audit}).
The conventional GDP methods have no evaluation budget, so both \texttt{Dynamic} panels show the same conventional runs.
Supplementary Figures~\ref{fig:joint-remaining-completion} and~\ref{fig:joint-remaining-completion-dynamic} show the other eight conditions.
Uniform: uniform sampling without replacement.}
\label{fig:process-completion}
\end{figure}

Figure~\ref{fig:process-completion} places the conventional GDP methods on the same time axis, and Supplementary Table~\ref{tab:full-gdp} lists their objectives; Supplementary Figures~\ref{fig:joint-remaining-completion} and~\ref{fig:joint-remaining-completion-dynamic} show the other eight conditions.
Hull reaches the target on seven of the nine models, RIC, LOA, and Big-\(M\) on five each, and GLOA only on \texttt{Small Batch}.
Big-\(M\) and Hull solve \texttt{Small Batch} in about 0.5~s, faster than any DFO run, and they are the only conventional methods that reach the \texttt{CSTR30} target, both stopping at the 1,800~s time limit.
No conventional method reaches the \texttt{Column} target.

All six \texttt{CSTR} GLOA runs report \texttt{optimal} termination, yet none comes within the relative target tolerance of \(10^{-5}\) (Equation~\eqref{eqn:target}); their objectives exceed the references by 2.2\% on \texttt{CSTR5} and by 32\% on \texttt{CSTR30}.
Their reported lower bounds exceed independently checked feasible objectives.
An earlier configuration that used \texttt{SCIP} as a global solver for both master problems and subproblems returned the same six objectives within a relative difference of \(3\times10^{-5}\), so these misses do not come from the local NLP solver.
Supplementary Section~\ref{sec:bound-audit} reports these and the other bound inconsistencies, which affect the interpretation of reported optimality but not the target counts.

\subsection{Simulation-based lomustine case}
\label{sec:lomustine-results}
Table~\ref{tab:lomustine} reports the \texttt{Lomustine} results.
\texttt{RBFOpt} and \texttt{MISO} reach the reference configuration in eight and seven of ten runs with ten evaluations and in all ten runs with 20, and their successful runs find it at a median of the fifth to sixth and the sixth to seventh evaluation.
\texttt{NOMAD} reaches it in six runs at every budget because, when the budget allows more, it stops after a median of 12.5 evaluations.
\texttt{D-SDA} reaches it once with ten evaluations: its unit \(L_\infty\) neighborhood in three dimensions contains up to 26 configurations, so ten evaluations do not complete the first neighborhood search.
With the separable neighborhood, which has at most six neighbors, \texttt{D-SDA} reaches the reference in two, five, and five runs at the three budgets (Supplementary Section~\ref{sec:supp-dsda-neighborhood}).
\texttt{BAM} misses every target with ten evaluations, and all ten of its runs end at the same PFR--PFR configuration, 11\% above the reference.
Of the 96 misses over all methods and budgets, 41 end at this PFR--PFR configuration and 27 at the CSTR--batch configuration, 0.8\% above the reference.

In this case, an evaluation costs minutes of simulation, so the completion time is governed by the number of evaluations.
Median search times stay below 13~s, whereas the median completion times with ten evaluations range from 1.0 to 1.5~h, compared with 5.3~h for evaluating all 36 configurations.
Search adds less than 0.7~s per evaluation in every run, so the methods differ in completion time through the evaluations they choose; \texttt{MISO}, for example, completes its ten-evaluation runs in a median of 1.13~h, about as fast as any method.

Because the order of the categorical alternatives is a convention, we repeated the 210 runs with five other orders, including the reversed order and an order by degree of back-mixing (Supplementary Section~\ref{sec:supp-lomustine}).
The order changes the results of several methods.
With 20 evaluations, \texttt{BAM} reaches the reference in 0 to 10 of the ten runs depending on the order, \texttt{D-SDA} in 3 to 8, and \texttt{MISO} in 6 to 10, whereas \texttt{RBFOpt} reaches it in 9 or 10 runs in every order.
When the objective is the capital cost of each flowsheet at its optimized operation, the all-batch flowsheet is the reference and \texttt{BAM} reaches it in every run.

\begin{table}[htbp]
\centering
\caption{\texttt{Lomustine}: target hits and median full-run completion times of the six DFO methods and uniform sampling, with ten runs per method and budget.
The objective is \(-(6\times10^{4}\,\dot m+d_{10})\), and a target hit is a run whose best objective lies within \(10^{-5}\max\{1,|f_{\mathrm{ref}}|\}\) of the reference objective \(f_{\mathrm{ref}}=-487.842\) (Equation~\eqref{eqn:target}).
The reference is the best of the 36 simulated configurations, in which a plug-flow reactor (PFR) is followed by a batch reactor and the batch crystallizer.
The completion time of a run adds its search time to the recorded wall time of the \texttt{PharmaPy} simulation-based optimization of every evaluated configuration; medians are over all ten runs, including target misses, and evaluating all 36 configurations takes 5.3~h.}
\label{tab:lomustine}
\footnotesize
\setlength{\tabcolsep}{3.5pt}
\begin{threeparttable}
\begin{tabular}{@{}rlrrrrrrr@{}}\toprule
\(B\) & Measure & \texttt{D-SDA} & \texttt{NOMAD} & \texttt{RBFOpt} & \texttt{MISO} & \texttt{BAM} & \texttt{FMQA-SA} & Uniform\\\midrule
\multirow{2}{*}{10} & Target hits (of 10) & 1 & 6 & \textbf{8} & 7 & 0 & 2 & 2 \\
 & Median completion time (h) & 1.15 & \textbf{1.05} & 1.24 & 1.13 & 1.36 & 1.47 & 1.48 \\
\midrule
\multirow{2}{*}{15} & Target hits (of 10) & 4 & 6 & 9 & \textbf{10} & 4 & 3 & 4 \\
 & Median completion time (h) & 1.68 & \textbf{1.22} & 1.70 & 1.82 & 1.90 & 2.39 & 2.11 \\
\midrule
\multirow{2}{*}{20} & Target hits (of 10) & 4 & 6 & \textbf{10} & \textbf{10} & 6 & 7 & 5 \\
 & Median completion time (h) & 1.97 & \textbf{1.22} & 2.44 & 2.66 & 2.48 & 3.09 & 2.91 \\
\bottomrule
\end{tabular}
\begin{tablenotes}[flushleft]\footnotesize
\item \(\dot m\): production rate (kg/s); \(d_{10}\): number-weighted mean crystal size (\(\mu\)m); \(f_{\mathrm{ref}}\): reference objective; \(B\): evaluation budget; Uniform: uniform sampling without replacement.
Bold: most target hits and shortest median completion time for each budget, including ties.
\end{tablenotes}
\end{threeparttable}
\end{table}
\FloatBarrier

\subsection{Catalytic distillation column}
\label{sec:catdist-results}
Table~\ref{tab:catdist} reports the \texttt{Catalytic Column} results.
No method reaches the reference with 25 evaluations; with 100 evaluations, \texttt{D-SDA} and \texttt{MISO} reach it in three of ten runs and \texttt{RBFOpt} in one.
Targets are rare because the 54,159 admissible configurations of the search box are more than 500 times the largest budget and the seven initial configurations are spread over the whole space.
The eight best designs accepted in the campaign share the ethanol feed and the first catalytic stage at stage 8, the butenes feed at the last catalytic stage, and the boil-up return at stage 20 or 21, so a run must first reach this region and then refine three stage positions within it.

Most runs that miss the reference with 100 evaluations end close to it.
Nine of them, three of \texttt{D-SDA}, five of \texttt{NOMAD}, and one of \texttt{RBFOpt}, end at \((8,8,13,17,17,21)\), whose objective is 1.29~\$/yr (0.005\%) above the reference.
This design differs from the reference only in the second catalytic stage, but moving that stage from 13 to 12 changes two gaps, \((g_3,g_4)\) from \((4,3)\) to \((3,4)\).
The separable neighborhood of \texttt{D-SDA} therefore does not contain the reference, and the descent stops one stage away from it; the unit \(L_\infty\) neighborhood contains this diagonal move, but with 728 neighbors it does not fit the budgets.
The choice of external variables and neighborhood thus decides which designs are neighbors, as the order of the alternatives did for \texttt{Lomustine}: stage positions as coordinates would make this move a single step, but most of their box violates the order of the stages.
\texttt{MISO} misses by 9.37 to 653~\$/yr, and its median gap of 34.5~\$/yr is close to that of \texttt{NOMAD}.
\texttt{BAM} evaluates, after the shared initial configurations, an initial design of its own that is the same in every seed, and its last point, \((6,7,10,13,14,18)\), lies 4.7\% above the reference; in six of the ten runs with 100 evaluations, \texttt{BAM} finds nothing better.

An evaluation solves an NLP in a median of 4.8~s, whereas the median search overhead per evaluation ranges from 2~ms for uniform sampling to 0.31~s for \texttt{MISO} and 0.65~s for \texttt{FMQA-SA}, so the completion times are governed by the number of evaluations.
With 100 evaluations, the median completion times range from 8.25~min for uniform sampling to 11.4~min for \texttt{NOMAD}, and \texttt{D-SDA} and \texttt{MISO} need 8.71 and 10.5~min.
Solving every admissible configuration of the box would take about 119~h, so these runs evaluate 0.2\% of them, and the median runs of \texttt{D-SDA}, \texttt{NOMAD}, and \texttt{MISO} end within 0.14\% of the reference.

No MINLP solver reaches the reference or comes within 15\% of it (Supplementary Table~\ref{tab:supp-catdist-minlp}).
On the published MINLP form, \texttt{SBB} stops after 10~s at \((4,4,5,6,6,9)\), next to the initial configuration and 15.4\% above the reference.
\texttt{DICOPT} ends without a feasible design on all three forms after 211 to 902~s, and \texttt{SBB} finds none on the Big-\(M\) and Hull reformulations.
\texttt{BARON} and \texttt{SCIP} declare the published form infeasible after 255 and 9~s, and \texttt{BARON} declares the Big-\(M\) and Hull reformulations infeasible after 214 and 362~s.
\texttt{SCIP} finds no feasible design on the Big-\(M\) and Hull reformulations and stops at the 3,600~s time limit with a lower bound of \(-5.7\times10^{7}\)~\$/yr.
These infeasibility claims are wrong: with the binary variables fixed at the reference design, \texttt{BARON} still declares both reformulations infeasible, although a \texttt{CONOPT} solution satisfies every constraint and bound of the Big-\(M\) reformulation within \(10^{-6}\).
When started from a feasible design, both solvers instead certify that design as the global optimum of the published form, and they do so both for the design of the pilot, at 25,271.33~\$/yr, and for the reference, at 24,895.55~\$/yr, so the first certificate is contradicted by the second design.
In every MINLP form, the equations of stages that do not exist and of stages without catalyst stay in the model, multiplied by zero or relaxed, so the solvers must bound terms such as Arrhenius rates over the whole temperature range and inverse powers of vanishing liquid flows; in the Big-\(M\) reformulation, the constraint coefficients span 35 orders of magnitude.
The GDP evaluation used by the search methods solves only the equations of the selected configuration, and every one of the 210 runs returned accepted designs.

\begin{table}[htbp]
\centering
\caption{\texttt{Catalytic Column}: target hits, gaps to the reference, and median full-run completion times of the six DFO methods and uniform sampling, with ten runs per method and budget.
A target hit is a run whose best objective, the negative annual profit, lies within \(10^{-5}\max\{1,|f_{\mathrm{ref}}|\}\) of \(f_{\mathrm{ref}}=24{,}895.55\)~\$/yr, the best known design, which places the ethanol feed and the first catalytic stage at stage 8, the other catalytic stages at stages 12 and 17, the butenes feed at stage 17, and the boil-up return at stage 21 (Equation~\eqref{eqn:target}).
The gap of a run is its best objective minus \(f_{\mathrm{ref}}\), and its completion time adds its search time to the recorded \texttt{CONOPT} solve time of every evaluated configuration; both are medians over all ten runs, including target misses.
Solving all 54,159 admissible configurations of the search box would take about 119~h.}
\label{tab:catdist}
\footnotesize
\setlength{\tabcolsep}{3.5pt}
\begin{threeparttable}
\begin{tabular}{@{}rlrrrrrrr@{}}\toprule
\(B\) & Measure & \texttt{D-SDA} & \texttt{NOMAD} & \texttt{RBFOpt} & \texttt{MISO} & \texttt{BAM} & \texttt{FMQA-SA} & Uniform\\\midrule
\multirow{3}{*}{25} & Target hits (of 10) & 0 & 0 & 0 & 0 & 0 & 0 & 0 \\
 & Median gap (\$/yr) & 1{,}720 & \textbf{843} & 1{,}970 & 1{,}140 & 1{,}170 & 1{,}900 & 1{,}400 \\
 & Median completion time (min) & 2.82 & 2.76 & 2.38 & 2.51 & 2.52 & 2.94 & \textbf{2.32} \\
\midrule
\multirow{3}{*}{50} & Target hits (of 10) & 0 & 0 & \textbf{1} & \textbf{1} & 0 & 0 & 0 \\
 & Median gap (\$/yr) & 1{,}020 & \textbf{351} & 1{,}110 & 946 & 1{,}170 & 1{,}250 & 1{,}190 \\
 & Median completion time (min) & 5.75 & 5.35 & 4.43 & 5.05 & 5.13 & 5.04 & \textbf{4.00} \\
\midrule
\multirow{3}{*}{100} & Target hits (of 10) & \textbf{3} & 0 & 1 & \textbf{3} & 0 & 0 & 0 \\
 & Median gap (\$/yr) & \textbf{1.29} & 30.4 & 417 & 34.5 & 1{,}170 & 1{,}020 & 1{,}120 \\
 & Median completion time (min) & 8.71 & 11.4 & 9.07 & 10.5 & 9.52 & 10.1 & \textbf{8.25} \\
\bottomrule
\end{tabular}
\begin{tablenotes}[flushleft]\footnotesize
\item \(B\): evaluation budget; Uniform: uniform sampling without replacement.
Bold: most target hits when any run reaches the reference, smallest median gap, and shortest median completion time for each budget, including ties.
\end{tablenotes}
\end{threeparttable}
\end{table}
\FloatBarrier

\subsection{Controlled search behavior on Grid-GDP}
\label{sec:grid-search-results}
Grid-GDP isolates search behavior on explicit two-dimensional functions whose configurations need no NLP solve.
Table~\ref{tab:grid-hits} reports the separate experiments at the five budgets.
\texttt{MISO} reaches the largest target fraction at \(B=30\), 70, and 100, while \texttt{RBFOpt} leads at \(B=50\) and \texttt{BAM} at \(B=10\).
At \(B=100\), \texttt{MISO} and \texttt{RBFOpt} reach 86.2\% and 81.4\% of the targets, \texttt{NOMAD} 60.0\%, and \texttt{D-SDA} 48.1\%, with \texttt{BAM} and \texttt{FMQA-SA} in between and uniform sampling at 21.4\%.
The \texttt{NOMAD} fraction at \(B=100\) includes 16 targets reached before a reproducible solver crash; counting those 21 crashed runs as misses would give 52.4\%.
The paired bootstrap intervals in Supplementary Section~\ref{sec:supp-additional} exclude zero for \texttt{RBFOpt} against \texttt{D-SDA} at both \(B=10\) and \(B=100\).

\begin{figure}[htbp]
\centering\includegraphics[width=\linewidth]{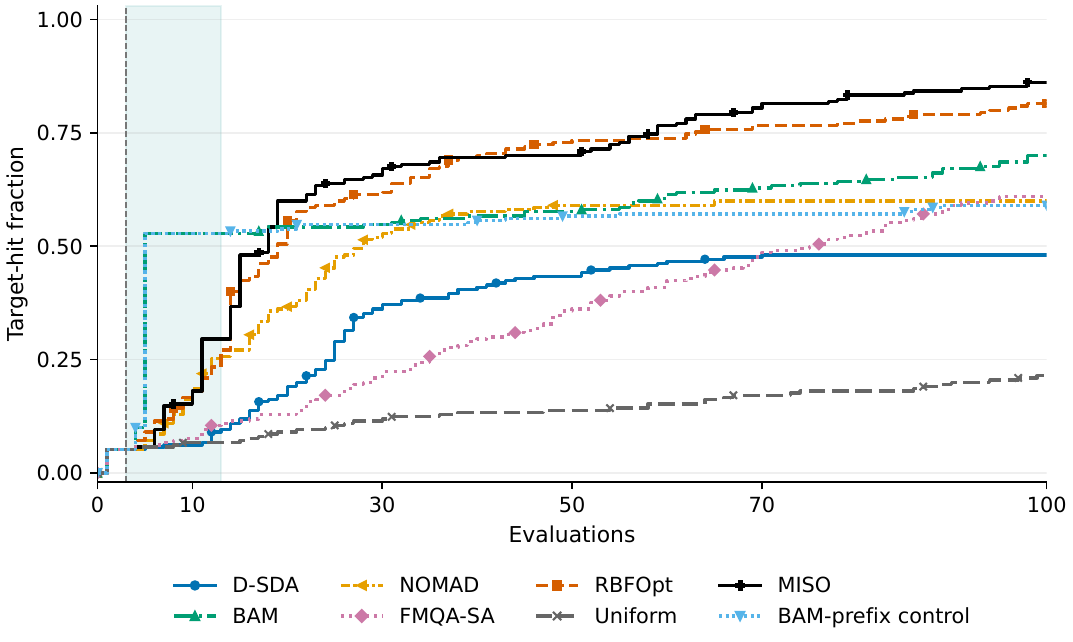}
\caption{Grid-GDP: fraction of the 210 runs (21 test functions and ten seeds) with the evaluation budget \(B=100\) that have reached the reference target of Equation~\eqref{eqn:target}, a configuration whose normalized objective is within \(10^{-5}\) of zero, within a given number of evaluations (Equation~\eqref{eqn:absolute-evaluations}).
Initial and repeated evaluations count toward the budget, and runs that miss the target remain in the denominator.
The runs with \(B=10\), 30, 50, and 70 are separate experiments; at those numbers of evaluations, this profile equals their target-hit percentages in Table~\ref{tab:grid-hits} for every method except \texttt{RBFOpt} at \(B=10\) (18.1\% against 20.0\%).
The dashed line marks the end of the three shared initial evaluations, and the shaded band the ten predetermined points of \texttt{BAM}, which the \texttt{BAM}-prefix control repeats before it samples uniformly (Section~\ref{sec:sampling-controls}).
Uniform: uniform sampling without replacement.
At 100 evaluations, \texttt{MISO} reaches the largest fraction, followed by \texttt{RBFOpt}.}
\label{fig:grid-profiles}
\end{figure}

At \(B=10\), \texttt{BAM} and its prefix control each reach 52.9\% because nine of the 21 functions have their target at the domain center, which \texttt{BAM} evaluates as the fifth point of its Sobol design.
By \(B=100\), the model-based search of \texttt{BAM} adds 23 hits over the prefix control (Figure~\ref{fig:grid-profiles}).
\texttt{D-SDA} levels off at 48.1\% from \(B=70\) onward, once the descents from its three initial configurations end, and its median number of evaluations stays at 37 for every budget from \(B=50\) onward.
Excluding Xin-She Yang~3, whose target set covers 870 of the 961 lattice points, and Bohachevsky~2, which nearly duplicates Sum of Squares on the lattice, preserves this ordering.
\begin{table}[htbp]
\centering
\caption{Grid-GDP: target hits of the six DFO methods, uniform sampling, and the \texttt{BAM}-prefix control at five evaluation budgets, as percentages of 210 runs (21 test functions and ten seeds).
A target hit is a run that evaluates a configuration whose normalized objective is within \(10^{-5}\) of zero, the normalized minimum over the lattice (Equation~\eqref{eqn:target}).
Each budget is a separate experiment and counts the three shared initial evaluations and repeated configurations.
At \(B=10\), all 111 \texttt{BAM} hits precede its surrogate-guided selection, 11 at the shared initial configurations and 100 at two of its predetermined points, and the \texttt{NOMAD} percentage at \(B=100\) includes 16 targets reached before a solver crash that recurred on every attempt.}
\label{tab:grid-hits}
\footnotesize
\setlength{\tabcolsep}{3.5pt}
\begin{threeparttable}
\begin{tabular}{@{}rrrrrrrrr@{}}\toprule
 & \multicolumn{8}{c}{Target hits (\% of 210 runs)}\\\cmidrule(l){2-9}
\(B\) & \texttt{D-SDA} & \texttt{NOMAD} & \texttt{RBFOpt} & \texttt{MISO} & \texttt{BAM} & \texttt{FMQA-SA} & Uniform & \makecell{\texttt{BAM}-prefix\\control}\\\midrule
10 & 6.2 & 18.6 & 20.0 & 18.1 & \textbf{52.9} & 7.6 & 6.7 & \textbf{52.9} \\
30 & 37.1 & 52.9 & 61.9 & \textbf{67.1} & 54.8 & 22.4 & 11.9 & 54.8 \\
50 & 43.3 & 59.0 & \textbf{73.3} & 70.0 & 57.6 & 36.2 & 13.8 & 56.7 \\
70 & 48.1 & 60.0 & 76.7 & \textbf{81.4} & 62.9 & 48.6 & 17.1 & 57.1 \\
100 & 48.1 & 60.0 & 81.4 & \textbf{86.2} & 70.0 & 61.4 & 21.4 & 59.0 \\
\bottomrule
\end{tabular}
\begin{tablenotes}[flushleft]\footnotesize
\item \(B\): evaluation budget; Uniform: uniform sampling without replacement; \texttt{BAM}-prefix control: the ten predetermined points of \texttt{BAM} after the three shared initial configurations, then uniform sampling (Section~\ref{sec:sampling-controls}).
Bold: largest percentage in each row, including ties.
\end{tablenotes}
\end{threeparttable}
\end{table}

Because the explicit expression is available, the direct Big-\(M\) and Hull reformulations solved with \texttt{SCIP} reach all 21 Grid targets in a median of 0.66 and 0.53~s per function, less than the median completion time of every DFO method (Supplementary Table~\ref{tab:grid-cost}).
The \texttt{GDPopt} methods reach 14 or 15 of the 21 targets, and in diagnostic runs their preprocessing rejected the feasible reference configurations of the three functions that all of them miss.
Grid-GDP therefore tests search behavior; for these explicit functions, a direct reformulation is the method of choice.

\FloatBarrier

\subsection{Discussion}
\label{sec:cross-benchmark}
No single method reaches the most targets in every case.
The ranking depends on the size of the configuration space, on where the target lies relative to the initial configurations and to the configurations without accepted objectives, and on how each method uses its evaluation budget.

On \texttt{CSTR}, the recycle/feed ordering makes roughly half of each configuration box logically inadmissible, and the target lies at the corner \((N,N)\), on the edge of the admissible region.
On \texttt{CSTR10} through \texttt{CSTR30}, 31 of the 36 \texttt{D-SDA} misses end at the configuration \((5,1)\), whose volume of 3.13020 is the same for every \(N\).
Because the reference volume decreases with \(N\), this local minimum lies 8.3\% above the reference for \texttt{CSTR10} and 13.0\% above it for \texttt{CSTR30}, which explains why the \texttt{D-SDA} gaps grow with size.
Restarted descents mostly end at \((5,1)\) or at the target, so the misses reflect the large basin of this local minimum, and additional budget does not remove them.
\texttt{NOMAD}, which polls configurations around its incumbent, ends 32 of its 33 misses on these sizes at \((5,1)\) as well.
\texttt{RBFOpt} and \texttt{MISO} can propose configurations far from their incumbents and reach the corner in nine and ten of the ten \texttt{CSTR30} runs.
A global surrogate does not guarantee success on small budgets, however: \texttt{MISO} reaches 13 of the 30 targets on \texttt{CSTR5} through \texttt{CSTR15}, with 8 to 22 evaluations, compared with 21 for \texttt{RBFOpt}.
On \texttt{Column}, four of the five successful \texttt{D-SDA} runs reach the target in one or two improving moves, and the fifth starts at it.
With only 29 and 12 accepted configurations and budgets of 13 and 10 evaluations, \texttt{Column} and \texttt{Small Batch} give a surrogate model few observations to learn from, and local descent remains competitive with \texttt{RBFOpt}.
\texttt{MISO} nevertheless reaches all ten \texttt{Column} targets; in the nine runs whose initial design does not contain the target, it finds it between the 7th and 11th of its 13 evaluations.
On Grid-GDP, the early lead of \texttt{BAM} reflects the symmetry of the test functions with respect to its initial design, and \texttt{FMQA-SA} improves on uniform sampling there, but in the small process domains, where more than half of its proposals come from its uniform fallback, it reaches only five more targets than uniform sampling (42 against 37 of 140).
From \(B=30\) onward, the two radial-basis-function methods, \texttt{MISO} and \texttt{RBFOpt}, reach the most Grid targets.

The methods also use additional evaluations differently.
\texttt{D-SDA} terminates when its descents reach configurations without an improving neighbor, and \texttt{NOMAD} when its mesh reaches the integer resolution, so additional budget helps them only until then.
On Grid, \texttt{NOMAD} stops after a median of 31.5 evaluations from \(B=50\) onward, and its target fraction stays at 60.0\% from \(B=70\) to \(B=100\).
The surrogate-based methods and the sampling controls have no comparable termination test and use the full budget in every Grid run, so their Grid target-hit rates keep increasing between \(B=70\) and \(B=100\), by 4.8 percentage points for \texttt{RBFOpt} and \texttt{MISO}, 7.1 for \texttt{BAM}, and 12.9 for \texttt{FMQA-SA}, whereas the \texttt{D-SDA} rate stays at 48.1\%.
Once a method has found the target, however, the remaining evaluations add cost without improving the result: \texttt{RBFOpt} and \texttt{MISO} reach all ten \texttt{Dynamic} targets at \(B=12\), yet their median completion times grow from 12.29 and 11.80~s at \(B=12\) to 65.14 and 195.36~s at \(B=50\), whereas \texttt{NOMAD}, which stops early, needs 21.99~s at \(B=50\).
Additional search time also does not guarantee more targets: on Grid at \(B=100\), \texttt{FMQA-SA} spends a median of 62.04~s on search and reaches fewer targets than \texttt{BAM}, which spends 0.43~s.
For \texttt{MISO}, it does: with a median of 32.70~s, about four times the search time of \texttt{RBFOpt}, it reaches 181 targets against 171 (Supplementary Table~\ref{tab:grid-cost}).
How much search time is worthwhile depends on the cost of an evaluation.
\texttt{FMQA-SA} spends about 0.5~s on search per evaluation, comparable to the evaluations in these models, while \texttt{RBFOpt} spends 17 to 98~ms in median but more as observations accumulate.
\texttt{MISO} spends 52~ms to 1.86~s per evaluation, and on \texttt{Dynamic} its search per evaluation exceeds the median evaluation cost of 0.84~s; \texttt{NOMAD} spends 3 to 97~ms (Supplementary Table~\ref{tab:supp-overhead}).
When each evaluation requires a longer optimization or a process simulation, this overhead becomes small, and the number of subproblem solves becomes the main cost.

Because the initial configurations change with the seed and the surrogate-based methods are randomized, repeated runs can return different solutions, and the spread of these solutions differs between methods.
\texttt{RBFOpt} returns at most three distinct final objectives across the ten seeds of each \texttt{CSTR}, \texttt{Small Batch}, and \texttt{Dynamic} condition; on \texttt{CSTR10} through \texttt{CSTR30}, its runs either reach the target or stop 8.1 to 12.5\% above it.
\texttt{MISO} and \texttt{NOMAD} are similarly consistent, with at most five and three distinct final objectives per condition, and the \texttt{MISO} misses on \texttt{CSTR10} and \texttt{CSTR15} lie 6.0 to 10.6\% above the reference.
\texttt{BAM} reaches the exact target in only 18 of the 60 \texttt{CSTR} runs, but all of its final objectives lie within 0.7\% of the reference.
\texttt{FMQA-SA} and uniform sampling spread more widely, with up to ten distinct final objectives in one condition and gaps of up to 17\%.
With a relative tolerance of 1\%, \texttt{BAM} would succeed in all 60 \texttt{CSTR} runs, compared with 45 for \texttt{RBFOpt}, 43 for \texttt{MISO}, 21 for \texttt{NOMAD}, and 19 for \texttt{D-SDA} (Supplementary Table~\ref{tab:supp-spread}), so the strict target separates exact attainment from near-optimal solutions that may suffice in a design study.

The comparison with enumeration and conventional GDP methods shows how the cost of each approach grows with the configuration space.
For small models such as \texttt{Small Batch} and the Grid functions, enumeration or a direct reformulation finds the optimum in seconds or less.
As the configuration space grows, logic-based methods that solve one subproblem per admissible configuration and direct reformulations that reach the time limit become expensive, while surrogate search keeps its number of subproblem solves near the budget, as on \texttt{CSTR30}.
\texttt{Dynamic} lies between these cases: RIC reaches its target in 149.39~s after solving all 37 admissible subproblems, and \texttt{RBFOpt} and \texttt{MISO} reach it in every run at \(B=12\) in medians of 12.29 and 11.80~s.

The \texttt{Lomustine} case shows that the configuration search carries over when a simulator defines each subproblem and conventional GDP methods do not apply.
Its ranking agrees with the \texttt{CSTR} and Grid results: the two radial-basis-function methods reach the reference most often, \texttt{NOMAD} stops early, and the size of the \texttt{D-SDA} neighborhood limits local descent on small budgets.
Because each evaluation costs minutes, the search overhead that slows \texttt{MISO} on \texttt{Dynamic} is negligible here.
Unordered alternatives can still be encoded as external variables, but their order is then part of the problem the search sees: with 20 evaluations, only \texttt{RBFOpt} reached the reference in at least nine of ten runs in every order.

The neighborhood of \texttt{D-SDA} trades the number of evaluations per step against the moves it can make.
Repeating every stored-response \texttt{D-SDA} run with the separable neighborhood, which changes one external variable at a time, shows no overall winner (Supplementary Section~\ref{sec:supp-dsda-neighborhood}).
The separable neighborhood reaches more targets when the budget barely covers the \(L_\infty\) neighborhood, as on Grid at \(B=10\) (22 against 13 of 210 runs) and on \texttt{Lomustine} (99 against 90 of 180 runs over the six orders).
It reaches fewer when a descent must move diagonally: on \texttt{CSTR}, 7 against 16 of 60 runs, because 18 of its misses stop on the edge of the admissible region that leads to the target.
The number of \(L_\infty\) neighbors grows as \(3^d-1\) with the number \(d\) of external variables, which is why the \texttt{Catalytic Column}, with six external variables, uses the separable neighborhood.

The \texttt{Catalytic Column} tests the search where enumeration is impractical and every evaluation is a live NLP solve.
With 100 evaluations, 0.2\% of the configurations in the search box, the median runs of \texttt{D-SDA}, \texttt{MISO}, and \texttt{NOMAD} end within 0.14\% of the best known design, and local descent is competitive with the surrogate methods here because the good designs form one region: six of the ten \texttt{D-SDA} descents end at the reference or one stage from it.
Within that region, the external variables and the neighborhood decide which designs are neighbors: in the gap coordinates, the runner-up and the reference are two coordinates apart although they differ in one stage, so the separable neighborhood stops local descent 1.29~\$/yr short of the reference, as the order of categorical alternatives changed the results of \texttt{Lomustine}.
The MINLP forms of the same model keep the equations of every possible stage and catalyst position, and the solvers fail in different ways: \texttt{SBB} stops at its limit on function evaluations outside their domains, \texttt{DICOPT} when its NLP subproblems fail or are all infeasible, \texttt{SCIP} at the time limit on the Big-\(M\) and Hull forms, and \texttt{BARON} and \texttt{SCIP} with infeasibility claims or certificates that the feasible designs contradict.
A plausible common cause is that terms harmless in a solved configuration, such as reaction rates on stages without catalyst and hydraulics with vanishing flows, make the relaxations poorly scaled, with constraint coefficients that span 35 orders of magnitude in the Big-\(M\) form; the size of the full MINLP, which contains all 100,947 configurations, may contribute as well.
Reported bounds were also unreliable in the process models, where GLOA terminated as optimal on every \texttt{CSTR} network while missing the reference by 2.2 to 32\%.
The GDP model avoids these terms because each evaluation contains only the equations of the selected configuration, so the external-variable representation lets a black-box method optimize a model that its direct MINLP reformulations cannot solve reliably.

Except for the \texttt{Catalytic Column}, these conclusions rest on configuration spaces of 25 to 961 configurations, on NLP responses obtained with one initialization per model, and on completion times that combine recorded search paths with one cost measurement per configuration.
The \texttt{Catalytic Column} reference is the best design found and has no reliable optimality certificate, because the global solvers that certify it also certify a worse design, and its responses also come from one initialization.
Because the \texttt{CSTR} budgets grow with \(N\), the \texttt{CSTR} series does not isolate the effect of size.
With ten seeds per condition, the per-condition comparisons that involve \texttt{NOMAD} or \texttt{MISO} are not significant after a Holm correction over the 70 exploratory tests, and the per-condition tests of \texttt{RBFOpt} against \texttt{D-SDA} are not adjusted for multiple comparisons, so rankings within a single condition are descriptive.
The \texttt{NOMAD} and \texttt{MISO} search times were measured in a later campaign on a shared server, so small timing differences between them and the other methods should not be overinterpreted.

\section{Conclusions}
This study modeled process design alternatives as GDP models and searched their configurations with derivative-free methods, using the external-variable reformulation as the interface between the two, and compared this pairing with local descent, conventional GDP methods, and full enumeration.
Five lessons follow from the experiments.

First, separating the model from the search extends GDP to models that conventional GDP methods cannot handle.
A configuration search needs only the value of each configuration, so the same GDP model can be solved when its subproblems are simulations, as for the \texttt{Lomustine} flowsheet, or when its MINLP reformulations fail, as for the \texttt{Catalytic Column}, whose evaluations contain only the equations of the selected design.

Second, surrogate search spends subproblem solves economically.
RIC and LOA solved one subproblem for each admissible configuration on most process models, whereas the two radial-basis-function methods reached the best configuration of the largest reactor network with a small fraction of those solves.

Third, the landscape and its representation decide whether a global search pays off.
Local descent can be competitive when the good configurations form one region, as it was in reaching designs close to the reference of the \texttt{Catalytic Column}, and surrogate methods that interpolate all observations are the safer choice when large basins of local minima lie between the initial configurations and the optimum.
The external variables, the order of categorical alternatives, and the neighborhood determine which configurations are close, so they shape the landscape that every method sees.

Fourth, exact target attainment is an incomplete measure for stochastic methods.
A method that rarely matches the reference can still end consistently close to it, so the spread of final objectives across runs belongs next to the target counts when methods are compared for design studies.

Fifth, the value of search effort depends on the cost of an evaluation.
Additional search effort per evaluation can improve target attainment, but its overhead can dominate the run time when evaluations are cheap, and methods that stop at a local criterion leave part of the budget unused; as evaluations become more expensive, the number of subproblem solves becomes the cost that matters.

For future GDP instances, these results suggest enumerating the configurations or solving a direct Big-\(M\) or Hull reformulation when the model is small, as for Grid-GDP and \texttt{Small Batch}.
When the configuration space is too large to enumerate and contains large basins of local minima, as for the larger \texttt{CSTR} networks, surrogate-based configuration search with radial basis functions is preferable; in these experiments, \texttt{MISO} reached slightly more targets when its search overhead was small relative to the evaluation cost, and \texttt{RBFOpt} was the less sensitive to the order of categorical alternatives.
For local descent, the unit \(L_\infty\) neighborhood suits two external variables and budgets that cover its neighbors, and the separable neighborhood suits more external variables or smaller budgets.
When the subproblems are simulations, or when the MINLP reformulations keep equations that are ill-defined for absent units, the GDP model is best solved by configuration search, in which local descent can be competitive when the good designs form one region and surrogate search is the safer choice when the landscape has several basins.
Because the wall-time advantage of searching fewer configurations grows with the cost of each evaluation, testing these recommendations on live simulation-based subproblems, in which every evaluation runs the inner optimization, is the next step.

\section*{CRediT authorship contribution statement}
\textbf{Albert Joon Lee}: Software, Investigation, Formal analysis, Validation, Visualization, Writing -- original draft, Writing -- review and editing.
\textbf{Woosik Kim}: Methodology, Software (factorization-machine surrogate optimization with simulated annealing), Validation, Visualization.
\textbf{David E. Bernal Neira}: Conceptualization, Methodology, Software, Validation, Formal analysis, Investigation, Resources, Data curation, Visualization, Supervision, Project administration, Funding acquisition, Writing -- original draft, Writing -- review and editing.

\section*{Declaration of competing interests}
The authors declare that they have no competing financial interests or personal relationships that could have influenced the work reported in this paper.

\section*{Funding}
This work was supported by the Davidson School of Chemical Engineering at Purdue University.
The funder had no role in study design, data collection, analysis or interpretation, preparation of the manuscript, or the decision to submit it for publication.

\section*{Acknowledgements}
The authors thank the Davidson School of Chemical Engineering at Purdue University for providing a supportive research environment.

\section*{Data and code availability}
Sanitized code and data supporting this study are available to reviewers on request from the corresponding author.
Public release is planned through the \href{https://github.com/SECQUOIA/Black-Box-for-GDP/tree/ea65c5339486b576008df68645ea68d22f7c7213}{Black-Box-for-GDP repository}~\citep{rbfoptTestArchive}.
Supplementary Section~\ref{sec:supp-data} identifies the archived inputs, including the CSTR updates, Column logical-domain runs, conventional GDP comparisons, \texttt{D-SDA} reruns, the \texttt{Lomustine} runs, and the \texttt{Catalytic Column} runs and MINLP comparisons.

\section*{Declaration of generative AI and AI-assisted technologies in the manuscript preparation process}
The authors used \texttt{ChatGPT} and \texttt{Codex} (OpenAI), including the GPT-6 Astra model, to assist with drafting, organizing, and revising the manuscript and developing analysis and plotting code.
The authors also used \texttt{Claude Code} (Anthropic), with the Claude Opus 5.5 model, to review the manuscript, revise the Discussion and Conclusions, and develop analysis, plotting, and experiment code.
The resulting text and code were checked and revised by the authors, who remain responsible for the article's content.

\bibliographystyle{elsarticle-harv}
\bibliography{reference}
\ifdefined\cacesubmission
\else
\def\manuscriptwithsupplement{1}
\ifdefined\papertitle\else
\newcommand{\papertitle}{Generalized Disjunctive Programming with Derivative-Free Configuration Search}
\newcommand{\paperfirstauthor}{Albert Joon Lee}
\newcommand{\papersecondauthor}{Woosik Kim}
\newcommand{\papercorrespondingauthor}{David E. Bernal Neira}
\newcommand{\papercorrespondingemail}{dbernaln@purdue.edu}
\newcommand{\paperauthors}{\paperfirstauthor, \papersecondauthor, and \papercorrespondingauthor}
\fi
\ifdefined\manuscriptwithsupplement
\clearpage
\begingroup
\normalsize
\setcounter{section}{0}
\setcounter{table}{0}
\setcounter{figure}{0}
\setcounter{equation}{0}
\renewcommand{\thesection}{S\arabic{section}}
\renewcommand{\thetable}{S\arabic{table}}
\renewcommand{\thefigure}{S\arabic{figure}}
\renewcommand{\theequation}{S\arabic{equation}}
\renewcommand{\theHsection}{supp.\arabic{section}}
\renewcommand{\theHtable}{supp.\arabic{table}}
\renewcommand{\theHfigure}{supp.\arabic{figure}}
\renewcommand{\theHequation}{supp.\arabic{equation}}
\setlength{\emergencystretch}{2em}
\clubpenalty=10000
\widowpenalty=10000
\begin{center}
{\LARGE Supplementary Material:\\Benchmark Definitions, Execution Protocol, and Uncertainty\par}
\medskip
{\large \papertitle\par}
\medskip
\paperauthors
\end{center}
\else
\documentclass{article}
\usepackage[T1]{fontenc}
\usepackage[utf8]{inputenc}
\usepackage[letterpaper,margin=1in]{geometry}
\usepackage{amsmath,amssymb,booktabs,placeins,xurl}
\usepackage{longtable,pdflscape,graphicx,makecell,multirow}
\usepackage{threeparttable}
\renewcommand{\thefigure}{S\arabic{figure}}
\usepackage[authoryear,round,sort]{natbib}
\usepackage[hidelinks]{hyperref}
\usepackage{xcolor}
\providecommand{\david}[1]{\textcolor{red}{D: #1}}
\providecommand{\albert}[1]{\textcolor{blue}{A: #1}}
\providecommand{\claude}[1]{\textcolor{teal}{C: #1}}
\renewcommand{\thetable}{S\arabic{table}}
\renewcommand{\theequation}{S\arabic{equation}}
\renewcommand{\thesection}{S\arabic{section}}
\setlength{\emergencystretch}{2em}
\clubpenalty=10000
\widowpenalty=10000
\title{Supplementary Material for\\\papertitle\\[1ex]\large Benchmark Definitions, Execution Protocol, and Uncertainty}
\author{\paperauthors}
\date{}
\setlength{\emergencystretch}{2em}
\begin{document}
\maketitle
\vspace{-2em}
\fi
The Grid-GDP functions are adapted from the Quantum Hamiltonian Descent collection; each instance is defined by the expression and domain implemented there, which can differ from other functions known by the same name~\citep{qhdBenchmarks,rbfoptTestArchive}.

\section{Coordinate domains and normalization}\label{sec:supp-domains}
The coordinate mapping and objective normalization follow Section~\ifdefined\manuscriptwithsupplement\ref{sec:grid-gdp-formulation}\else3.4\fi{} of the main paper.
Both coordinate endpoints are included, and \(L\) and \(U\) are the minimum and maximum over the 961 lattice points.
All trigonometric arguments are in radians.
Full-precision normalization constants and target coordinates are archived in \nolinkurl{grid_inventory.json} (Section~\ref{sec:supp-data}).

\begin{table}[htbp]
\centering
\caption{Grid-GDP: coordinate domains and target-set sizes of the 21 test functions of Section~\ref{sec:supp-objectives}.
Each coordinate range is sampled at 31 equally spaced points that include both endpoints, which gives a \(31\times31\) lattice of 961 configurations and 62 Boolean indicators.
The target set of a function contains the lattice points whose normalized objective \(F=(f-L)/(U-L)\) is at most \(10^{-5}\), and a run reaches the reference target of Equation~\ifdefined\manuscriptwithsupplement\eqref{eqn:target}\else(8)\fi{} of the main paper when it evaluates one of these points.}
\label{tab:grid-domains}
\footnotesize
\begin{threeparttable}
\begin{tabular}{@{}lllr@{}}\toprule
Function & \(x\) domain & \(y\) domain & \(n_{\rm target}\)\\\midrule
Ackley & \([-5,5]\) & \([-5,5]\) & 1\\
Alpine 1 & \([0,10]\) & \([0,10]\) & 1\\
Alpine 2 & \([0,1000]\) & \([0,1000]\) & 1\\
Bohachevsky 2 & \([-100,100]\) & \([-100,100]\) & 1\\
Three-hump camel & \([-5,5]\) & \([-5,5]\) & 1\\
Csendes (squared variant) & \([-1,1]\) & \([-1,1]\) & 36\\
Deflected corrugated spring & \([0,10]\) & \([0,10]\) & 1\\
Drop-Wave & \([-5.12,5.12]\) & \([-5.12,5.12]\) & 1\\
Easom & \([-10,10]\) & \([-10,10]\) & 1\\
Griewank & \([-600,600]\) & \([-600,600]\) & 1\\
H\"older table & \([-10,10]\) & \([-10,10]\) & 4\\
Hosaki & \([0,5]\) & \([0,6]\) & 1\\
L\'evy & \([-10,10]\) & \([-10,10]\) & 2\\
L\'evy 13 & \([-10,10]\) & \([-10,10]\) & 4\\
Michalewicz & \([0,\pi]\) & \([0,\pi]\) & 1\\
Rastrigin & \([-5.12,5.12]\) & \([-5.12,5.12]\) & 1\\
Rosenbrock & \([-2,2]\) & \([-1,3]\) & 2\\
Shubert (variant) & \([0,10]\) & \([0,10]\) & 1\\
Styblinski--Tang & \([-5,5]\) & \([-5,5]\) & 1\\
Sum of squares & \([-10,10]\) & \([-10,10]\) & 1\\
Xin-She Yang 3 (variant) & \([-20,20]\) & \([-20,20]\) & 870\\\bottomrule
\end{tabular}
\begin{tablenotes}[flushleft]\footnotesize
\item \(x\) domain, \(y\) domain: coordinate ranges; \(n_{\rm target}\): number of lattice points in the target set; \(F\): normalized objective; \(L\), \(U\): minimum and maximum of the objective \(f\) over the lattice.
\end{tablenotes}
\end{threeparttable}
\end{table}

\FloatBarrier 
\section{Grid-GDP objective functions}\label{sec:supp-objectives}
The implemented functions below use the domains in Table~\ref{tab:grid-domains}~\citep{rbfoptTestArchive}.

\subsection*{1. Ackley}
\begin{equation}\label{eq:supp-ackley}
\begin{aligned}
f(x,y)={}&-20\exp\!\left(-0.2\sqrt{\frac{x^2+y^2}{2}}\right)\\
&-\exp\!\left(\frac{\cos(2\pi x)+\cos(2\pi y)}{2}\right)+20+e.
\end{aligned}
\end{equation}

\subsection*{2. Alpine 1}
\begin{equation}\label{eq:supp-alpine1}
f(x,y)=|x\sin x+0.1x|+|y\sin y+0.1y|.
\end{equation}

\subsection*{3. Alpine 2}
\begin{equation}\label{eq:supp-alpine2}
f(x,y)=-\sqrt{xy}\,\sin x\sin y.
\end{equation}

\subsection*{4. Bohachevsky 2}
\begin{equation}\label{eq:supp-bohachevsky2}
f(x,y)=x^2+2y^2-0.3\cos(3\pi x)\cos(4\pi y)+0.3.
\end{equation}

\subsection*{5. Three-hump camel}
\begin{equation}\label{eq:supp-camel3}
f(x,y)=2x^2-1.05x^4+\frac{x^6}{6}+xy+y^2.
\end{equation}

\subsection*{6. Csendes (squared variant)}
\begin{equation}\label{eq:supp-csendes}
f(x,y)=c(x)+c(y),\qquad
c(t)=\begin{cases}
t^6\bigl[2+\sin(1/t)\bigr]^2,&t\ne0,\\
0,&t=0.
\end{cases}
\end{equation}
The GDP uses branch-fixed auxiliaries \(r_t=1/t\) for \(t\ne0\) and \(r_t=0\) otherwise, avoiding division by zero without an epsilon shift.

\subsection*{7. Deflected corrugated spring}
With \(r=\sqrt{(x-5)^2+(y-5)^2}\),
\begin{equation}\label{eq:supp-spring}
f(x,y)=0.1r^2-\cos(5r).
\end{equation}

\subsection*{8. Drop-Wave}
\begin{equation}\label{eq:supp-dropwave}
f(x,y)=-\frac{1+\cos\!\left(12\sqrt{x^2+y^2}\right)}{0.5(x^2+y^2)+2}.
\end{equation}

\subsection*{9. Easom}
\begin{equation}\label{eq:supp-easom}
f(x,y)=-\cos x\cos y\,\exp\!\left(-[(x-\pi)^2+(y-\pi)^2]\right).
\end{equation}

\subsection*{10. Griewank}
\begin{equation}\label{eq:supp-griewank}
f(x,y)=\frac{x^2+y^2}{4000}-\cos x\cos\!\left(\frac{y}{\sqrt{2}}\right)+1.
\end{equation}

\subsection*{11. H\"older table}
\begin{equation}\label{eq:supp-holder}
f(x,y)=-\left|\sin x\cos y\,
\exp\!\left(\left|1-\frac{\sqrt{x^2+y^2}}{\pi}\right|\right)\right|.
\end{equation}

\subsection*{12. Hosaki}
\begin{equation}\label{eq:supp-hosaki}
f(x,y)=\left(1-8x+7x^2-\frac{7}{3}x^3+\frac14x^4\right)y^2\exp(-y).
\end{equation}

\FloatBarrier 
\subsection*{13. L\'evy}
Define \(w_x=1+(x-1)/4\) and \(w_y=1+(y-1)/4\).
Then
\begin{equation}\label{eq:supp-levy}
\begin{aligned}
f(x,y)={}&\sin^2(\pi w_x)\\
&+(w_x-1)^2\bigl[1+10\sin^2(\pi w_x+1)\bigr]\\
&+(w_y-1)^2\bigl[1+\sin^2(2\pi w_y)\bigr].
\end{aligned}
\end{equation}

\subsection*{14. L\'evy 13}
\begin{equation}\label{eq:supp-levy13}
\begin{aligned}
f(x,y)={}&\sin^2(3\pi x)+(x-1)^2\bigl[1+\sin^2(3\pi y)\bigr]\\
&+(y-1)^2\bigl[1+\sin^2(2\pi y)\bigr].
\end{aligned}
\end{equation}

\subsection*{15. Michalewicz}
\begin{equation}\label{eq:supp-michalewicz}
f(x,y)=-\sin x\,\sin^{20}\!\left(\frac{x^2}{\pi}\right)
-\sin y\,\sin^{20}\!\left(\frac{2y^2}{\pi}\right).
\end{equation}

\subsection*{16. Rastrigin}
\begin{equation}\label{eq:supp-rastrigin}
f(x,y)=x^2-10\cos(2\pi x)+y^2-10\cos(2\pi y)+20.
\end{equation}

\subsection*{17. Rosenbrock}
\begin{equation}\label{eq:supp-rosenbrock}
f(x,y)=(1-x)^2+100(y-x^2)^2.
\end{equation}

\subsection*{18. Shubert (variant)}
\begin{equation}\label{eq:supp-shubert}
\begin{aligned}
f(x,y)={}&\bigl[\cos(2x+1)+2\cos(3x+2)+3\cos(4x+3)\bigr]\\
&\times\bigl[\cos(2y+1)+\cos(y+2)\bigr].
\end{aligned}
\end{equation}

\FloatBarrier 
\subsection*{19. Styblinski--Tang}
\begin{equation}\label{eq:supp-styblinski}
f(x,y)=\frac{1}{78}\,\frac12
\left(x^4-16x^2+5x+y^4-16y^2+5y\right).
\end{equation}

\subsection*{20. Sum of squares}
\begin{equation}\label{eq:supp-squares}
f(x,y)=x^2+2y^2.
\end{equation}

\subsection*{21. Xin-She Yang 3 (variant)}
\begin{equation}\label{eq:supp-xin}
f(x,y)=\exp\!\left(-\left[\left(\frac{x}{15}\right)^6+y^6\right]\right)
+2\exp\!\left(-(x^2+y^2)\right)\cos^2x\cos^2y.
\end{equation}

\FloatBarrier 
\section{Executable protocol and implementation settings}\label{sec:supp-protocol}

\texttt{Python} and \texttt{NumPy} use the job seed, and \texttt{RBFOpt} also receives \nolinkurl{rand_seed}.
\texttt{BAM} exposes no separate native random-seed option: its process wrapper seeds \texttt{Python}/\texttt{NumPy}, while the Grid seed identifies the common initial design.
\texttt{FMQA-SA} uses three \texttt{NumPy} \nolinkurl{SeedSequence(seed)} streams for FM initialization, annealing, and uniform unseen-point fallback; each training or annealing invocation draws a seed from its stream.
\texttt{PyTorch} uses deterministic CPU operations and one thread.

Grid \texttt{D-SDA} enumerates \(L_\infty\) neighbors lexicographically from \((-1,0,1)\), excluding the zero direction and out-of-bounds points; ties select the first point.

The \texttt{RBFOpt} options are \nolinkurl{rbf=auto}, \nolinkurl{global_search_method=genetic}, \nolinkurl{do_init_strategy=False}, \nolinkurl{ga_base_population_size=400}, and \nolinkurl{ga_num_generations=20}.

\texttt{BAM} receives \nolinkurl{x0=None} and \nolinkurl{max_evals}=B-3; its local models use \texttt{ALAMO}/\texttt{BARON}, whose internal versions are not reported by the package.

The \texttt{FMQA-SA} time limit is checked between operations and through training/sampling callbacks.
Each FM uses identity activation, zero bias, and normal initial weights with mean zero and standard deviation \(0.01\).
Training uses Adam~\citep{kingma2015adam}, 1,000 full-batch epochs, mean squared error, and learning rate \(0.01\).
Coordinate \(\xi_j\) uses \(\lfloor\log_2 u_j\rfloor+1\) unsigned bits, where \(u_j>0\) is its upper bound; QUBO coefficients are divided by the largest absolute entry of the quadratic matrix.
D-Wave's \texttt{SimulatedAnnealingSampler} in \texttt{dwave-samplers} 1.6.0 uses 100 reads, 1,000 sweeps, and a geometric inverse-temperature schedule~\citep{dwaveSamplers}.
All process workers and the Grid \texttt{FMQA-SA} worker have a 2,100 s watchdog.
The Grid \texttt{RBFOpt}/\texttt{BAM} trace watchdogs are disabled (\nolinkurl{max_run_seconds=0}); archived Grid \texttt{D-SDA} uses budget, target, or local-search termination without an additional method-clock limit.

\subsection{Versions and independently timed process evaluations}\label{sec:supp-acceptance}
Software versions are given in Section~\ifdefined\manuscriptwithsupplement\ref{sec:protocol}\else3.6\fi{} of the main paper.

The conventional GDP process models retain their archived initializations, including constructor defaults for \texttt{CSTR30}; the updated DFO evaluation table is identified separately in Section~\ref{sec:supp-data}.

The acceptance rule in \nolinkurl{gdp_models/evaluation.py} requires all of the following: a finite objective; termination in \{\texttt{optimal}, \texttt{locallyOptimal}, \texttt{globallyOptimal}\}; solver status other than \texttt{error} or \texttt{aborted}; interface status \texttt{SOLVESTAT}=1 and \texttt{MODELSTAT}\(\in\{1,2\}\); and a finite maximum primal violation no greater than \(10^{-5}\).
Time-limit, iteration-limit, evaluation-limit, and error outcomes are not accepted, even if a finite feasible candidate is available.

The primal check uses unscaled violations of original constraints and initialized-variable bounds, including constraint rows deactivated during preprocessing.
Nonfinite objectives, constraint bodies, or initialized-variable values preclude acceptance; unused uninitialized variables are skipped.
Preprocessing rejection or worker failure supplies no accepted objective.

The \texttt{CSTR} initializer \nolinkurl{cstr_box_ipopt.py:initialize_feasible} runs after preprocessing and checks the mass-balance/zero-recycle start at \(10^{-8}\) without changing constraints or bounds.
All process evaluations use cold starts; NLP solutions are not transferred between configurations.
Objectives, acceptance status, and NLP point costs are linked to the same recorded attempt by model, \texttt{CSTR} size, and coordinate, without successful-retry selection or cross-size time substitution.

\subsection{Observed proposal and initialization behavior}\label{sec:supp-proposals}
Table~\ref{tab:supp-fmqa-fallback} pools \texttt{FMQA-SA} proposals across budgets and seeds, excluding the shared initialization and including terminal callbacks.
Fallback supplies 52.89\% of process proposals; the uniform comparison evaluates the complete procedure, without isolating the contributions of learning, annealing, and fallback.

\begin{table}[htbp]
\centering
\caption{Grid-GDP and process GDP models: proposals and uniform fallback proposals of \texttt{FMQA-SA}, pooled over all budgets and seeds, with 1,050 Grid-GDP runs (21 functions, five budgets, and ten seeds) and ten runs per process condition.
A proposal is a configuration that \texttt{FMQA-SA} selects after the three shared initial configurations, so a run with evaluation budget \(B\) makes \(B-3\) proposals.
A fallback is a proposal drawn uniformly from the unevaluated configurations because no annealed candidate is a new configuration within the coordinate bounds; the rate is the percentage of proposals that are fallbacks, and the process total pools the four process models.}
\label{tab:supp-fmqa-fallback}
\footnotesize
\begin{threeparttable}
\begin{tabular}{@{}lrrrr@{}}
\toprule
Benchmark & Runs & Proposals & Fallbacks & Rate (\%) \\
\midrule
Grid-GDP & 1,050 & 51,450 & 9,522 & 18.51 \\
\midrule
\texttt{CSTR} & 60 & 1,350 & 369 & 27.33 \\
\texttt{Column} & 10 & 100 & 86 & 86.00 \\
\texttt{Small Batch} & 10 & 70 & 58 & 82.86 \\
\texttt{Dynamic} & 60 & 1,460 & 1,063 & 72.81 \\
\addlinespace[2pt]
Process total & 140 & 2,980 & 1,576 & 52.89 \\
\bottomrule
\end{tabular}
\begin{tablenotes}[flushleft]\footnotesize
\item Runs: number of \texttt{FMQA-SA} runs.
\end{tablenotes}
\end{threeparttable}
\end{table}
\FloatBarrier

\section{Process model scale, units, and discrete decisions}\label{sec:supp-models}
Table~\ref{tab:supp-model-sizes} counts the original GDP before reformulation or presolve, including all active disjunct branches and fixed continuous variables.
\texttt{Dynamic} is discretized with 30 finite elements per stage and three Radau points.
Binary counts refer to the indicator variables already constructed by the model, before Boolean-to-binary conversion.
All nine models have zero additional nonbinary integer variables.

\begin{table}[htbp]
\centering
\caption{Process GDP models: numbers of variables, constraints, and disjunctions of the original GDP models.
The counts are taken before reformulation or presolve and include all active disjunct branches and fixed continuous variables, and \texttt{Dynamic} is discretized with 30 finite elements per stage and three Radau collocation points.}
\label{tab:supp-model-sizes}
\footnotesize
\begin{threeparttable}
\begin{tabular}{@{}lrrrrrr@{}}
\toprule
Model & \(N\) & \(n_c\) & \(n_b\) & \(m_a\) & \(m_l\) & \(n_d\)\\\midrule
\multirow{6}{*}{\texttt{CSTR}} & 5 & 56 & 20 & 100 & 12 & 10 \\
 & 10 & 106 & 40 & 195 & 22 & 20 \\
 & 15 & 156 & 60 & 290 & 32 & 30 \\
 & 20 & 206 & 80 & 385 & 42 & 40 \\
 & 25 & 256 & 100 & 480 & 52 & 50 \\
 & 30 & 306 & 120 & 575 & 62 & 60 \\
\midrule
\texttt{Column} & -- & 446 & 30 & 637 & 16 & 15 \\
\midrule
\texttt{Small Batch} & -- & 19 & 18 & 34 & 3 & 9 \\
\midrule
\texttt{Dynamic} & -- & 1{,}647 & 27 & 3{,}276 & 31 & 9 \\
\bottomrule
\end{tabular}
\begin{tablenotes}[flushleft]\footnotesize
\item \(N\): number of potential reactors of \texttt{CSTR}; --: not applicable.
\(n_c\): continuous variables; \(n_b\): binary indicator variables constructed by the model; \(m_a\): algebraic constraints; \(m_l\): logical constraints; \(n_d\): disjunctions.
\end{tablenotes}
\end{threeparttable}
\end{table}

\FloatBarrier

The \texttt{Column} objective is \(10^3(Q_c+Q_b)+10^3 n_{\rm active}\), the weighted cost proxy of the source model, which we did not convert to an annualized cost.
The \texttt{Dynamic} objective is \(-\int_0^9 x(t)^2\,dt\).

\section{Lomustine simulation-based benchmark}\label{sec:supp-lomustine}
The \texttt{Lomustine} configurations and responses come from the corrected-model tables of \citet{secquoiaPdIsing}, generated with \texttt{PharmaPy} at commit 3194efd.
For each superstructure option, the source ran \texttt{PyNomadBBO} from a stored initial point on the normalized operating inputs with 60 simulations and stored the best objective, the constraint values, the equipment sizes, and the wall time of the inner optimization.
The inner search accepts a point when the solubility path constraint holds.
The heptane mass-fraction specification of the source is not enforced, because the fixed three-hour vaporizer operation cannot reach it within the pressure bounds.
Table~\ref{tab:supp-lomustine} lists the 36 two-reactor configurations; the four single-reactor semibatch options of the source do not fit the three external variables and are excluded.
The configurations with a batch crystallizer have the lowest objective for every reactor pair, and the choice of the first reactor changes the objective by up to a factor of 2.9.
The vendored table, case module, run records, and analysis are listed in Section~\ref{sec:supp-data}.

\paragraph{Option order and objective.}
The external variables index categorical alternatives, so we repeated the 210 runs with the five other orders of Table~\ref{tab:supp-lomustine-orders}, using the same seeds, initial designs in index space, and budgets.
Table~\ref{tab:supp-lomustine-order-hits} reports the hits in the source order and their range over the six orders.
\texttt{NOMAD} ended with a segmentation fault after identical evaluations on both attempts in 12 runs over all campaigns; as for the Grid runs of Section~\ref{sec:supp-nomad-miso}, these runs are scored on the evaluations completed before the crash.
Table~\ref{tab:supp-lomustine-capex} repeats the source-order runs with the bare-module capital cost of each flowsheet at its optimized operation as the objective, whose reference is the all-batch flowsheet at \$82{,}766.
This cost is evaluated at the operating point that maximizes production and crystal size, and it is not itself optimized over the operating inputs.

\begin{table}[htbp]
\centering
\caption{\texttt{Lomustine}: the six orders of the reactor and crystallizer alternatives used as external-variable indices.
An order assigns consecutive index values to the alternatives of each external variable in the listed sequence.
The source order is used in the main comparison, whose runs are repeated in the five other orders with the same seeds, initial designs in index space, and budgets.}
\label{tab:supp-lomustine-orders}
\footnotesize
\begin{threeparttable}
\begin{tabular}{@{}lll@{}}\toprule
Order & Reactors & Crystallizers\\\midrule
Source & batch, CSTR, PFR & batch, 1, 2, 3 \\
Reversed & PFR, CSTR, batch & 3, 2, 1, batch \\
Back-mixing & batch, PFR, CSTR & batch, 3, 2, 1 \\
Random 1 & PFR, batch, CSTR & 1, 2, batch, 3 \\
Random 2 & batch, PFR, CSTR & batch, 3, 1, 2 \\
Random 3 & CSTR, PFR, batch & 2, 1, batch, 3 \\
\bottomrule
\end{tabular}
\begin{tablenotes}[flushleft]\footnotesize
\item Source: order of the source enumeration; Reversed: source order backward; Back-mixing: alternatives ordered by degree of back-mixing; Random 1 to 3: three further permutations.
CSTR: continuous stirred-tank reactor; PFR: plug-flow reactor; 1, 2, 3: number of mixed-suspension, mixed-product-removal (MSMPR) crystallizer stages.
\end{tablenotes}
\end{threeparttable}
\end{table}

\begin{table}[htbp]
\centering
\caption{\texttt{Lomustine}: target hits of the six DFO methods and uniform sampling in the source order and over the six option orders of Table~\ref{tab:supp-lomustine-orders}, with ten runs per method, budget, and order.
A target hit is a run that reaches the reference target of Equation~\ifdefined\manuscriptwithsupplement\eqref{eqn:target}\else(8)\fi{} of the main paper, here the reference configuration of Table~\ref{tab:supp-lomustine}.
All orders share the seeds, initial designs in index space, and budgets of the source order, and the rows for the fewest and most target hits give the smallest and largest count over the six orders.}
\label{tab:supp-lomustine-order-hits}
\footnotesize
\setlength{\tabcolsep}{3.5pt}
\begin{threeparttable}
\begin{tabular}{@{}rlrrrrrrr@{}}\toprule
 & & \multicolumn{7}{c}{Target hits (of 10)}\\\cmidrule(l){3-9}
\(B\) & & \texttt{D-SDA} & \texttt{NOMAD} & \texttt{RBFOpt} & \texttt{MISO} & \texttt{BAM} & \texttt{FMQA-SA} & Uniform\\\midrule
\multirow{3}{*}{10} & Source order & 1 & 6 & 8 & 7 & 0 & 2 & 2 \\
 & Fewest over orders & 1 & 5 & 2 & 4 & 0 & 2 & 0 \\
 & Most over orders & 6 & 6 & 8 & 7 & 2 & 7 & 4 \\
\midrule
\multirow{3}{*}{15} & Source order & 4 & 6 & 9 & 10 & 4 & 3 & 4 \\
 & Fewest over orders & 3 & 5 & 7 & 5 & 0 & 3 & 2 \\
 & Most over orders & 8 & 9 & 10 & 10 & 10 & 10 & 4 \\
\midrule
\multirow{3}{*}{20} & Source order & 4 & 6 & 10 & 10 & 6 & 7 & 5 \\
 & Fewest over orders & 3 & 5 & 9 & 6 & 0 & 7 & 4 \\
 & Most over orders & 8 & 9 & 10 & 10 & 10 & 10 & 7 \\
\bottomrule
\end{tabular}
\begin{tablenotes}[flushleft]\footnotesize
\item \(B\): evaluation budget; Uniform: uniform sampling without replacement; Source order: order of the source enumeration, used in the main comparison.
\end{tablenotes}
\end{threeparttable}
\end{table}

\begin{table}[htbp]
\centering
\caption{\texttt{Lomustine} with the capital-cost objective: target hits of the six DFO methods and uniform sampling in the source order, with ten runs per method and budget.
The capital-cost objective is the bare-module capital cost of each flowsheet at its optimized operation, and a target hit is a run that reaches the reference of this objective, the all-batch configuration at \$82{,}766.}
\label{tab:supp-lomustine-capex}
\footnotesize
\setlength{\tabcolsep}{3.5pt}
\begin{threeparttable}
\begin{tabular}{@{}rrrrrrrr@{}}\toprule
 & \multicolumn{7}{c}{Target hits (of 10)}\\\cmidrule(l){2-8}
\(B\) & \texttt{D-SDA} & \texttt{NOMAD} & \texttt{RBFOpt} & \texttt{MISO} & \texttt{BAM} & \texttt{FMQA-SA} & Uniform\\\midrule
10 & 3 & 9 & 7 & 8 & \textbf{10} & 2 & 4 \\
15 & 6 & 9 & 8 & \textbf{10} & \textbf{10} & 3 & 7 \\
20 & 9 & 9 & \textbf{10} & \textbf{10} & \textbf{10} & 5 & 8 \\
\bottomrule
\end{tabular}
\begin{tablenotes}[flushleft]\footnotesize
\item \(B\): evaluation budget; Uniform: uniform sampling without replacement.
Bold: most target hits for each budget, including ties.
\end{tablenotes}
\end{threeparttable}
\end{table}

\begin{table}[htbp]
\centering
\caption{\texttt{Lomustine}: objective of every configuration after the simulation-based optimization of its operating inputs, and the wall time of that optimization, with one optimization per configuration.
The objective is \(-(6\times10^{4}\,\dot m+d_{10})\), and each configuration combines a first reactor, a second reactor, and a crystallizer.
The reference configuration, a plug-flow reactor followed by a batch reactor and the batch crystallizer, has the lowest objective.}
\label{tab:supp-lomustine}
\footnotesize
\setlength{\tabcolsep}{3.5pt}
\begin{threeparttable}
\begin{tabular}{@{}llrrrrrrrr@{}}\toprule
 & & \multicolumn{4}{c}{Objective} & \multicolumn{4}{c}{Wall time (s)}\\\cmidrule(lr){3-6}\cmidrule(l){7-10}
R01 & R02 & batch & 1 MSMPR & 2 MSMPR & 3 MSMPR & batch & 1 MSMPR & 2 MSMPR & 3 MSMPR\\\midrule
\multirow{3}{*}{batch} & batch & -170.68 & -131.70 & -146.86 & -148.92 & 91 & 310 & 538 & 835 \\
 & CSTR & -230.52 & -175.25 & -175.27 & -166.13 & 146 & 348 & 650 & 970 \\
 & PFR & -268.46 & -209.36 & -213.29 & -202.26 & 159 & 348 & 651 & 975 \\
\addlinespace[2pt]
\multirow{3}{*}{CSTR} & batch & -483.85 & -360.16 & -365.91 & -345.28 & 152 & 331 & 656 & 936 \\
 & CSTR & -331.90 & -244.39 & -245.34 & -239.50 & 159 & 367 & 660 & 1{,}019 \\
 & PFR & -429.81 & -315.22 & -321.46 & -304.62 & 206 & 408 & 717 & 1{,}017 \\
\addlinespace[2pt]
\multirow{3}{*}{PFR} & batch & \textbf{-487.84} & -362.45 & -361.55 & -348.06 & 208 & 347 & 655 & 886 \\
 & CSTR & -351.71 & -249.82 & -259.47 & -244.65 & 183 & 388 & 659 & 983 \\
 & PFR & -433.41 & -318.34 & -324.18 & -306.21 & 211 & 393 & 629 & 907 \\
\bottomrule
\end{tabular}
\begin{tablenotes}[flushleft]\footnotesize
\item \(\dot m\): production rate (kg/s); \(d_{10}\): number-weighted mean crystal size (\(\mu\)m); R01, R02: first and second reactor; CSTR: continuous stirred-tank reactor; PFR: plug-flow reactor; batch, 1 to 3 MSMPR: batch crystallizer or one to three mixed-suspension, mixed-product-removal (MSMPR) stages.
Bold: reference configuration.
\end{tablenotes}
\end{threeparttable}
\end{table}
\FloatBarrier

\section{Catalytic distillation column}\label{sec:supp-catdist}

\paragraph{Corrected model.}
We checked every equation block, parameter table, and bound of the published model against standard forms and literature values, and the \texttt{Catalytic Column} benchmark uses the model with the errors fixed.
The corrections concern the vapor--liquid equilibrium and the equilibrium constant of the reaction, the mixture liquid density (Rackett equation with the Spencer--Danner mixing rule) and surface tension (Macleod--Sugden), the tray hydraulics, and the objective, which counted the annualization factor twice and credited the impurities in the bottoms at the ETBE price.
The pure-component data of isobutene, 1-butene, and ethanol are fitted to open reference equations of state, and the ETBE data, for which no reference equation exists, agree with open measurements of its vapor pressure and density.
The bottoms must contain at least 95~mol\% ETBE, the fuel-grade specification, the catalytic stages stay at or below 120~\(^\circ\)C, the maximum operating temperature of the Amberlyst~15 catalyst, and the prices are 2010 market prices.
The repository lists each finding with its effect on the published optimum; with the 95~mol\% specification, the published optimal configuration is infeasible, because its column is too short.

\paragraph{Evaluation campaign.}
The 210 runs made 12,015 evaluations: 6,921 returned accepted objectives, 2,789 were logical rejections without an NLP solve, and 2,305 returned unaccepted NLP responses.
Eight \texttt{D-SDA} runs stopped before their budgets because no neighbor improved, and all other runs used their full budgets.
The runs requested 3,499 distinct configurations, of which 2,652 were accepted, and solving them took 7.7~h.
The completion times of the runs sum to 20.2~h, but the campaign took 2.8~h on three cores because a configuration requested by several runs was solved once.
Each run evaluated distinct configurations, so the reuse comes from different runs requesting the same configuration.

\paragraph{Search box.}
The GDP admits 100,947 configurations, and the gap box \(g_i\le6\) of the search methods keeps 54,159 of them, including every design within 5\% of the best of a pilot of 202 configurations.
The reference has \(g_1=f_1-2=6\), on the bound of the box, so a better design could lie just outside it.
We therefore solved, with the initialization and acceptance rule of the campaign, all 747 configurations outside the box within two stages of the reference in every stage position.
Of these, 255 are accepted and 492 are locally infeasible, and none is better than the reference: the best, \((9,9,13,17,17,21)\) at 24,930.96~\$/yr, lies 0.14\% above it.
This check ran on the server of the campaign with its \texttt{GAMS} and \texttt{CONOPT} versions, and the longest solve took 40~s, below the campaign's 300~s limit; an earlier run of the same check on a workstation with \texttt{GAMS} 54.2.2 returned the same status for every configuration and the same objective for every accepted one.

\paragraph{MINLP forms.}
The published MINLP form makes the stage decisions binary variables, multiplies switched equations by these variables or by the stage-existence indicator, and keeps the property and kinetic equations on every stage.
The disjunctive form behind the Big-\(M\) and Hull reformulations has five disjunctions for each internal stage, one per decision, whose terms collect the equations that the GDP model switches on or off (Table~\ref{tab:supp-catdist-disjunctions}).
The property and feed-enthalpy equations hold on every stage, as in the published form, but the kinetics stay in the catalyst disjunction; with the kinetics on every stage, the bounds of the reaction rate restricted the temperature of a stage without catalyst, and with the binary variables of \((3,4,5,11,11,14)\) fixed, the Big-\(M\) reformulation stopped at 32,075~\$/yr, a point from which the GDP form continues to 30,822~\$/yr.
Each nonlinear disjunct equation \(h(\mathbf x)=0\) or \(h(\mathbf x)\le 0\) is replaced by the definition \(s\,r=h(\mathbf x)\), which holds on every stage, and by the linear disjunct constraint \(r=0\) or \(r\le 0\); a linear disjunct equation is divided by \(s\).
The scale \(s\) is 1.01 times the largest absolute value of \(h\) over the variable bounds, computed by interval arithmetic after feasibility-based bound tightening with the constraints that hold for every configuration.
Every scaled disjunct constraint therefore lies in \([-1,1]\) at every feasible point, so \(M=1\) is a valid Big-\(M\) value for all of them, and the bound tightening adds 101 implied variable bounds that are written into the model.
The scales range from 0.10 to \(1.3\times10^{18}\); the largest belongs to the rate equation, whose products of activities, adsorption constant, and rate span this range over the variable bounds.
The interval bounds of the flooding and weeping correlations are finite only if the liquid is denser than the vapor, which the source bounds do not imply, so the disjunctive form bounds the liquid and vapor mass densities by \(\rho_L\ge200\) and \(5\le\rho_V\le150\)~kg/m\(^3\); every validated solution has \(\rho_L\ge529\) and \(\rho_V\le34\)~kg/m\(^3\).
These bounds belong only to the Big-\(M\) and Hull reformulations; the GDP model, its fixed-configuration evaluations, and the published MINLP form do not impose them.
Within the model's bounds (pressure 9.5 to 10~bar, temperature 200 to 445~K, compressibility factor 0.5 to 1.3) and for any vapor composition, whose molar mass lies between 46.07 and 102.18~g/mol, the vapor density lies between 9 and 123~kg/m\(^3\), so the vapor bounds cannot bind and only the liquid bound restricts the reformulations.
The \texttt{GAMS} convex hull reformulation accepts only inequalities, so each disjunct equality is written as two inequalities, and Big-\(M\) of the pair equals Big-\(M\) of the equality.
Table~\ref{tab:supp-catdist-validation} compares the forms with the binary variables fixed.

\begin{table}[htbp]
\centering
\caption{\texttt{Catalytic Column}: the five disjunctions of each internal stage \(k\) in the disjunctive form behind the Big-\(M\) and Hull reformulations.
Each row is one disjunction, whose Boolean indicator enforces the equations of the second column on stage \(k\) when it is true and those of the third column when it is false.
The property and feed-enthalpy equations hold on every stage and belong to no disjunction.}
\label{tab:supp-catdist-disjunctions}
\footnotesize
\begin{threeparttable}
\begin{tabular}{@{}lll@{}}\toprule
Indicator & Indicator true & Indicator false\\\midrule
stage \(k\) exists & \makecell[l]{vapor--liquid equilibrium, tray hydraulics,\\pressure drop, flooding and weeping limits,\\tray height} & \makecell[l]{vapor composition and flow pass through,\\no pressure drop, no tray height} \\
catalyst on stage \(k\) & \makecell[l]{kinetics, catalyst hydraulics and flooding,\\reaction terms, 120~\(^\circ\)C limit} & no reaction, no catalyst pressure drop \\
ethanol feed at \(k\) & feed flow & no feed flow \\
butenes feed at \(k\) & feed flow & no feed flow \\
boil-up return at \(k\) & boil-up flow & no boil-up flow \\
\bottomrule
\end{tabular}
\begin{tablenotes}[flushleft]\footnotesize
\item \(k\): internal stage of the column; Indicator: Boolean indicator of the disjunction.
\end{tablenotes}
\end{threeparttable}
\end{table}

\begin{table}[htbp]
\centering
\caption{\texttt{Catalytic Column}: objectives of four model forms with the binary variables fixed at nine configurations, with one \texttt{CONOPT} solve per form and configuration from the initialization of the published model.
The four forms are the GDP model evaluated by the search methods, the MINLP form published with the original \texttt{D-SDA} study of this column, and the Big-\(M\) and Hull reformulations of the disjunctive form of Table~\ref{tab:supp-catdist-disjunctions}.
The first row is the reference configuration.}
\label{tab:supp-catdist-validation}
\footnotesize
\begin{threeparttable}
\begin{tabular}{@{}lrrrr@{}}\toprule
 & \multicolumn{4}{c}{Objective (\$/yr)}\\\cmidrule(l){2-5}
\((f_1,c_1,c_2,c_3,f_2,b)\) & GDP & Published MINLP & Big-\(M\) & Hull\\\midrule
\((8,8,12,17,17,21)\) & 24{,}895.55 & 24{,}895.55 & 24{,}895.55 & 24{,}895.55 \\
\((7,8,11,16,16,20)\) & 25{,}271.33 & 25{,}271.33 & 25{,}271.33 & 25{,}271.33 \\
\((3,4,5,11,11,14)\) & 30{,}822.02 & 30{,}822.02 & 30{,}822.02 & 30{,}822.02 \\
\((4,4,9,10,14,17)\) & 28{,}373.98 & 28{,}373.98 & 28{,}373.98 & 28{,}373.98 \\
\((3,6,12,14,19,21)\) & 31{,}116.60 & 32{,}320.55 & 31{,}116.60 & 31{,}116.60 \\
\((2,3,9,10,20,21)\) & 35{,}745.43 & 35{,}194.65 & 35{,}745.43 & 35{,}745.43 \\
\((4,4,6,8,8,10)\) & \multicolumn{4}{c}{infeasible} \\
\((3,4,5,6,8,9)\) & \multicolumn{4}{c}{infeasible} \\
\((13,13,15,16,17,19)\) & \multicolumn{4}{c}{infeasible} \\
\bottomrule
\end{tabular}
\begin{tablenotes}[flushleft]\footnotesize
\item \(f_1\), \(f_2\): stages of the ethanol and butenes feeds; \(c_1<c_2<c_3\): catalytic stages; \(b\): stage of the boil-up return.
infeasible: every form reports the configuration as locally infeasible.
\end{tablenotes}
\end{threeparttable}
\end{table}

\paragraph{MINLP solvers.}
Table~\ref{tab:supp-catdist-minlp} lists the results of the MINLP solvers.
\texttt{SBB} stops on the Big-\(M\) and Hull reformulations when the NLP relaxations of its nodes exceed the limit on function evaluations outside their domains.
Table~\ref{tab:supp-catdist-certificates} lists the global optimality certificates that \texttt{SCIP} and \texttt{BARON} return when they start from a feasible design of the published form; the \texttt{GAMS} \texttt{Examiner} confirms that both starting designs satisfy every constraint and bound within \(10^{-6}\).

\begin{table}[htbp]
\centering
\caption{\texttt{Catalytic Column}: outcomes and solve times of four MINLP solvers on three model forms, with one run per solver and form started from the initial configuration \((3,4,5,6,8,9)\) of the published model.
The outcome is the objective and configuration of the design found, or the reason no feasible design was returned; the reference design has an objective of 24,895.55~\$/yr.
The time is the wall-clock time of the \texttt{GAMS} job, and each run uses one core and a 3,600~s time limit.}
\label{tab:supp-catdist-minlp}
\footnotesize
\begin{threeparttable}
\begin{tabular}{@{}lllr@{}}\toprule
Form & Solver & Outcome & Time (s)\\\midrule
\multirow{4}{*}{Published} & \texttt{DICOPT} & no feasible design; every NLP subproblem infeasible, iteration limit & 686 \\
 & \texttt{SBB} & 28{,}717.22~\$/yr at \((4,4,5,6,6,9)\) & 10 \\
 & \texttt{BARON} & declared infeasible & 255 \\
 & \texttt{SCIP} & declared infeasible & 9 \\
\midrule
\multirow{4}{*}{Big-\(M\)} & \texttt{DICOPT} & no feasible design; NLP subproblem failed & 902 \\
 & \texttt{SBB} & no feasible design; domain violation limit & 38 \\
 & \texttt{BARON} & declared infeasible & 214 \\
 & \texttt{SCIP} & no feasible design; lower bound \(-5.7\times10^{7}\)~\$/yr & 3{,}600+ \\
\midrule
\multirow{4}{*}{Hull} & \texttt{DICOPT} & no feasible design; NLP subproblem failed & 211 \\
 & \texttt{SBB} & no feasible design; domain violation limit & 130 \\
 & \texttt{BARON} & declared infeasible & 362 \\
 & \texttt{SCIP} & no feasible design; lower bound \(-5.7\times10^{7}\)~\$/yr & 3{,}600+ \\
\bottomrule
\end{tabular}
\begin{tablenotes}[flushleft]\footnotesize
\item Published: published MINLP form; Big-\(M\), Hull: reformulations of the disjunctive form; configurations \((f_1,c_1,c_2,c_3,f_2,b)\) as in Table~\ref{tab:supp-catdist-validation}; 3,600+: stopped at the time limit.
\end{tablenotes}
\end{threeparttable}
\end{table}

\begin{table}[htbp]
\centering
\caption{\texttt{Catalytic Column}: global optimality certificates that \texttt{SCIP} and \texttt{BARON} return for the published MINLP form, with one run per solver from each of two feasible starting designs.
The first starting design comes from the pilot of 202 configurations and the second is the reference configuration.
The status and objective are those reported by each solver, and the time is its solve time.}
\label{tab:supp-catdist-certificates}
\footnotesize
\begin{threeparttable}
\begin{tabular}{@{}llrrlrr@{}}\toprule
 & \multicolumn{3}{c}{\texttt{SCIP}} & \multicolumn{3}{c}{\texttt{BARON}}\\\cmidrule(lr){2-4}\cmidrule(l){5-7}
Starting design & Status & Objective (\$/yr) & Time (s) & Status & Objective (\$/yr) & Time (s)\\\midrule
\((7,8,11,16,16,20)\) & optimal & 25,271.33 & 8 & optimal & 25,271.33 & 172\\
\((8,8,12,17,17,21)\) & optimal & 24,895.55 & 10 & optimal & 24,895.55 & 270\\
\bottomrule
\end{tabular}
\begin{tablenotes}[flushleft]\footnotesize
\item Starting design: configuration \((f_1,c_1,c_2,c_3,f_2,b)\) as in Table~\ref{tab:supp-catdist-validation}.
\end{tablenotes}
\end{threeparttable}
\end{table}
\FloatBarrier

\section{Completion costs}\label{sec:supp-uncertainty}

Table~\ref{tab:process-search} reports the search times in all 14 process conditions; the completion times are in Table~\ifdefined\manuscriptwithsupplement\ref{tab:process-targets}\else2\fi{} and the enumeration times in Table~\ifdefined\manuscriptwithsupplement\ref{tab:process-cases}\else1\fi{} of the main paper.
Table~\ref{tab:grid-cost} reports the Grid-GDP costs at \(B=100\).
\texttt{MISO} reaches the most Grid targets, but its median search time of 32.70~s is about four times that of \texttt{RBFOpt}, whereas \texttt{D-SDA} and \texttt{NOMAD} stop after medians of 37 and 31.5 evaluations and have the lowest median completion times, 18.38 and 18.50~s (\texttt{NOMAD} over its 189 timed runs).
\begin{table}[p]
\centering
\caption{Process GDP models: search times of the six DFO methods and uniform sampling, with ten runs per method and condition.
The search time of a run is the time the method spends selecting configurations, which excludes the cost of evaluating them, and the medians and percentiles are over the ten runs.
Table~\ifdefined\manuscriptwithsupplement\ref{tab:process-targets}\else2\fi{} of the main paper gives the full-run completion times and target hits.}
\label{tab:process-search}
\footnotesize
\setlength{\tabcolsep}{3.5pt}
\begin{threeparttable}
\begin{tabular}{@{}lrrlrrrrrrr@{}}\toprule
Model & \(N\) & \(B\) & Search time & \texttt{D-SDA} & \texttt{NOMAD} & \texttt{RBFOpt} & \texttt{MISO} & \texttt{BAM} & \texttt{FMQA-SA} & Uniform \\\midrule
\multirow{18}{*}{\texttt{CSTR}} & \multirow{3}{*}{5} & \multirow{3}{*}{8} & Median (s) & 0.0105 & 0.0367 & 0.0596 & 0.659 & 0.0360 & 3.66 & 0.00460 \\
 & &  & 25th percentile (s) & 0.0102 & 0.0323 & 0.0590 & 0.658 & 0.0151 & 3.64 & 0.00387 \\
 & &  & 75th percentile (s) & 0.0111 & 1.13 & 0.0597 & 0.669 & 0.0816 & 3.67 & 0.00478 \\
 & \multirow{3}{*}{10} & \multirow{3}{*}{15} & Median (s) & 0.0137 & 1.83 & 0.302 & 0.677 & 0.137 & 8.00 & 0.00807 \\
 & &  & 25th percentile (s) & 0.0130 & 0.469 & 0.288 & 0.675 & 0.126 & 7.98 & 0.00772 \\
 & &  & 75th percentile (s) & 0.0147 & 1.88 & 0.315 & 0.683 & 0.158 & 8.08 & 0.00812 \\
 & \multirow{3}{*}{15} & \multirow{3}{*}{22} & Median (s) & 0.0169 & 0.825 & 1.00 & 3.29 & 0.150 & 12.4 & 0.0113 \\
 & &  & 25th percentile (s) & 0.0156 & 0.0949 & 0.492 & 1.21 & 0.140 & 12.3 & 0.0108 \\
 & &  & 75th percentile (s) & 0.0172 & 1.82 & 3.26 & 7.50 & 0.162 & 12.4 & 0.0116 \\
 & \multirow{3}{*}{20} & \multirow{3}{*}{29} & Median (s) & 0.0207 & 1.63 & 2.82 & 10.8 & 0.167 & 16.9 & 0.0139 \\
 & &  & 25th percentile (s) & 0.0189 & 0.521 & 0.632 & 9.02 & 0.165 & 16.8 & 0.0121 \\
 & &  & 75th percentile (s) & 0.0218 & 1.92 & 4.01 & 15.8 & 0.313 & 16.9 & 0.0148 \\
 & \multirow{3}{*}{25} & \multirow{3}{*}{36} & Median (s) & 0.0258 & 1.08 & 2.47 & 17.6 & 0.401 & 21.3 & 0.0161 \\
 & &  & 25th percentile (s) & 0.0238 & 0.170 & 1.86 & 16.4 & 0.392 & 21.3 & 0.0147 \\
 & &  & 75th percentile (s) & 0.0299 & 2.32 & 4.43 & 18.8 & 0.449 & 21.4 & 0.0186 \\
 & \multirow{3}{*}{30} & \multirow{3}{*}{43} & Median (s) & 0.0585 & 1.90 & 4.99 & 18.1 & 0.425 & 25.8 & 0.0241 \\
 & &  & 25th percentile (s) & 0.0556 & 1.76 & 2.95 & 14.0 & 0.404 & 25.7 & 0.0224 \\
 & &  & 75th percentile (s) & 0.0610 & 7.65 & 10.3 & 19.3 & 0.451 & 26.0 & 0.0257 \\
\midrule
\multirow{3}{*}{\texttt{Column}} & \multirow{3}{*}{--} & \multirow{3}{*}{13} & Median (s) & 0.0120 & 0.0423 & 0.226 & 0.672 & 0.0177 & 6.85 & 0.00716 \\
 & &  & 25th percentile (s) & 0.0116 & 0.0358 & 0.205 & 0.669 & 0.0168 & 6.81 & 0.00711 \\
 & &  & 75th percentile (s) & 0.0124 & 0.0461 & 0.239 & 0.681 & 0.155 & 6.89 & 0.00748 \\
\midrule
\multirow{3}{*}{\texttt{Small Batch}} & \multirow{3}{*}{--} & \multirow{3}{*}{10} & Median (s) & 0.0102 & 0.0475 & 0.170 & 0.668 & 0.0169 & 4.87 & 0.00536 \\
 & &  & 25th percentile (s) & 0.00971 & 0.0333 & 0.159 & 0.666 & 0.0161 & 4.86 & 0.00506 \\
 & &  & 75th percentile (s) & 0.0103 & 0.0540 & 0.180 & 0.670 & 0.0191 & 4.89 & 0.00545 \\
\midrule
\multirow{18}{*}{\texttt{Dynamic}} & \multirow{3}{*}{--} & \multirow{3}{*}{12} & Median (s) & 0.0119 & 0.0419 & 0.186 & 0.671 & 0.0173 & 6.16 & 0.00589 \\
 & &  & 25th percentile (s) & 0.0113 & 0.0377 & 0.182 & 0.670 & 0.0172 & 6.14 & 0.00540 \\
 & &  & 75th percentile (s) & 0.0125 & 0.0630 & 0.189 & 0.676 & 0.0174 & 6.20 & 0.00658 \\
 &  & \multirow{3}{*}{15} & Median (s) & 0.0123 & 0.0896 & 0.296 & 0.679 & 0.139 & 8.04 & 0.00740 \\
 & &  & 25th percentile (s) & 0.0115 & 0.0627 & 0.291 & 0.676 & 0.128 & 8.01 & 0.00622 \\
 & &  & 75th percentile (s) & 0.0131 & 1.09 & 0.321 & 0.683 & 0.156 & 8.09 & 0.00760 \\
 &  & \multirow{3}{*}{20} & Median (s) & 0.0133 & 1.49 & 0.367 & 25.0 & 0.143 & 11.2 & 0.00987 \\
 & &  & 25th percentile (s) & 0.0119 & 0.447 & 0.355 & 24.6 & 0.130 & 11.1 & 0.00927 \\
 & &  & 75th percentile (s) & 0.0154 & 6.90 & 0.381 & 30.1 & 0.156 & 11.2 & 0.00991 \\
 &  & \multirow{3}{*}{30} & Median (s) & 0.0149 & 1.77 & 7.61 & 79.7 & 0.355 & 17.4 & 0.0141 \\
 & &  & 25th percentile (s) & 0.0124 & 1.44 & 4.88 & 79.4 & 0.204 & 17.2 & 0.0135 \\
 & &  & 75th percentile (s) & 0.0157 & 9.24 & 13.7 & 80.5 & 0.364 & 17.5 & 0.0145 \\
 &  & \multirow{3}{*}{37} & Median (s) & 0.0134 & 1.77 & 10.1 & 81.0 & 0.376 & 21.7 & 0.0173 \\
 & &  & 25th percentile (s) & 0.0118 & 1.43 & 5.25 & 79.7 & 0.364 & 21.6 & 0.0171 \\
 & &  & 75th percentile (s) & 0.0160 & 9.23 & 14.2 & 82.9 & 0.383 & 21.8 & 0.0174 \\
 &  & \multirow{3}{*}{50} & Median (s) & 0.0135 & 1.77 & 26.3 & 158 & 0.732 & 29.8 & 0.0233 \\
 & &  & 25th percentile (s) & 0.0129 & 1.42 & 17.0 & 157 & 0.685 & 29.7 & 0.0229 \\
 & &  & 75th percentile (s) & 0.0153 & 9.21 & 37.2 & 160 & 0.751 & 30.1 & 0.0235 \\
\bottomrule
\end{tabular}
\begin{tablenotes}[flushleft]\footnotesize
\item \(N\): number of potential reactors of \texttt{CSTR}; \(B\): evaluation budget; --: not applicable; Uniform: uniform sampling without replacement.
\end{tablenotes}
\end{threeparttable}
\end{table}

\begin{table}[htbp]
\centering
\caption{Grid-GDP: target hits, full-run completion times, summed evaluation costs, and search times of the six DFO methods and the two sampling controls with a budget of 100 evaluations, with 210 runs per method (21 test functions and ten seeds).
A target hit is a run that reaches the reference target of Equation~\ifdefined\manuscriptwithsupplement\eqref{eqn:target}\else(8)\fi{} of the main paper, a normalized objective within \(10^{-5}\) of the normalized reference objective zero.
The summed evaluation cost of a run adds the separately measured costs of all its evaluations, including repeated configurations, the search time is the method time excluding this cost, and the full-run completion time is their sum (Equation~\ifdefined\manuscriptwithsupplement\eqref{eqn:reconstructed-cost}\else(10)\fi{} of the main paper).
Medians and linearly interpolated 25th and 75th percentiles are taken over all runs with recorded times, including target misses, which leaves out the times of the 21 \texttt{NOMAD} runs marked \(a\), and separately for each quantity, so the percentiles of the two components need not sum to those of the completion time.}
\label{tab:grid-cost}
\footnotesize
\setlength{\tabcolsep}{3.5pt}
\begin{threeparttable}
\begin{tabular}{@{}lrrrrrrrrrr@{}}\toprule
 & & \multicolumn{3}{c}{Full-run completion time (s)} & \multicolumn{3}{c}{Summed evaluation cost (s)} & \multicolumn{3}{c}{Search time (s)}\\\cmidrule(lr){3-5}\cmidrule(lr){6-8}\cmidrule(l){9-11}
Method & \makecell{Target hits\\(of 210)} & Median & 25th pct. & 75th pct. & Median & 25th pct. & 75th pct. & Median & 25th pct. & 75th pct.\\\midrule
\texttt{D-SDA} & 101 & 18.38 & 12.42 & 37.06 & 18.38 & 12.41 & 37.05 & 0.0024 & 0.0023 & 0.0026 \\
\texttt{NOMAD}\tnote{a} & 126 & 18.50 & 13.46 & 35.00 & 15.20 & 11.09 & 30.82 & 2.04 & 1.79 & 2.28 \\
\texttt{RBFOpt} & 171 & 70.81 & 54.83 & 93.09 & 60.84 & 44.31 & 82.15 & 7.58 & 4.12 & 14.82 \\
\texttt{MISO} & 181 & 94.82 & 79.22 & 116.05 & 63.21 & 44.03 & 82.37 & 32.70 & 25.19 & 39.38 \\
\texttt{BAM} & 147 & 64.42 & 48.74 & 81.67 & 63.97 & 48.34 & 81.22 & 0.43 & 0.37 & 0.51 \\
\texttt{FMQA-SA} & 129 & 124.13 & 105.57 & 139.26 & 62.34 & 43.98 & 77.52 & 62.04 & 61.58 & 62.53 \\
Uniform & 45 & 60.51 & 43.33 & 78.60 & 60.51 & 43.33 & 78.59 & 0.0026 & 0.0026 & 0.0027 \\
\makecell[l]{\texttt{BAM}-prefix\\control} & 124 & 61.81 & 40.96 & 78.33 & 61.80 & 40.95 & 78.33 & 0.0026 & 0.0026 & 0.0027 \\
\bottomrule
\end{tabular}
\begin{tablenotes}[flushleft]\footnotesize
\item pct.: percentile; Uniform: uniform sampling without replacement after the three shared initial configurations; \texttt{BAM}-prefix control: the ten predetermined sampling points of \texttt{BAM} after the three shared initial configurations, followed by uniform sampling without replacement.
\item[a] \texttt{NOMAD} times cover the 189 runs with recorded times; the 21 runs that crashed on every attempt have none, but their target hits are counted (Section~\ref{sec:supp-nomad-miso}).
\end{tablenotes}
\end{threeparttable}
\end{table}

Per-function Grid results at \(B=100\) and complete seed-level outcomes are available in the analysis files linked in Section~\ref{sec:supp-data}.

\FloatBarrier
\begingroup
\section{Data sources and reproducibility}\label{sec:supp-data}

The \href{https://github.com/SECQUOIA/Black-Box-for-GDP/tree/ea65c5339486b576008df68645ea68d22f7c7213}{publication snapshot} contains the source code, experimental data, analysis results, and reproduction instructions for every experiment in this article.
The \texttt{CSTR} updates are the penalty of \(10^6\) for unaccepted responses, analyzed in \nolinkurl{results/cstr_penalty_1e6_analysis_20261022}, and the \texttt{CSTR30} evaluation table rebuilt with the mass-balance/zero-recycle initialization, in \nolinkurl{results/cstr30_publication_20261025}; the \texttt{Column} runs on its 64-configuration logical domain with \(B=13\) are in \nolinkurl{results/column_logical_b13_analysis_20261025}.
The updated \href{https://github.com/SECQUOIA/Black-Box-for-GDP/tree/329c6f93f4940436425582ae60323c0c9f95b581/results/gdpopt_scip_ipopt_analysis_20261030}{\texttt{GDPopt} analysis and execution records} contain the 63 Grid and 27 process runs configured with \texttt{SCIP} for the mixed-integer master problems and \texttt{IPOPT} for the NLP subproblems.
The three presolve-disabled LOA diagnostic runs are excluded.
The same records contain the paired comparison with the earlier configuration that used \texttt{SCIP} for both the mixed-integer master problems and the NLP subproblems, in which GLOA returned the same \texttt{CSTR} objectives within a relative difference of \(3\times10^{-5}\).
Seed spread, tolerance sensitivity, search overhead, and \texttt{GDPopt} enumeration counts are generated by \nolinkurl{review_comparison/revision_diagnostics.py} from the retained run records and execution logs.
The \texttt{D-SDA} rerun, including the process restart variants, is \nolinkurl{experiments/dsda_no_revisit/run.py}, which replays the sealed lookup package of the adopted \texttt{D-SDA} rerun.
The \texttt{NOMAD} and \texttt{MISO} runs were executed at \href{https://github.com/SECQUOIA/Black-Box-for-GDP/tree/c4d7db6b671bbe92a314be524a650d75cd49529c}{source revision c4d7db6}; the \href{https://github.com/SECQUOIA/Black-Box-for-GDP/tree/841ca986f175dda4d91c784b195e60279cf94b99/results/pr56_dantzig_return_20261031}{returned records} contain the original archive, its SHA-256 checksum, and hash-checked analysis tables of every run.
All \texttt{D-SDA} runs, including the restart variants, were rerun at \href{https://github.com/SECQUOIA/Black-Box-for-GDP/tree/c186408728c6f95ac33a5465edd53035d85c5620}{source revision c186408}, where no configuration is evaluated twice; with re-evaluation switched back on, the same code reproduces all 1,190 published \texttt{D-SDA} runs exactly, and the \href{https://github.com/SECQUOIA/Black-Box-for-GDP/tree/f3b57ac44cccb27552b3a8618ce6abd463e540a4/results/dsda_no_revisit_20261031}{returned records} hold every run.
The \texttt{Lomustine} case is \nolinkurl{experiments/lomustine_simulation}, which vendors the corrected-model table of \citet{secquoiaPdIsing} and runs the seven methods through the harness in \nolinkurl{experiments/catalytic_distillation_live}; the run records and analysis are in \nolinkurl{results/lomustine_simulation_20261031}.
The \texttt{Catalytic Column} model and its review are in \nolinkurl{gdp_models/catalytic_distillation}, the live runs are \nolinkurl{experiments/catalytic_distillation_live}, and their returned records, with the original archive, its SHA-256 checksum, and an audit of every run, are in \nolinkurl{results/catalytic_distillation_live_20261101}.
The MINLP forms, their generators, the fixed-configuration validation, and the solver runs are \nolinkurl{experiments/catalytic_distillation_minlp}; the runs reported here are in \nolinkurl{results/catalytic_distillation_minlp_dantzig_20261101}, and an earlier run of the same experiments on a workstation with \texttt{GAMS} 54.2.2, in \nolinkurl{results/catalytic_distillation_minlp_20261101}, returned the same outcome for every solver, form, and configuration.
The \texttt{D-SDA} neighborhood comparison of Section~\ref{sec:supp-dsda-neighborhood} is \nolinkurl{experiments/dsda_separable}, with results in \nolinkurl{results/dsda_separable_20261101}, and the check of the \texttt{Catalytic Column} configurations outside the search box is in \nolinkurl{results/catalytic_distillation_outside_box_dantzig_20261103}, and its earlier workstation run in \nolinkurl{results/catalytic_distillation_outside_box_20261102}.
The seven-method tables and figures are generated from these tables and the published records by \nolinkurl{review_comparison/seven_methods.py} and the \nolinkurl{review_comparison/joint_process_*.py} scripts, without new solver calls.

\section{Feasibility and bound consistency of conventional GDP results}\label{sec:bound-audit}
Of the 45 conventional GDP runs in Table~\ref{tab:full-gdp}, 21 have flagged bounds (Table~\ref{tab:bound-audit}) relative to the checked solutions in Table~\ref{tab:witnesses}.
These results combine the 27 new \texttt{SCIP}--\texttt{IPOPT} \texttt{GDPopt} runs with the 18 retained direct \texttt{SCIP} runs.
Unrounded values, residuals, and source records are linked in Section~\ref{sec:supp-data}.

For each minimization model, let \(f_w\) be the recorded objective of its best independently checked saved feasible solution, and let \(L\) denote a reported lower bound.
We apply the same rule to all five methods:
\begin{equation}
 \tau_w=10^{-5}\max\{1,|f_w|\},\qquad
 \Delta=L-f_w,\qquad \rho=\Delta/\tau_w,\qquad
 \dagger\ \Longleftrightarrow\ \rho>1.
\end{equation}
Missing lower bounds are left unreported.

Saved vectors are loaded into the original GDP to reevaluate objectives and check active constraints, variable bounds, integrality, disjunctions, and logical validity without optimization.
Violations are computed from the unscaled residuals of the original equations.
All 43 saved feasible vectors satisfy the \(10^{-5}\) threshold; the largest violation is \(9.90\times10^{-7}\).
The two Big-\(M\) runs without accepted feasible solutions remain in the 45-run denominator.
This numerical check establishes neither exact feasibility nor global optimality.

Eight of the 21 discrepancies are only 1.09--3.54 times the tolerance and are sensitive to numerical feasibility conventions.
These conflicting bounds and \nolinkurl{optimal} reports provide no independent global certificate or intrinsic method ranking; their cause in \texttt{GDPopt} or its subsolvers remains unresolved.

\par\endgroup
\begin{table}[htbp]
\centering
\caption{Process GDP models: the 21 conventional GDP runs whose reported lower bound exceeds the recorded objective of the comparison solution of the same model by more than the audit tolerance, from one run per method and model.
The comparison solution of each model is the checked feasible solution of Table~\ref{tab:witnesses}, with recorded objective \(f_w\), and each row gives the reported lower bound of the run, the largest violation at its saved solution, and the excess of the bound over \(f_w\) in multiples of the audit tolerance (Section~\ref{sec:bound-audit}).
Displayed values are rounded, and the flags use unrounded values.}
\label{tab:bound-audit}
\footnotesize
\begin{threeparttable}
\begin{tabular}{@{}lrllrlrlr@{}}
\toprule
Model & \(N\) & Method & ID & Reported \(L\) & Unit & \(v_{\max}\) & C & \(\rho\)\\\midrule
\multirow{17}{*}{\texttt{CSTR}} & \multirow{4}{*}{5} & GLOA & S02 & 3.13020 & L & \(9.85\times10^{-12}\) & C01 & \(2.23\times10^{3}\)\(^{\dagger}\) \\
 & & LOA & S03 & 3.06201 & L & \(2.60\times10^{-11}\) & C01 & 1.29\(^{\dagger}\) \\
 & & RIC & S01 & 3.06201 & L & \(2.60\times10^{-11}\) & C01 & 1.29\(^{\dagger}\) \\
 & & Hull & S05 & 3.06201 & L & \(9.34\times10^{-7}\) & C01 & 1.09\(^{\dagger}\) \\
\addlinespace[2pt]
 & \multirow{4}{*}{10} & GLOA & S07 & 3.44796 & L & \(3.11\times10^{-15}\) & C02 & \(1.93\times10^{4}\)\(^{\dagger}\) \\
 & & LOA & S08 & 2.96605 & L & \(1.62\times10^{-10}\) & C02 & \(2.65\times10^{3}\)\(^{\dagger}\) \\
 & & RIC & S06 & 2.96605 & L & \(1.62\times10^{-10}\) & C02 & \(2.65\times10^{3}\)\(^{\dagger}\) \\
 & & Hull & S10 & 2.88951 & L & \(9.37\times10^{-7}\) & C02 & 3.54\(^{\dagger}\) \\
\addlinespace[2pt]
 & \multirow{3}{*}{15} & GLOA & S12 & 3.56196 & L & \(1.01\times10^{-15}\) & C03 & \(2.59\times10^{4}\)\(^{\dagger}\) \\
 & & LOA & S13 & 2.82939 & L & \(1.33\times10^{-15}\) & C03 & 1.30\(^{\dagger}\) \\
 & & RIC & S11 & 2.82939 & L & \(1.33\times10^{-15}\) & C03 & 1.30\(^{\dagger}\) \\
\addlinespace[2pt]
 & \multirow{3}{*}{20} & GLOA & S17 & 3.61612 & L & \(4.60\times10^{-11}\) & C04 & \(2.92\times10^{4}\)\(^{\dagger}\) \\
 & & LOA & S18 & 2.79941 & L & \(1.58\times10^{-15}\) & C04 & 1.67\(^{\dagger}\) \\
 & & RIC & S16 & 2.79941 & L & \(1.58\times10^{-15}\) & C04 & 1.67\(^{\dagger}\) \\
\addlinespace[2pt]
 & \multirow{2}{*}{25} & GLOA & S22 & 3.64786 & L & \(2.22\times10^{-14}\) & C05 & \(3.12\times10^{4}\)\(^{\dagger}\) \\
 & & RIC & S21 & 2.79125 & L & \(1.38\times10^{-10}\) & C05 & \(3.53\times10^{2}\)\(^{\dagger}\) \\
\addlinespace[2pt]
 & 30 & GLOA & S27 & 3.66872 & L & \(7.11\times10^{-15}\) & C06 & \(3.25\times10^{4}\)\(^{\dagger}\) \\
\midrule
\multirow{2}{*}{\texttt{Column}} & \multirow{2}{*}{--} & GLOA & S32 & 20{,}355.1 & none & \(1.15\times10^{-13}\) & C07 & \(4.65\times10^{3}\)\(^{\dagger}\) \\
 & & Hull & S35 & 20{,}044.2 & none & \(1.00\times10^{-8}\) & C07 & \(3.06\times10^{3}\)\(^{\dagger}\) \\
\midrule
\multirow{2}{*}{\texttt{Dynamic}} & \multirow{2}{*}{--} & GLOA & S42 & -197.140 & none & \(9.85\times10^{-12}\) & C09 & \(7.01\times10^{3}\)\(^{\dagger}\) \\
 & & Hull & S45 & -196.931 & none & \(1.00\times10^{-7}\) & C09 & \(7.11\times10^{3}\)\(^{\dagger}\) \\
\bottomrule
\end{tabular}
\begin{tablenotes}[flushleft]\footnotesize
\item \(f_w\): recorded objective of the comparison solution; \(N\): number of potential reactors of \texttt{CSTR}; --: not applicable; GLOA: global logic-based outer approximation; LOA: logic-based outer approximation; RIC: relaxation with integer cuts; all three are \texttt{GDPopt} methods; Hull: the reformulation solved directly by \texttt{SCIP}.
ID: run identifier; Reported \(L\): lower bound reported by the run; Unit: unit of the objective, as in Table~\ifdefined\manuscriptwithsupplement\ref{tab:process-cases}\else1\fi{} of the main paper, with none for no physical unit; C: identifier of the comparison solution in Table~\ref{tab:witnesses}.
\(v_{\max}\): largest unscaled violation of the original GDP constraints and variable bounds at the saved solution of the run; \(\rho=(L-f_w)/\tau_w\), with the audit tolerance \(\tau_w=10^{-5}\max\{1,|f_w|\}\); \(\dagger\): \(\rho>1\).
\end{tablenotes}
\end{threeparttable}
\end{table}
\FloatBarrier 
\begingroup
\begin{table}[htbp]
\centering
\caption{Process GDP models: comparison solutions of Table~\ref{tab:bound-audit}, one for each model, from the conventional GDP runs with one run per method and model.
The comparison solution of a model is the saved solution with the lowest recorded objective \(f_w\) among the conventional GDP runs whose saved solutions pass the independent feasibility check in the original GDP at \(10^{-5}\), and the method is the conventional GDP method that returned it.
Displayed values are rounded.}
\label{tab:witnesses}
\footnotesize
\begin{threeparttable}
\begin{tabular}{@{}lrllrlr@{}}\toprule
Model & \(N\) & Method & ID & Recorded \(f_w\) & Unit & \(v_{\max}\)\\\midrule
\multirow{6}{*}{\texttt{CSTR}} & 5 & Big-\(M\) & C01 & 3.06198 & L & \(9.40\times10^{-7}\) \\
 & 10 & Big-\(M\) & C02 & 2.88941 & L & \(9.76\times10^{-7}\) \\
 & 15 & Hull & C03 & 2.82936 & L & \(9.47\times10^{-7}\) \\
 & 20 & Hull & C04 & 2.79936 & L & \(9.90\times10^{-7}\) \\
 & 25 & Hull & C05 & 2.78143 & L & \(9.90\times10^{-7}\) \\
 & 30 & Big-\(M\) & C06 & 2.76957 & L & \(9.68\times10^{-9}\) \\
\midrule
\texttt{Column} & -- & RIC & C07 & 19{,}449.9 & none & \(1.01\times10^{-13}\) \\
\midrule
\texttt{Small Batch} & -- & Big-\(M\) & C08 & 167{,}428 & \$ & \(9.00\times10^{-7}\) \\
\midrule
\texttt{Dynamic} & -- & Big-\(M\) & C09 & -212.009 & none & \(1.00\times10^{-7}\) \\
\bottomrule
\end{tabular}
\begin{tablenotes}[flushleft]\footnotesize
\item \(f_w\): recorded objective of the comparison solution; \(N\): number of potential reactors of \texttt{CSTR}; --: not applicable; RIC: relaxation with integer cuts of \texttt{GDPopt}; Big-\(M\), Hull: reformulations solved directly by \texttt{SCIP}.
ID: identifier used in column C of Table~\ref{tab:bound-audit}; Unit: unit of the objective, as in Table~\ifdefined\manuscriptwithsupplement\ref{tab:process-cases}\else1\fi{} of the main paper, with none for no physical unit; \(v_{\max}\): largest unscaled violation of the original GDP constraints and variable bounds at the solution.
\end{tablenotes}
\end{threeparttable}
\end{table}

\paragraph{Additional tightened-feasibility check for \texttt{CSTR15}.}
An \texttt{IPOPT} solve kept the discrete configuration and starting values of an archived LOA solution and targeted an original-GDP residual of \(10^{-8}\)~\citep{rbfoptTestArchive}.
That solution differs from the comparison solution C03 of Table~\ref{tab:witnesses}.
Only exactly satisfied constant rows were omitted from the computational clone; the independent check retained all original constraints.
The returned objective was \(2.82939\), with maximum unscaled violation \(1.11\times10^{-16}\).
The new GLOA lower bound exceeds this objective by \(0.732565\), or approximately \(2.59\times10^4\) times the audit tolerance.
The discrepancy persists under the tighter check, but its cause remains unresolved.

\paragraph{Additional formulation check for \texttt{Dynamic} Hull.}
The saved Big-\(M\) point has objective \(-212.009\), approximately \(15.0776\) below the Hull lower bound \(-196.931\).
It satisfies both the original GDP and a reconstructed Hull formulation at \(10^{-5}\), with maximum violation approximately \(10^{-7}\).
The cause of the discrepancy remains unresolved.

\FloatBarrier
\par\endgroup

\FloatBarrier 
\section{Additional candidate-search and sampling comparisons}\label{sec:supp-additional}

\begin{figure}[htbp]
\centering\includegraphics[width=\linewidth,height=0.6\textheight,keepaspectratio]{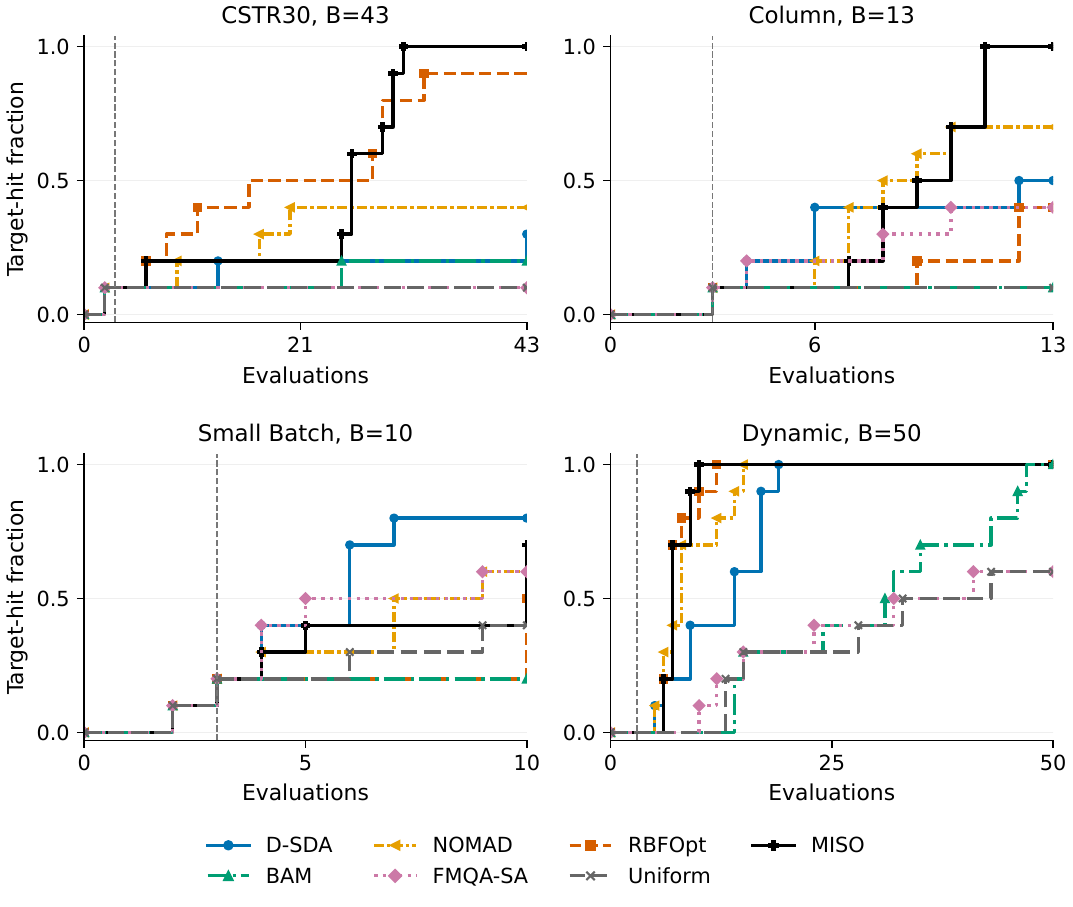}
\caption{Process GDP models: fraction of the ten runs per method that have reached the reference target of Equation~\ifdefined\manuscriptwithsupplement\eqref{eqn:target}\else(8)\fi{} of the main paper, an objective within \(10^{-5}\max\{1,|f_{\mathrm{ref}}|\}\) of the reference objective \(f_{\mathrm{ref}}\) of Table~\ifdefined\manuscriptwithsupplement\ref{tab:process-cases}\else1\fi{} of the main paper, within a given number of evaluations, for four conditions and the six DFO methods and uniform sampling.
Panels: \texttt{CSTR30} at \(B=43\), \texttt{Column} at \(B=13\), \texttt{Small Batch} at \(B=10\), and \texttt{Dynamic} at \(B=50\), where \(B\) is the evaluation budget.
For each method and seed, a \texttt{Dynamic} run with \(B=12\), 15, 20, 30, or 37 reaches the target at the same evaluation as the run with \(B=50\) when that evaluation falls within its budget and misses it otherwise, so the \texttt{Dynamic} panel also gives the profiles of the smaller budgets up to their budgets.
Every evaluation counts, including the three shared initial configurations, unaccepted evaluations, and repeated configurations, and the dashed vertical line at three evaluations marks the end of the shared initial evaluations.
Each curve steps up at the evaluation in which a run first observes an accepted objective that meets the target; runs that miss the target remain in the denominator.
Methods are distinguished by color, line style, and marker.
Uniform: uniform sampling without replacement.
Target hits in all 14 process conditions are listed in Table~\ifdefined\manuscriptwithsupplement\ref{tab:process-targets}\else2\fi{} of the main paper.}
\label{fig:process-evaluations}
\end{figure}

Figure~\ref{fig:process-evaluations} shows the evaluation at which each process run that reaches the target first reaches it, for \texttt{CSTR30}, \texttt{Column}, \texttt{Small Batch}, and \texttt{Dynamic} at its largest budget.

\FloatBarrier

\subsection{Paired differences and sampling controls}
The paired comparisons quantify differences in target-hit rates under shared initial designs, matching case, budget, and seed.
Positive differences favor the first method and are expressed in percentage points.
For Grid, the descriptive percentile bootstrap resamples the 21 intact function vectors, each containing ten paired seed differences, because the ten seeds of a function share its landscape and are not independent draws.
For Process, paired seeds are resampled within each of the 14 conditions; the six related \texttt{Dynamic} budgets are not pooled for inference.
The seed is 20261020 with 2,000 bootstrap replications and deterministic case/seed ordering.
Intervals are exploratory and unadjusted for multiple comparisons; ten-seed boundary cases may yield degenerate intervals, and zero inclusion does not establish equivalence.

An interval entirely above zero favors the first method under this resampling; an interval spanning zero leaves the direction uncertain.
There is no preset cutoff for practical importance.
For example, \texttt{BAM} matches its prefix control at \(B=10\) but exceeds it by 10.95 percentage points at \(B=100\), so initial sampling alone does not account for its final target-hit rate.

\begin{table}[htbp]
\centering
\caption{Grid-GDP and process GDP models: seed-paired differences in the target-hit percentage of selected pairs of methods and their 95\% percentile bootstrap intervals, with 210 paired runs per Grid-GDP budget (21 functions and ten seeds) and ten per process condition.
The difference is the percentage of runs of the first method that reach the reference target of Equation~\ifdefined\manuscriptwithsupplement\eqref{eqn:target}\else(8)\fi{} of the main paper minus that of the second, so positive values favor the first method.
The intervals come from 2,000 bootstrap replications that resample the 21 functions for Grid-GDP and the ten seeds within the stated condition only for the process models, and they are not adjusted for multiple comparisons.
For the process conditions, the seed outcomes count the seeds in which both methods, only the first, only the second, or neither reaches the target.}
\label{tab:supp-paired-differences}
\footnotesize
\setlength{\tabcolsep}{3.5pt}
\begin{threeparttable}
\begin{tabular}{@{}lrlrrrrrrr@{}}\toprule
 & & & & \multicolumn{2}{c}{95\% interval (pp)} & \multicolumn{4}{c}{Seed outcomes}\\\cmidrule(lr){5-6}\cmidrule(l){7-10}
Benchmark & \(B\) & Contrast & \makecell{Difference\\(pp)} & Lower & Upper & Both & \makecell{First\\only} & \makecell{Second\\only} & Neither\\\midrule
\multirow{8}{*}{Grid-GDP} & 10 & \texttt{D-SDA} \(-\) \texttt{RBFOpt} & -13.81 & -24.29 & -6.19 & & & & \\
 & 100 & \texttt{D-SDA} \(-\) \texttt{RBFOpt} & -33.33 & -48.10 & -19.52 & & & & \\
 & 10 & \texttt{RBFOpt} \(-\) Uniform & 13.33 & 5.71 & 24.29 & & & & \\
 & 100 & \texttt{RBFOpt} \(-\) Uniform & 60.00 & 42.86 & 75.71 & & & & \\
 & 10 & \texttt{BAM} \(-\) \texttt{BAM}-prefix control & 0.00 & 0.00 & 0.00 & & & & \\
 & 100 & \texttt{BAM} \(-\) \texttt{BAM}-prefix control & 10.95 & 0.48 & 23.81 & & & & \\
 & 10 & \texttt{FMQA-SA} \(-\) Uniform & 0.95 & -0.95 & 2.86 & & & & \\
 & 100 & \texttt{FMQA-SA} \(-\) Uniform & 40.00 & 28.57 & 52.38 & & & & \\
\midrule
\texttt{Small Batch} & 10 & \texttt{FMQA-SA} \(-\) Uniform & 20 & -30 & 60 & 2 & 4 & 2 & 2 \\
\midrule
\texttt{Dynamic} & 12 & \texttt{RBFOpt} \(-\) \texttt{D-SDA} & 60 & 30 & 90 & 4 & 6 & 0 & 0 \\
\bottomrule
\end{tabular}
\begin{tablenotes}[flushleft]\footnotesize
\item \(B\): evaluation budget; pp: percentage points; Contrast: first method \(-\) second method; Uniform: uniform sampling without replacement; \texttt{BAM}-prefix control: the ten predetermined sampling points of \texttt{BAM} after the three shared initial configurations, followed by uniform sampling without replacement.
\end{tablenotes}
\end{threeparttable}
\end{table}
\FloatBarrier

\subsection{Exact paired tests for process conditions}
Table~\ref{tab:supp-paired} lists the seed-paired target outcomes of \texttt{RBFOpt} and \texttt{D-SDA} in every process condition and of \texttt{BAM} and \texttt{RBFOpt} on \texttt{CSTR5}.
The \(p\)-values come from the two-sided exact McNemar test, which uses only the pairs in which exactly one method reaches the target.
They are descriptive and not adjusted for multiple comparisons.
In the 60 paired \texttt{CSTR} runs, \texttt{RBFOpt} reaches the target in 30 of the 31 discordant pairs.

\begin{table}[htbp]
\centering
\caption{Process GDP models: seed-paired target outcomes of \texttt{RBFOpt} and \texttt{D-SDA} and of \texttt{BAM} and \texttt{RBFOpt}, with ten seed pairs per condition.
The counts give the seeds in which both methods, only the first, only the second, or neither reaches the reference target of Equation~\ifdefined\manuscriptwithsupplement\eqref{eqn:target}\else(8)\fi{} of the main paper.
The \(p\)-value uses only the seeds in which exactly one method reaches the target, and the \(p\)-values are descriptive and not adjusted for multiple comparisons.
The pooled \texttt{CSTR} row combines the six sizes, with 60 seed pairs.}
\label{tab:supp-paired}
\footnotesize
\begin{threeparttable}
\begin{tabular}{@{}lrrrrrrr@{}}\toprule
Model & \(N\) & \(B\) & Both & First only & Second only & Neither & \(p\)\\\midrule
\multicolumn{8}{@{}l}{First \texttt{RBFOpt}, second \texttt{D-SDA}}\\\midrule
\multirow{7}{*}{\texttt{CSTR}} & 5 & 8 & 2 & 6 & 0 & 2 & 0.0312 \\
 & 10 & 15 & 2 & 4 & 1 & 3 & 0.375 \\
 & 15 & 22 & 2 & 5 & 0 & 3 & 0.0625 \\
 & 20 & 29 & 3 & 4 & 0 & 3 & 0.125 \\
 & 25 & 36 & 3 & 5 & 0 & 2 & 0.0625 \\
 & 30 & 43 & 3 & 6 & 0 & 1 & 0.0312 \\
 & all & -- & 15 & 30 & 1 & 14 & \(2.98\times10^{-8}\) \\
\midrule
\texttt{Column} & -- & 13 & 3 & 1 & 2 & 4 & 1 \\
\midrule
\texttt{Small Batch} & -- & 10 & 4 & 1 & 4 & 1 & 0.375 \\
\midrule
\multirow{6}{*}{\texttt{Dynamic}} & -- & 12 & 4 & 6 & 0 & 0 & 0.0312 \\
 & & 15 & 6 & 4 & 0 & 0 & 0.125 \\
 & & 20 & 10 & 0 & 0 & 0 & 1 \\
 & & 30 & 10 & 0 & 0 & 0 & 1 \\
 & & 37 & 10 & 0 & 0 & 0 & 1 \\
 & & 50 & 10 & 0 & 0 & 0 & 1 \\
\midrule
\multicolumn{8}{@{}l}{First \texttt{BAM}, second \texttt{RBFOpt}}\\\midrule
\texttt{CSTR} & 5 & 8 & 8 & 2 & 0 & 0 & 0.500 \\
\bottomrule
\end{tabular}
\begin{tablenotes}[flushleft]\footnotesize
\item \(N\): number of potential reactors of \texttt{CSTR}; \(B\): evaluation budget; all: the six sizes pooled; --: not applicable in the \(N\) column and several budgets in the \(B\) column; \(p\): two-sided exact McNemar \(p\)-value.
\end{tablenotes}
\end{threeparttable}
\end{table}
\FloatBarrier

\subsection{Grid suite exclusions}
Excluding Xin-She Yang~3 and Bohachevsky~2 preserves the main descriptive pattern.
For the sampled Bohachevsky~2 lattice, \(\cos(3\pi a_i)=1\) at every coordinate, leaving an objective nearly redundant with sampled Sum of Squares after normalization.
Without both functions, \(B=100\) rates are 84.7\% for \texttt{MISO}, 79.5\% for \texttt{RBFOpt}, 66.8\% for \texttt{BAM}, 57.9\% for \texttt{FMQA-SA}, 55.8\% for \texttt{NOMAD}, 54.7\% for the \texttt{BAM}-prefix control, 42.6\% for \texttt{D-SDA}, and 17.9\% for uniform.
At \(B=10\), \texttt{BAM} and its prefix control each retain 47.9\%, compared with 15.8\% for \texttt{RBFOpt}, 14.7\% for \texttt{MISO}, and 12.6\% for \texttt{NOMAD}.


\FloatBarrier
\section{Seed spread, tolerance sensitivity, and search overhead}\label{sec:supp-spread-overhead}
Table~\ref{tab:supp-spread} compares exact target attainment with relaxed relative tolerances for the 980 primary process runs.
The relative gap is \((f_{\mathrm{best}}-f_{\mathrm{ref}})/\max\{1,|f_{\mathrm{ref}}|\}\), the normalization of the target rule in Equation~\ifdefined\manuscriptwithsupplement\eqref{eqn:target}\else(8)\fi{} of the main paper.
The exact target-hit rows apply the target tolerance of \(10^{-5}\).
Distinct finals is the largest number of different final objectives among the ten seeds of any condition in the family, and the smallest and largest miss-gap rows cover every run that misses the exact target.
Table~\ref{tab:supp-overhead} compares the median configuration-evaluation cost with the median search time per evaluation.
Both are medians over runs of per-run averages, namely the summed point cost or the stored search time divided by the number of evaluations, over all seven methods for the evaluation cost.
The 21 \texttt{NOMAD} Grid runs without a recorded search time are left out of its search median.

\begin{table}[htbp]
\centering
\caption{Process GDP models: exact and relaxed target attainment and the spread of final objectives of the six DFO methods and uniform sampling, with ten runs per method and condition, pooled over the six sizes of \texttt{CSTR} and the six budgets of \texttt{Dynamic}.
An exact target hit is a run that reaches the reference target of Equation~\ifdefined\manuscriptwithsupplement\eqref{eqn:target}\else(8)\fi{} of the main paper, with relative tolerance \(10^{-5}\), and the relaxed rows count the runs whose relative gap \((f_{\mathrm{best}}-f_{\mathrm{ref}})/\max\{1,|f_{\mathrm{ref}}|\}\) is at most the stated value.
The miss gaps are the smallest and largest relative gaps over the runs that miss the exact target, and the distinct-finals row gives the largest number of different final objectives among the ten seeds of any condition of the model.}
\label{tab:supp-spread}
\footnotesize
\setlength{\tabcolsep}{3.5pt}
\begin{threeparttable}
\begin{tabular}{@{}llrrrrrrr@{}}\toprule
Model & & \texttt{D-SDA} & \texttt{NOMAD} & \texttt{RBFOpt} & \texttt{MISO} & \texttt{BAM} & \texttt{FMQA-SA} & Uniform\\\midrule
\multirow{7}{*}{\texttt{CSTR}} & Exact target hits (of 60) & 16 & 21 & 45 & 43 & 18 & 9 & 11 \\
 & Runs with gap \(\leq0.1\%\) & 16 & 21 & 45 & 43 & 18 & 9 & 11 \\
 & Runs with gap \(\leq1\%\) & 19 & 21 & 45 & 43 & 60 & 23 & 25 \\
 & Runs with gap \(\leq5\%\) & 27 & 27 & 47 & 49 & 60 & 41 & 48 \\
 & Distinct finals & 5 & 3 & 3 & 5 & 3 & 10 & 10 \\
 & Smallest miss gap (\%) & 0.224 & 2.23 & 2.23 & 2.23 & 0.224 & 0.169 & 0.224 \\
 & Largest miss gap (\%) & 13.0 & 13.0 & 12.5 & 10.6 & 0.692 & 17.2 & 14.3 \\
\midrule
\multirow{7}{*}{\texttt{Column}} & Exact target hits (of 10) & 5 & 7 & 4 & 10 & 1 & 4 & 1 \\
 & Runs with gap \(\leq0.1\%\) & 5 & 7 & 4 & 10 & 1 & 4 & 1 \\
 & Runs with gap \(\leq1\%\) & 8 & 9 & 5 & 10 & 1 & 4 & 2 \\
 & Runs with gap \(\leq5\%\) & 10 & 10 & 9 & 10 & 10 & 10 & 9 \\
 & Distinct finals & 4 & 3 & 6 & 1 & 2 & 4 & 8 \\
 & Smallest miss gap (\%) & 0.536 & 0.536 & 0.536 & -- & 1.05 & 1.05 & 0.536 \\
 & Largest miss gap (\%) & 1.14 & 1.14 & 5.92 & -- & 1.05 & 2.76 & 5.22 \\
\midrule
\multirow{7}{*}{\texttt{Small Batch}} & Exact target hits (of 10) & 8 & 6 & 5 & 7 & 2 & 6 & 4 \\
 & Runs with gap \(\leq0.1\%\) & 8 & 6 & 5 & 7 & 2 & 6 & 4 \\
 & Runs with gap \(\leq1\%\) & 8 & 6 & 5 & 7 & 2 & 6 & 4 \\
 & Runs with gap \(\leq5\%\) & 8 & 6 & 5 & 7 & 2 & 6 & 4 \\
 & Distinct finals & 2 & 2 & 3 & 2 & 3 & 5 & 4 \\
 & Smallest miss gap (\%) & 22.2 & 6.64 & 6.64 & 6.64 & 6.64 & 6.64 & 6.64 \\
 & Largest miss gap (\%) & 22.2 & 6.64 & 22.2 & 6.64 & 8.23 & 22.2 & 22.2 \\
\midrule
\multirow{7}{*}{\texttt{Dynamic}} & Exact target hits (of 60) & 50 & 58 & 60 & 60 & 27 & 23 & 21 \\
 & Runs with gap \(\leq0.1\%\) & 50 & 58 & 60 & 60 & 27 & 23 & 21 \\
 & Runs with gap \(\leq1\%\) & 50 & 58 & 60 & 60 & 27 & 23 & 21 \\
 & Runs with gap \(\leq5\%\) & 54 & 60 & 60 & 60 & 46 & 42 & 38 \\
 & Distinct finals & 5 & 3 & 1 & 1 & 4 & 7 & 6 \\
 & Smallest miss gap (\%) & 2.07 & 2.07 & -- & -- & 2.07 & 2.07 & 2.07 \\
 & Largest miss gap (\%) & 36.3 & 2.66 & -- & -- & 7.01 & 36.3 & 30.7 \\
\bottomrule
\end{tabular}
\begin{tablenotes}[flushleft]\footnotesize
\item \(f_{\mathrm{best}}\): best accepted objective of the run; \(f_{\mathrm{ref}}\): reference objective of Table~\ifdefined\manuscriptwithsupplement\ref{tab:process-cases}\else1\fi{} of the main paper; Uniform: uniform sampling without replacement; --: no run misses the exact target.
\end{tablenotes}
\end{threeparttable}
\end{table}

\begin{table}[htbp]
\centering
\caption{Grid-GDP and process GDP models: median evaluation cost and median search time per evaluation of the six DFO methods and uniform sampling, with 1,050 Grid-GDP runs per method (21 functions, five budgets, and ten seeds) and ten runs per method and process condition.
The evaluation cost of a run is the summed recorded cost of its evaluated configurations divided by its number of evaluations, and its median pools the runs of all seven methods.
The search time per evaluation of a run is its search time, the method time excluding the evaluation cost, divided by its number of evaluations.
Medians are over runs, and the 21 \texttt{NOMAD} Grid-GDP runs without a recorded search time are left out of the \texttt{NOMAD} median.}
\label{tab:supp-overhead}
\footnotesize
\setlength{\tabcolsep}{3.5pt}
\begin{threeparttable}
\begin{tabular}{@{}lrrrrrrrr@{}}\toprule
 & & \multicolumn{7}{c}{Median search time per evaluation (ms)}\\\cmidrule(l){3-9}
Benchmark & \makecell{Median evaluation\\cost (s)} & \texttt{D-SDA} & \texttt{NOMAD} & \texttt{RBFOpt} & \texttt{MISO} & \texttt{BAM} & \texttt{FMQA-SA} & Uniform\\\midrule
Grid-GDP & 0.390 & 0.0777 & 50.5 & 31.1 & 231 & 3.77 & 604 & 0.0336 \\
\midrule
\texttt{CSTR} & 1.46 & 1.24 & 72.2 & 24.2 & 287 & 9.01 & 572 & 0.517 \\
\texttt{Column} & 2.50 & 0.941 & 3.25 & 17.4 & 51.7 & 1.36 & 527 & 0.551 \\
\texttt{Small Batch} & 0.422 & 1.02 & 4.75 & 17.0 & 66.8 & 1.69 & 487 & 0.536 \\
\texttt{Dynamic} & 0.836 & 0.893 & 96.9 & 98.4 & 1{,}860 & 9.34 & 567 & 0.471 \\
\bottomrule
\end{tabular}
\begin{tablenotes}[flushleft]\footnotesize
\item Uniform: uniform sampling without replacement.
\end{tablenotes}
\end{threeparttable}
\end{table}

Every NLP subproblem solved by \texttt{GDPopt} in the process comparisons fixes a different configuration, as recorded in its execution logs before each subproblem.
RIC and LOA solve one subproblem for each of the \(N(N+1)/2\) logically admissible configurations of \texttt{CSTR5} to \texttt{CSTR20} and for each of the 37 admissible \texttt{Dynamic} configurations, and RIC also covers all 325 admissible configurations of \texttt{CSTR25}.
\FloatBarrier

\section{\texttt{NOMAD} and \texttt{MISO} runs}\label{sec:supp-nomad-miso}
The \texttt{NOMAD} and \texttt{MISO} runs use the stored responses, shared initial configurations, seeds, budgets, penalties, and reference targets of the other five methods.
Neither method receives the target as a stopping criterion.
They were executed after the other methods, on a shared server with six workers, each restricted to one logical CPU on a different physical core and to one OpenMP and BLAS thread; the cores were not reserved exclusively.
The server used \texttt{NumPy} 2.4.6 and \texttt{SciPy} 1.17.1; the other versions are given in Section~\ifdefined\manuscriptwithsupplement\ref{sec:protocol}\else3.6\fi{} of the main paper.
As for the other methods, completion times add the measured search time of each run to the recorded cost of every evaluated configuration, so they do not depend on the speed of the server's NLP solves.

\paragraph{Termination.}
All 1,190 \texttt{MISO} runs use their full budget.
\texttt{NOMAD} stops before the budget in 77 of the 140 process runs and in 697 of the 1,050 Grid runs, including all 210 runs at \(B=100\), where it stops after a median of 31.5 evaluations.
Its Grid target fraction is therefore 60.0\% at both \(B=70\) and \(B=100\).

\paragraph{\texttt{NOMAD} crashes and reruns.}
Of the 1,190 \texttt{NOMAD} runs, 1,046 completed on the first attempt.
The other 144 ended with a segmentation fault after their last new evaluation.
Of these, 110 completed on the first rerun and 13 on the second, each after evaluating the same configurations in the same order as the crashed attempt, and 21 crashed on all three attempts after identical evaluations.
All 21 are Grid runs at \(B=100\), and 16 of them had reached the target before the crash.
Scoring them on the evaluations completed before the crash gives \texttt{NOMAD} 526 of 1,050 Grid targets and 92 of 140 process targets; counting all 21 as misses gives 510 Grid targets, or 52.4\% instead of 60.0\% at \(B=100\).
These 21 runs have no recorded search time, so the \texttt{NOMAD} timing summaries at Grid \(B=100\) use the other 189 runs.
For a completed rerun, the search time is that of the completed attempt; the time spent in crashed attempts is not recorded.
No \texttt{MISO} run failed.

\paragraph{Exploratory paired tests.}
Table~\ref{tab:supp-nomad-miso-tests} compares \texttt{NOMAD} and \texttt{MISO} with \texttt{D-SDA}, \texttt{RBFOpt}, and each other in every process condition, a family of 70 exact McNemar tests defined after the runs were completed.
With a Holm correction over the 70 tests, the smallest adjusted \(p\)-value is 0.547, for \texttt{MISO} against \texttt{NOMAD} on \texttt{CSTR20}, so no contrast is significant at the 0.05 level.
Ten seeds give each test little power, and a nonsignificant test does not establish equivalence.
The six \texttt{Dynamic} budgets reuse the same seeds and are not independent replications.
The last row pools the six \texttt{CSTR} sizes for the descriptive comparisons in the main text; it is not part of the corrected family.

\begin{table}[htbp]
\centering
\caption{Process GDP models: seed-paired target outcomes of \texttt{NOMAD} and \texttt{MISO} against \texttt{D-SDA}, \texttt{RBFOpt}, and each other, with ten seed pairs per condition.
The discordant seeds are those in which only one method of the column reaches the reference target of Equation~\ifdefined\manuscriptwithsupplement\eqref{eqn:target}\else(8)\fi{} of the main paper, and the \(p\)-values are not adjusted for multiple comparisons.
With a Holm correction over the 70 per-condition tests, the smallest adjusted \(p\)-value is 0.547.
The pooled \texttt{CSTR} row combines the six sizes, with 60 seed pairs, and is not part of the 70 per-condition tests.}
\label{tab:supp-nomad-miso-tests}
\footnotesize
\setlength{\tabcolsep}{3.5pt}
\begin{threeparttable}
\begin{tabular}{@{}lrrlrrrrr@{}}\toprule
Model & \(N\) & \(B\) & & \makecell{\texttt{NOMAD}\\vs \texttt{D-SDA}} & \makecell{\texttt{NOMAD}\\vs \texttt{RBFOpt}} & \makecell{\texttt{MISO}\\vs \texttt{D-SDA}} & \makecell{\texttt{MISO}\\vs \texttt{RBFOpt}} & \makecell{\texttt{MISO}\\vs \texttt{NOMAD}}\\\midrule
\multirow{12}{*}{\texttt{CSTR}} & \multirow{2}{*}{5} & \multirow{2}{*}{8} & Discordant seeds & 2:0 & 1:5 & 2:0 & 0:4 & 1:1 \\
 & & & \(p\) & 0.500 & 0.219 & 0.500 & 0.125 & 1 \\
\addlinespace[2pt]
 & \multirow{2}{*}{10} & \multirow{2}{*}{15} & Discordant seeds & 0:0 & 1:4 & 0:1 & 0:4 & 0:1 \\
 & & & \(p\) & 1 & 0.375 & 1 & 0.125 & 1 \\
\addlinespace[2pt]
 & \multirow{2}{*}{15} & \multirow{2}{*}{22} & Discordant seeds & 0:0 & 0:5 & 5:0 & 1:1 & 5:0 \\
 & & & \(p\) & 1 & 0.0625 & 0.0625 & 1 & 0.0625 \\
\addlinespace[2pt]
 & \multirow{2}{*}{20} & \multirow{2}{*}{29} & Discordant seeds & 0:1 & 0:5 & 7:0 & 3:0 & 8:0 \\
 & & & \(p\) & 1 & 0.0625 & 0.0156 & 0.250 & 0.00781 \\
\addlinespace[2pt]
 & \multirow{2}{*}{25} & \multirow{2}{*}{36} & Discordant seeds & 4:1 & 2:4 & 7:0 & 2:0 & 4:0 \\
 & & & \(p\) & 0.375 & 0.688 & 0.0156 & 0.500 & 0.125 \\
\addlinespace[2pt]
 & \multirow{2}{*}{30} & \multirow{2}{*}{43} & Discordant seeds & 1:0 & 1:6 & 7:0 & 1:0 & 6:0 \\
 & & & \(p\) & 1 & 0.125 & 0.0156 & 1 & 0.0312 \\
\midrule
\multirow{2}{*}{\texttt{Column}} & \multirow{2}{*}{--} & \multirow{2}{*}{13} & Discordant seeds & 3:1 & 4:1 & 5:0 & 6:0 & 3:0 \\
 & & & \(p\) & 0.625 & 0.375 & 0.0625 & 0.0312 & 0.250 \\
\midrule
\multirow{2}{*}{\texttt{Small Batch}} & \multirow{2}{*}{--} & \multirow{2}{*}{10} & Discordant seeds & 1:3 & 2:1 & 1:2 & 2:0 & 3:2 \\
 & & & \(p\) & 0.625 & 1 & 1 & 0.500 & 1 \\
\midrule
\multirow{12}{*}{\texttt{Dynamic}} & \multirow{2}{*}{--} & \multirow{2}{*}{12} & Discordant seeds & 4:0 & 0:2 & 6:0 & 0:0 & 2:0 \\
 & & & \(p\) & 0.125 & 0.500 & 0.0312 & 1 & 0.500 \\
\addlinespace[2pt]
 &  & \multirow{2}{*}{15} & Discordant seeds & 4:0 & 0:0 & 4:0 & 0:0 & 0:0 \\
 & & & \(p\) & 0.125 & 1 & 0.125 & 1 & 1 \\
\addlinespace[2pt]
 &  & \multirow{2}{*}{20} & Discordant seeds & 0:0 & 0:0 & 0:0 & 0:0 & 0:0 \\
 & & & \(p\) & 1 & 1 & 1 & 1 & 1 \\
\addlinespace[2pt]
 &  & \multirow{2}{*}{30} & Discordant seeds & 0:0 & 0:0 & 0:0 & 0:0 & 0:0 \\
 & & & \(p\) & 1 & 1 & 1 & 1 & 1 \\
\addlinespace[2pt]
 &  & \multirow{2}{*}{37} & Discordant seeds & 0:0 & 0:0 & 0:0 & 0:0 & 0:0 \\
 & & & \(p\) & 1 & 1 & 1 & 1 & 1 \\
\addlinespace[2pt]
 &  & \multirow{2}{*}{50} & Discordant seeds & 0:0 & 0:0 & 0:0 & 0:0 & 0:0 \\
 & & & \(p\) & 1 & 1 & 1 & 1 & 1 \\
\midrule
\multirow{2}{*}{\texttt{CSTR}} & \multirow{2}{*}{all} & \multirow{2}{*}{--} & Discordant seeds & 7:2 & 5:29 & 28:1 & 7:9 & 24:2 \\
 & & & \(p\) & 0.180 & \(3.86\times10^{-5}\) & \(1.12\times10^{-7}\) & 0.804 & \(1.05\times10^{-5}\) \\
\bottomrule
\end{tabular}
\begin{tablenotes}[flushleft]\footnotesize
\item \(a{:}b\): discordant seeds in which only the first or only the second method of the column reaches the target; \(p\): unadjusted two-sided exact McNemar \(p\)-value; \(N\): number of potential reactors of \texttt{CSTR}; \(B\): evaluation budget; all: the six sizes pooled; --: not applicable in the \(N\) column and several budgets in the \(B\) column.
\end{tablenotes}
\end{threeparttable}
\end{table}
\FloatBarrier

\section{Process \texttt{D-SDA} restarts}\label{sec:supp-dsda-restarts}
The process \texttt{D-SDA} of the main comparison descends once from the best of the three shared initial configurations (S).
Variant R applies the Grid rule and descends from all three shared initial configurations, best first, while budget remains.
Variant RR continues after R with descents from uniformly drawn unvisited configurations until the budget is exhausted, and it uses the full budget in all 140 runs.
No variant receives the reference target.
All variants read the stored responses and, like the main comparison, never evaluate a configuration twice.
Table~\ref{tab:supp-dsda-restarts} shows that R and RR evaluate 1,055 and 1,249 more configurations than S and add one target, on seed 4 of \texttt{CSTR30}, so target attainment rises from 79 to 80 of the 140 runs.
Both restart variants improve the final objective in the same 14 \texttt{CSTR} runs.
In each of them, S ends at \((5,1)\), 8.3 to 13.0\% above the reference, and the restarts end at the target in the one new hit and 0.4 to 8.6\% above it in the other 13.
On \texttt{CSTR10} to \texttt{CSTR30}, 31 of the 36 single-descent misses end at \((5,1)\).

\begin{table}[htbp]
\centering
\caption{Process GDP models: target hits, full-budget runs, median evaluations, and median full-run completion times of three \texttt{D-SDA} restart variants, with ten runs per variant and condition.
A target hit is a run that reaches the reference target of Equation~\ifdefined\manuscriptwithsupplement\eqref{eqn:target}\else(8)\fi{} of the main paper, and a full-budget run uses all \(B\) evaluations, as every RR run does.
The completion time of a run adds its measured search time to the recorded cost of every evaluated configuration, and medians are over the ten runs.}
\label{tab:supp-dsda-restarts}
\footnotesize
\setlength{\tabcolsep}{3.5pt}
\begin{threeparttable}
\begin{tabular}{@{}lrrrrrrrrrrrrr@{}}\toprule
 & & & \multicolumn{3}{c}{\makecell{Target hits\\(of 10)}} & \multicolumn{2}{c}{\makecell{Full\\budget\\(of 10)}} & \multicolumn{3}{c}{\makecell{Median\\evaluations}} & \multicolumn{3}{c}{\makecell{Median completion\\time (s)}}\\\cmidrule(lr){4-6}\cmidrule(lr){7-8}\cmidrule(lr){9-11}\cmidrule(l){12-14}
Model & \(N\) & \(B\) & S & R & RR & S & R & S & R & RR & S & R & RR\\\midrule
\multirow{6}{*}{\texttt{CSTR}} & 5 & 8 & 2 & 2 & 2 & 8 & 10 & 8 & 8 & 8 & 10.67 & 10.71 & 10.71 \\
 & 10 & 15 & 3 & 3 & 3 & 3 & 10 & 11.5 & 15 & 15 & 17.61 & 21.57 & 21.57 \\
 & 15 & 22 & 2 & 2 & 2 & 2 & 10 & 15.5 & 22 & 22 & 25.68 & 35.19 & 35.19 \\
 & 20 & 29 & 3 & 3 & 3 & 0 & 10 & 19 & 29 & 29 & 32.31 & 47.49 & 47.49 \\
 & 25 & 36 & 3 & 3 & 3 & 2 & 9 & 23.5 & 36 & 36 & 45.49 & 60.78 & 60.78 \\
 & 30 & 43 & 3 & 4 & 4 & 1 & 8 & 24.5 & 43 & 43 & 47.17 & 78.88 & 78.88 \\
\midrule
\texttt{Column} & -- & 13 & 5 & 5 & 5 & 9 & 10 & 13 & 13 & 13 & 35.44 & 35.69 & 35.69 \\
\midrule
\texttt{Small Batch} & -- & 10 & 8 & 8 & 8 & 10 & 10 & 10 & 10 & 10 & 4.23 & 4.23 & 4.23 \\
\midrule
\multirow{6}{*}{\texttt{Dynamic}} & -- & 12 & 4 & 4 & 4 & 9 & 10 & 12 & 12 & 12 & 12.22 & 12.25 & 12.25 \\
 &  & 15 & 6 & 6 & 6 & 6 & 10 & 15 & 15 & 15 & 15.32 & 15.76 & 15.76 \\
 &  & 20 & 10 & 10 & 10 & 3 & 10 & 18 & 20 & 20 & 20.02 & 20.76 & 20.76 \\
 &  & 30 & 10 & 10 & 10 & 0 & 10 & 18 & 30 & 30 & 20.02 & 27.13 & 27.13 \\
 &  & 37 & 10 & 10 & 10 & 0 & 2 & 18 & 34.5 & 37 & 20.02 & 28.01 & 30.05 \\
 &  & 50 & 10 & 10 & 10 & 0 & 0 & 18 & 34.5 & 50 & 20.02 & 28.01 & 37.00 \\
\bottomrule
\end{tabular}
\begin{tablenotes}[flushleft]\footnotesize
\item S: single descent from the best shared initial configuration; R: descents from the three shared initial configurations, best first, while budget remains; RR: R followed by descents from uniformly drawn unvisited configurations until the budget is exhausted.
\(N\): number of potential reactors of \texttt{CSTR}; \(B\): evaluation budget; --: not applicable.
\end{tablenotes}
\end{threeparttable}
\end{table}
\FloatBarrier

\clearpage
\section{\texttt{D-SDA} neighborhood}\label{sec:supp-dsda-neighborhood}
The process, Grid, and \texttt{Lomustine} \texttt{D-SDA} runs use the unit \(L_\infty\) neighborhood, whose \(3^d-1\) neighbors number 8 for \(d=2\) external variables and 26 for \(d=3\).
The \texttt{Catalytic Column}, with \(d=6\) and 728 such neighbors, uses the separable neighborhood, which changes one external variable by \(\pm1\) and has at most \(2d\) neighbors.
To test the effect of this choice, we repeated every \texttt{D-SDA} run on stored responses with the separable neighborhood and all other settings unchanged: the 1,050 Grid runs, the 140 single-descent process runs, and the 180 \texttt{Lomustine} runs in the six option orders of Section~\ref{sec:supp-lomustine}.
Table~\ref{tab:supp-dsda-neighborhood} pairs each run with its \(L_\infty\) counterpart.

The separable neighborhood reaches more targets when the budget is small compared with the \(L_\infty\) neighborhood: on Grid at \(B=10\) (22 against 13 of 210 runs, with 9 and 0 discordant runs, exact McNemar \(p=0.004\)), on \texttt{Dynamic} at \(B=12\) and 15 (16 against 10 of 20, \(p=0.031\)), and on \texttt{Lomustine} over the six orders (99 against 90 of 180, \(p=0.022\)).
Its descents also stop earlier: from \(B=50\), the Grid runs end after a median of 25 evaluations in place of 37, and they reach 88 targets against 101 at \(B=70\) and 100 (\(p=0.011\)).
On \texttt{CSTR}, the separable neighborhood reaches 7 of the 60 targets against 16 (\(p=0.004\)).
Eighteen of its 53 misses stop on the edge \(\xi_1=\xi_2\) of the admissible region, where every one-coordinate move leaves the region or worsens the objective, whereas a diagonal move continues along the edge to the target at \((N,N)\).
Over the 140 process runs, the two neighborhoods reach 79 and 78 targets.
Neither neighborhood dominates: the \(L_\infty\) neighborhood suits two external variables with budgets that cover its eight neighbors, and the separable neighborhood suits more external variables or smaller budgets.

\begin{table}[htbp]
\centering
\caption{Grid-GDP, process GDP models, and \texttt{Lomustine}: target hits and median evaluations of \texttt{D-SDA} with the unit \(L_\infty\) neighborhood and with the separable neighborhood, with 210 runs of each neighborhood per Grid-GDP row (21 functions and ten seeds), ten per process row, and 60 per \texttt{Lomustine} row (six option orders and ten seeds).
Both neighborhoods run on stored responses with all other settings unchanged, the process runs use the single descent, and each separable run is paired with its \(L_\infty\) counterpart.
A target hit is a run that reaches the reference target of Equation~\ifdefined\manuscriptwithsupplement\eqref{eqn:target}\else(8)\fi{} of the main paper, and the median evaluations are medians of the number of evaluations per run.}
\label{tab:supp-dsda-neighborhood}
\footnotesize
\begin{threeparttable}
\begin{tabular}{@{}lrrrrrrrr@{}}\toprule
 & & & \multicolumn{4}{c}{Target hits} & \multicolumn{2}{c}{Median evaluations}\\\cmidrule(lr){4-7}\cmidrule(l){8-9}
Benchmark & \(N\) & \(B\) & \(L_\infty\) & Separable & Separable only & \(L_\infty\) only & \(L_\infty\) & Separable \\\midrule
\multirow{5}{*}{Grid-GDP} & -- & 10 & 13 & \textbf{22} & 9 & 0 & 10 & 10 \\
 &  & 30 & 78 & \textbf{79} & 14 & 13 & 30 & 25 \\
 &  & 50 & \textbf{91} & 87 & 10 & 14 & 37 & 25 \\
 &  & 70 & \textbf{101} & 88 & 5 & 18 & 37 & 25 \\
 &  & 100 & \textbf{101} & 88 & 5 & 18 & 37 & 25 \\
\midrule
\multirow{6}{*}{\texttt{CSTR}} & 5 & 8 & 2 & 2 & 0 & 0 & 8 & 8 \\
 & 10 & 15 & \textbf{3} & 1 & 0 & 2 & 11.5 & 8 \\
 & 15 & 22 & \textbf{2} & 1 & 0 & 1 & 15.5 & 10.5 \\
 & 20 & 29 & \textbf{3} & 1 & 0 & 2 & 19 & 14 \\
 & 25 & 36 & \textbf{3} & 1 & 0 & 2 & 23.5 & 15.5 \\
 & 30 & 43 & \textbf{3} & 1 & 0 & 2 & 24.5 & 17 \\
\midrule
\texttt{Column} & -- & 13 & 5 & \textbf{6} & 3 & 2 & 13 & 12 \\
\midrule
\texttt{Small Batch} & -- & 10 & 8 & \textbf{9} & 2 & 1 & 10 & 10 \\
\midrule
\multirow{6}{*}{\texttt{Dynamic}} & -- & 12 & 4 & \textbf{6} & 2 & 0 & 12 & 12 \\
 &  & 15 & 6 & \textbf{10} & 4 & 0 & 15 & 13 \\
 &  & 20 & 10 & 10 & 0 & 0 & 18 & 13 \\
 &  & 30 & 10 & 10 & 0 & 0 & 18 & 13 \\
 &  & 37 & 10 & 10 & 0 & 0 & 18 & 13 \\
 &  & 50 & 10 & 10 & 0 & 0 & 18 & 13 \\
\midrule
\multirow{3}{*}{\texttt{Lomustine}} & -- & 10 & 23 & \textbf{27} & 6 & 2 & 10 & 10 \\
 &  & 15 & 32 & \textbf{36} & 4 & 0 & 15 & 11 \\
 &  & 20 & 35 & \textbf{36} & 1 & 0 & 19 & 11 \\
\bottomrule
\end{tabular}
\begin{tablenotes}[flushleft]\footnotesize
\item \(L_\infty\): unit neighborhood, which changes every external variable by at most one; Separable: neighborhood that changes one external variable by \(\pm1\); Separable only, \(L_\infty\) only: paired runs in which only that neighborhood reaches the target.
\(N\): number of potential reactors of \texttt{CSTR}; \(B\): evaluation budget; --: not applicable.
Bold: larger number of target hits of the two neighborhoods in each row; equal counts are not bold.
\end{tablenotes}
\end{threeparttable}
\end{table}
\FloatBarrier

\section{Joint comparison with conventional GDP methods}
\label{sec:supp-joint-comparison}
\begin{figure}[htbp]
\centering
\includegraphics[width=\linewidth]{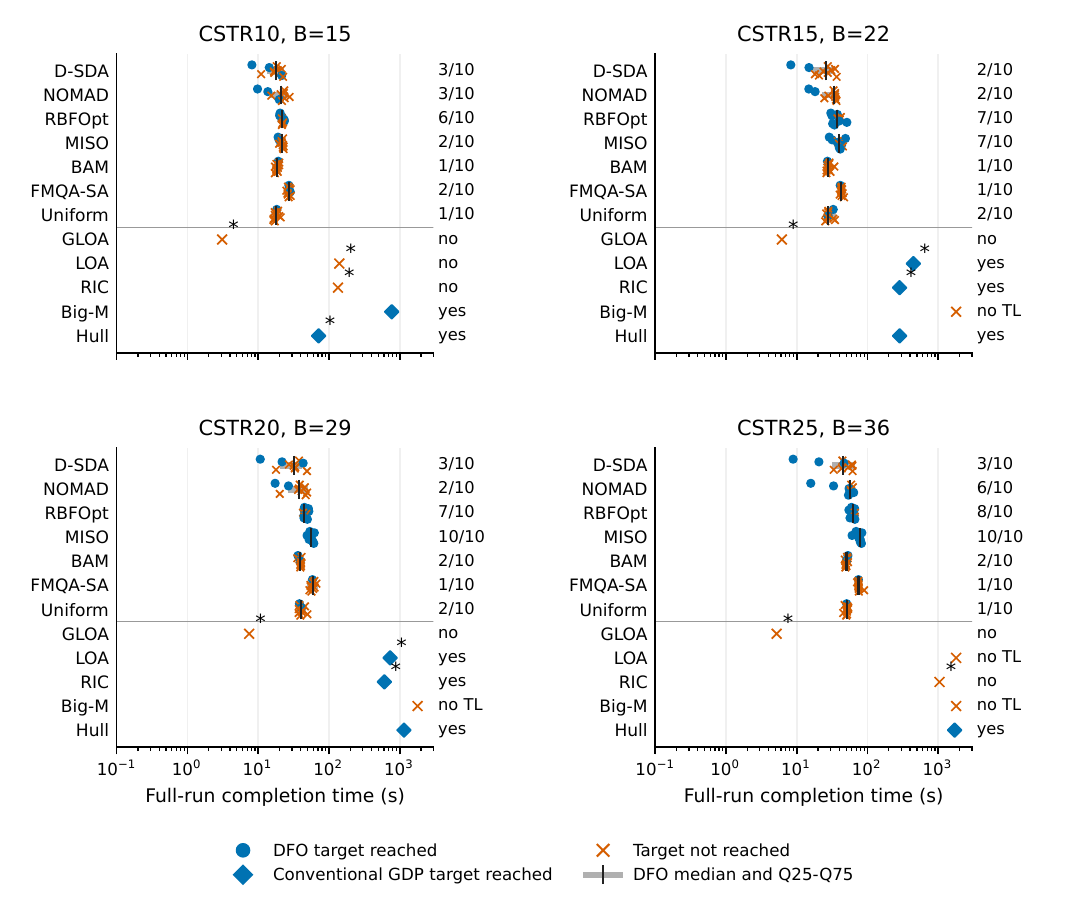}
\caption{Process GDP models: full-run completion times on \texttt{CSTR10} to \texttt{CSTR25}, which Figure~\ifdefined\manuscriptwithsupplement\ref{fig:process-completion}\else3\fi{} of the main paper does not show, of the six DFO methods and uniform sampling (ten runs each, at the evaluation budget \(B\) of each panel title) and of the single run of each conventional GDP method.
GLOA: global logic-based outer approximation; LOA: logic-based outer approximation; RIC: relaxation with integer cuts; Big-\(M\) and Hull: reformulations solved by \texttt{SCIP}.
DFO times add the search time to the measured cost of each evaluated configuration (Equation~\ifdefined\manuscriptwithsupplement\eqref{eqn:reconstructed-cost}\else(10)\fi{} of the main paper); conventional GDP times are measured and include model preparation.
Dots mark DFO runs that reach the reference target of Equation~\ifdefined\manuscriptwithsupplement\eqref{eqn:target}\else(8)\fi{} of the main paper, an objective within \(10^{-5}\max\{1,|f_{\mathrm{ref}}|\}\) of the reference objective \(f_{\mathrm{ref}}\) of Table~\ifdefined\manuscriptwithsupplement\ref{tab:process-cases}\else1\fi{} of the main paper; diamonds mark conventional GDP runs whose feasible objective is at most \(f_{\mathrm{ref}}+10^{-5}\max\{1,|f_{\mathrm{ref}}|\}\); crosses mark the other runs; gray bars span the 25th to 75th percentiles of the ten DFO runs, with a tick at the median.
The right-hand labels give the target hits out of ten of each DFO method and the outcome, yes or no, of each conventional run, with TL for runs that report termination at the 1,800~s time limit.
Asterisks mark conventional runs whose reported lower bound exceeds a checked feasible objective (Section~\ref{sec:bound-audit}).
Uniform: uniform sampling without replacement.}
\label{fig:joint-remaining-completion}
\end{figure}

\begin{figure}[htbp]
\centering
\includegraphics[width=\linewidth]{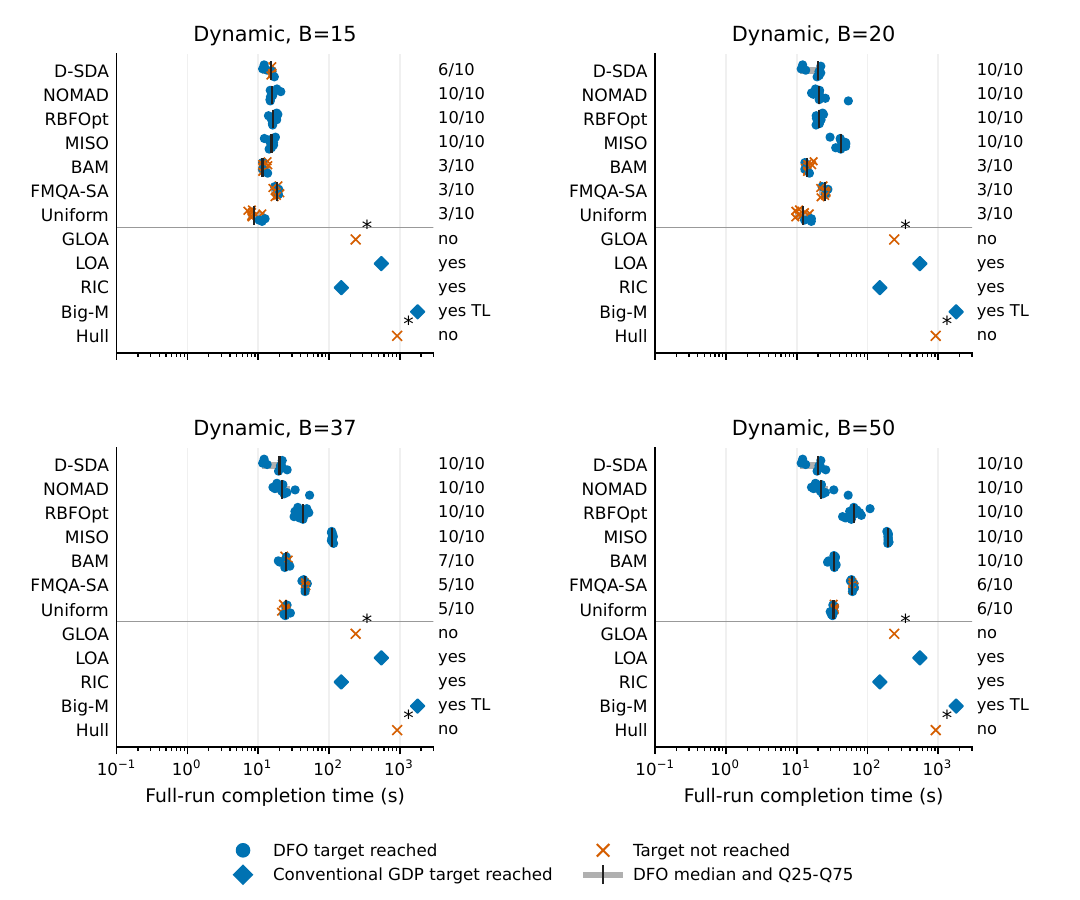}
\caption{Process GDP models: full-run completion times on \texttt{Dynamic} at \(B=15\), 20, 37, and 50, which Figure~\ifdefined\manuscriptwithsupplement\ref{fig:process-completion}\else3\fi{} of the main paper does not show, of the six DFO methods and uniform sampling (ten runs each, at the evaluation budget \(B\) of each panel title) and of the single run of each conventional GDP method.
GLOA: global logic-based outer approximation; LOA: logic-based outer approximation; RIC: relaxation with integer cuts; Big-\(M\) and Hull: reformulations solved by \texttt{SCIP}.
DFO times add the search time to the measured cost of each evaluated configuration (Equation~\ifdefined\manuscriptwithsupplement\eqref{eqn:reconstructed-cost}\else(10)\fi{} of the main paper); conventional GDP times are measured and include model preparation.
Dots mark DFO runs that reach the reference target of Equation~\ifdefined\manuscriptwithsupplement\eqref{eqn:target}\else(8)\fi{} of the main paper, an objective within \(10^{-5}\max\{1,|f_{\mathrm{ref}}|\}\) of the reference objective \(f_{\mathrm{ref}}\) of Table~\ifdefined\manuscriptwithsupplement\ref{tab:process-cases}\else1\fi{} of the main paper; diamonds mark conventional GDP runs whose feasible objective is at most \(f_{\mathrm{ref}}+10^{-5}\max\{1,|f_{\mathrm{ref}}|\}\); crosses mark the other runs; gray bars span the 25th to 75th percentiles of the ten DFO runs, with a tick at the median.
The right-hand labels give the target hits out of ten of each DFO method and the outcome, yes or no, of each conventional run, with TL for runs that report termination at the 1,800~s time limit.
Asterisks mark conventional runs whose reported lower bound exceeds a checked feasible objective (Section~\ref{sec:bound-audit}).
The conventional GDP methods have no evaluation budget, so the four panels show the same conventional runs.
Uniform: uniform sampling without replacement.}
\label{fig:joint-remaining-completion-dynamic}
\end{figure}

Figures~\ref{fig:joint-remaining-completion} and~\ref{fig:joint-remaining-completion-dynamic} compare the eight process conditions not shown in the main joint completion-time figure.
Each DFO condition has ten seeds, whereas each conventional GDP method has one run per problem.
Timing and target criteria follow Figure~\ifdefined\manuscriptwithsupplement\ref{fig:process-completion}\else3\fi{} of the main paper.
Table~\ref{tab:full-gdp} reports all nine models, including final objectives and termination flags.

\begin{table}[htbp]
\centering
\caption{Process GDP models: final objectives and full-run completion times of the conventional GDP methods, with one run per method and model.
The \texttt{GDPopt} methods GLOA, LOA, and RIC use \texttt{SCIP} for the master problems and \texttt{IPOPT} for the NLP subproblems, and the Big-\(M\) and Hull columns are the retained runs that solve these reformulations directly with \texttt{SCIP}.
The final objective is in the unit of the row label, as in Table~\ifdefined\manuscriptwithsupplement\ref{tab:process-cases}\else1\fi{} of the main paper, with no physical unit for \texttt{Column} and \texttt{Dynamic}, and the completion time is the measured full-run time, including model-template preparation.
The success threshold of a conventional GDP run is a final objective of at most \(f_{\mathrm{ref}}+10^{-5}\max\{1,|f_{\mathrm{ref}}|\}\), and feasibility means numerical feasibility under the original-GDP residual test at \(10^{-5}\), which does not establish exact feasibility.
Runs that do not stop at the time limit report \texttt{optimal} termination.}
\label{tab:full-gdp}
\footnotesize
\setlength{\tabcolsep}{3.5pt}
\begin{threeparttable}
\begin{tabular}{@{}lrlr@{}lr@{}lr@{}lr@{}lr@{}l@{}}
\toprule
Model & \(N\) & & \multicolumn{2}{c}{GLOA} & \multicolumn{2}{c}{LOA} & \multicolumn{2}{c}{RIC} & \multicolumn{2}{c}{Big-\(M\)} & \multicolumn{2}{c}{Hull} \\
\midrule
\multirow{12}{*}{\texttt{CSTR}} & \multirow{2}{*}{5} & Final objective (L) & 3.13020 & \({}^{\dagger\ddagger}\) & 3.06201 & \({}^{\dagger}\) & 3.06201 & \({}^{\dagger}\) & \textbf{3.06198} & & 3.06201 & \({}^{\dagger}\) \\
 & & Completion time (s) & \textbf{4.57} & & 24.58 & & 15.38 & & 53.68 & & 15.43 & \\
\addlinespace[2pt]
 & \multirow{2}{*}{10} & Final objective (L) & 3.44796 & \({}^{\dagger\ddagger}\) & 2.96605 & \({}^{\dagger\ddagger}\) & 2.96605 & \({}^{\dagger\ddagger}\) & \textbf{2.88941} & & 2.88951 & \({}^{\dagger}\) \\
 & & Completion time (s) & \textbf{3.10} & & 140.38 & & 134.31 & & 771.36 & & 71.63 & \\
\addlinespace[2pt]
 & \multirow{2}{*}{15} & Final objective (L) & 3.56196 & \({}^{\dagger\ddagger}\) & 2.82939 & \({}^{\dagger}\) & 2.82939 & \({}^{\dagger}\) & -- & \({}^{\mathrm{TL}}\) & \textbf{2.82936} & \\
 & & Completion time (s) & \textbf{6.17} & & 445.43 & & 283.77 & & 1800+ & & 284.63 & \\
\addlinespace[2pt]
 & \multirow{2}{*}{20} & Final objective (L) & 3.61612 & \({}^{\dagger\ddagger}\) & 2.79941 & \({}^{\dagger}\) & 2.79941 & \({}^{\dagger}\) & 2.80413 & \({}^{\mathrm{TL}}\) & \textbf{2.79936} & \\
 & & Completion time (s) & \textbf{7.46} & & 733.53 & & 607.72 & & 1800+ & & 1{,}150.67 & \\
\addlinespace[2pt]
 & \multirow{2}{*}{25} & Final objective (L) & 3.64786 & \({}^{\dagger\ddagger}\) & 2.79125 & \({}^{\mathrm{TL}}\) & 2.79125 & \({}^{\dagger\ddagger}\) & 2.78771 & \({}^{\mathrm{TL}}\) & \textbf{2.78143} & \\
 & & Completion time (s) & \textbf{5.20} & & 1800+ & & 1{,}042.79 & & 1800+ & & 1{,}701.75 & \\
\addlinespace[2pt]
 & \multirow{2}{*}{30} & Final objective (L) & 3.66872 & \({}^{\dagger\ddagger}\) & 2.87314 & \({}^{\mathrm{TL}}\) & 2.77162 & \({}^{\mathrm{TL}}\) & \textbf{2.76957} & \({}^{\mathrm{TL}}\) & \textbf{2.76957} & \({}^{\mathrm{TL}}\) \\
 & & Completion time (s) & \textbf{5.52} & & 1800+ & & 1800+ & & 1800+ & & 1800+ & \\
\midrule
\multirow{2}{*}{\texttt{Column}} & \multirow{2}{*}{--} & Final objective & 20{,}355.1 & \({}^{\dagger\ddagger}\) & \textbf{19{,}449.9} & & \textbf{19{,}449.9} & & -- & \({}^{\mathrm{TL}}\) & 22{,}271.5 & \({}^{\mathrm{TL}\dagger}\) \\
 & & Completion time (s) & \textbf{42.41} & & 199.22 & & 124.03 & & 1800+ & & 1800+ & \\
\midrule
\multirow{2}{*}{\texttt{Small Batch}} & \multirow{2}{*}{--} & Final objective (\$) & \textbf{167{,}428} & & \textbf{167{,}428} & & \textbf{167{,}428} & & \textbf{167{,}428} & & \textbf{167{,}428} & \\
 & & Completion time (s) & 10.32 & & 5.61 & & 13.15 & & 0.51 & & \textbf{0.50} & \\
\midrule
\multirow{2}{*}{\texttt{Dynamic}} & \multirow{2}{*}{--} & Final objective & -197.140 & \({}^{\dagger\ddagger}\) & \textbf{-212.009} & & \textbf{-212.009} & & \textbf{-212.009} & \({}^{\mathrm{TL}}\) & -196.931 & \({}^{\dagger\ddagger}\) \\
 & & Completion time (s) & 238.96 & & 550.49 & & \textbf{149.39} & & 1800+ & & 921.28 & \\
\bottomrule
\end{tabular}
\begin{tablenotes}[flushleft]\footnotesize
\item \(N\): number of potential reactors of \texttt{CSTR}; \(f_{\mathrm{ref}}\): reference objective of Table~\ifdefined\manuscriptwithsupplement\ref{tab:process-cases}\else1\fi{} of the main paper; GLOA: global logic-based outer approximation; LOA: logic-based outer approximation; RIC: relaxation with integer cuts.
1800+: stopped at the 1,800~s time limit; TL: stopped at the time limit, with reported \texttt{maxTimeLimit}; --: no accepted feasible solution in an objective cell and not applicable in the \(N\) column.
\(\dagger\): reported lower bound exceeds the best checked saved numerically feasible objective \(f\) of the model by more than \(10^{-5}\max\{1,|f|\}\) (Section~\ref{sec:bound-audit}); \(\ddagger\): final objective above the success threshold, despite reported \texttt{optimal} termination and a better checked numerically feasible solution.
Bold: lowest final objective of each model, including ties within \(10^{-5}\max\{1,|f_{\min}|\}\) of that lowest objective \(f_{\min}\), and shortest completion time of each model, including ties.
\end{tablenotes}
\end{threeparttable}
\end{table}
\FloatBarrier

\begin{table}[htbp]
\centering
\caption{Process GDP models: NLP subproblems solved per run by the six DFO methods and uniform sampling, with ten runs per method and condition, and by the single run of each \texttt{GDPopt} method, which uses \texttt{SCIP} for the master problems and \texttt{IPOPT} for the NLP subproblems.
The DFO entries are medians over the ten runs of the number of evaluations that called \texttt{IPOPT}, whether or not its result was accepted, without logical and preprocessing rejections, and the \texttt{GDPopt} entries are the \texttt{IPOPT} calls confirmed in the execution log of the run.
The counts are numbers of subproblem solves, which differ in cost, and do not include solver iterations; the internal NLP counts of the direct Big-\(M\) and Hull solves are unavailable and are not treated as zero.
Every \texttt{GDPopt} NLP call fixes a different configuration, out of \(N(N+1)/2\) logically admissible configurations for \texttt{CSTR} with \(N\) potential reactors and 37 for \texttt{Dynamic} (Section~\ref{sec:supp-spread-overhead}).
The \texttt{GDPopt} methods have no evaluation budget, so the single \texttt{Dynamic} run of each appears in every \texttt{Dynamic} row.}
\label{tab:joint-nlp-workload}
\footnotesize
\setlength{\tabcolsep}{3.5pt}
\begin{threeparttable}
\begin{tabular}{@{}lrrrrrrrrrrrr@{}}\toprule
 & & & \multicolumn{7}{c}{DFO methods and uniform sampling (median of ten)} & \multicolumn{3}{c}{\texttt{GDPopt} (one run)}\\\cmidrule(lr){4-10}\cmidrule(l){11-13}
Model & \(N\) & \(B\) & \texttt{D-SDA} & \texttt{NOMAD} & \texttt{RBFOpt} & \texttt{MISO} & \texttt{BAM} & \texttt{FMQA-SA} & Uniform & GLOA & LOA & RIC\\\midrule
\multirow{6}{*}{\texttt{CSTR}} & 5 & 8 & 6 & 6 & 6 & 6 & 6 & 5.5 & 4.5 & 4 & 15 & 15 \\
 & 10 & 15 & 10.5 & 10.5 & 12 & 11.5 & 8.5 & 9 & 7.5 & 2 & 55 & 55 \\
 & 15 & 22 & 14 & 15 & 16.5 & 16.5 & 11 & 15 & 11 & 4 & 120 & 120 \\
 & 20 & 29 & 18 & 20 & 20.5 & 19.5 & 17 & 17.5 & 15.5 & 4 & 210 & 210 \\
 & 25 & 36 & 22 & 25 & 25.5 & 25 & 18.5 & 23 & 21 & 2 & 315 & 325 \\
 & 30 & 43 & 23 & 25 & 31.5 & 29.5 & 21.5 & 30.5 & 21.5 & 2 & 243 & 455 \\
\midrule
\texttt{Column} & -- & 13 & 11 & 10.5 & 12 & 11 & 7.5 & 7 & 7 & 18 & 38 & 38 \\
\midrule
\texttt{Small Batch} & -- & 10 & 5 & 6 & 6 & 6 & 5 & 5 & 4 & 12 & 5 & 12 \\
\midrule
\multirow{6}{*}{\texttt{Dynamic}} & -- & 12 & 10 & 10 & 10 & 9 & 8 & 6.5 & 6 & 20 & 37 & 37 \\
 &  & 15 & 12.5 & 12 & 13 & 12 & 9 & 8 & 6.5 & 20 & 37 & 37 \\
 &  & 20 & 16.5 & 14.5 & 16 & 14 & 11 & 11 & 9.5 & 20 & 37 & 37 \\
 &  & 30 & 16.5 & 14.5 & 21 & 21 & 17 & 14.5 & 15 & 20 & 37 & 37 \\
 &  & 37 & 16.5 & 14.5 & 25 & 25 & 19 & 19 & 19 & 20 & 37 & 37 \\
 &  & 50 & 16.5 & 14.5 & 31 & 29 & 26 & 24.5 & 26 & 20 & 37 & 37 \\
\bottomrule
\end{tabular}
\begin{tablenotes}[flushleft]\footnotesize
\item \(N\): number of potential reactors of \texttt{CSTR}; \(B\): evaluation budget of the DFO methods; --: not applicable; Uniform: uniform sampling without replacement.
GLOA: global logic-based outer approximation; LOA: logic-based outer approximation; RIC: relaxation with integer cuts.
\end{tablenotes}
\end{threeparttable}
\end{table}
\FloatBarrier

\FloatBarrier

\ifdefined\manuscriptwithsupplement
\FloatBarrier 
\endgroup
\let\finishsupplement\relax
\else
\bibliographystyle{elsarticle-harv}
\bibliography{reference}
\def\finishsupplement{\end{document}}
\fi
\finishsupplement
\fi
\end{document}